\def\paperversion{2}
\ifnum\paperversion=1
\documentclass[opre, blindrev]{informs4}
\fi
\ifnum\paperversion=2
\documentclass[opre]{informs4}
\fi
\OneAndAHalfSpacedXI

\usepackage{silence}

\usepackage{endnotes}
\let\footnote=\endnote

\usepackage[sort]{natbib}
\bibpunct[, ]{(}{)}{,}{a}{}{,}%
 \def\bibfont{\small}%

\TheoremsNumberedThrough     %
\ECRepeatTheorems

\EquationsNumberedThrough    %
\usepackage{anyfontsize}
\usepackage{multirow}
\usepackage{booktabs}
\usepackage{array}
\usepackage{mathtools}
\usepackage{lipsum}
\usepackage{enumitem}
\usepackage{algorithm}
\usepackage{algorithmic}

\usepackage{comment}
\let\hat\widehat
\let\tilde\widetilde
\usepackage{tikz}
\usetikzlibrary{positioning, arrows}
\usepackage[small, margin=1cm]{caption}
\usepackage{color}
\makeatletter
\newcommand{\mylabel}[2]{#2\def\@currentlabel{#2}\label{#1}}
\makeatother

\usepackage{subfigure}
\usepackage[export]{adjustbox}
\usepackage{xcolor}

\usepackage{verbatim}

\usepackage{amsmath,amsfonts,bm}
\usepackage{amssymb}

\usepackage[utf8]{inputenc} 
\usepackage[T1]{fontenc}    
\usepackage{hyperref}       
\hypersetup{
    colorlinks=true,
    citecolor=blue,
    filecolor=blue,
    linkcolor=blue,
    urlcolor=blue
}
\usepackage{url}            
\usepackage{booktabs}       
\usepackage{amsfonts}       
\usepackage{nicefrac}       
\usepackage{microtype}      
\usepackage{xcolor}         

\usepackage{macros}
\usepackage{wrapfig}

\usepackage{enumitem}
\setlist[enumerate,1]{label=\normalfont{(\Roman*)},leftmargin=*}

\newcommand{\cv}[1]{\mathbf{#1}}

\newcommand{\LSGD}{L\text{-}\mathrm{SGD}}
\newcommand{\VMLMC}{\mathrm{V}\text{-}\mathrm{MLMC}}
\newcommand{\RTMLMC}{\mathrm{RT}\text{-}\mathrm{MLMC}}
\newcommand{\RRMLMC}{\mathrm{RR}\text{-}\mathrm{MLMC}}
\newcommand{\RUMLMC}{\mathrm{RU}\text{-}\mathrm{MLMC}}

\definecolor{longhorn}{rgb}{0.8, 0.33, 0.0}

\allowdisplaybreaks[4]
\newcommand*{\QED}{%
\leavevmode\unskip\penalty9999 \hbox{}\nobreak\hfill
    \quad\hbox{$\square$}%
}

\newtheorem{innercustomgeneric}{\customgenericname}
\providecommand{\customgenericname}{}
\newcommand{\newcustomtheorem}[2]{%
  \newenvironment{#1}[1]
  {%
   \renewcommand\customgenericname{#2}%
   \renewcommand\theinnercustomgeneric{##1}%
   \innercustomgeneric
  }
  {\endinnercustomgeneric}
}
\newcustomtheorem{customthm}{Assumption}

\begin{document}

\RUNAUTHOR{Hu et al.}

\RUNTITLE{Stochastic Opt with Biased Oracles}

\TITLE{Multi-level 
Monte-Carlo Gradient Methods for Stochastic Optimization with 
Biased Oracles
}

\ARTICLEAUTHORS{%
\AUTHOR{Yifan Hu}
\AFF{College of Management of Technology, EPFL, Switzerland, \EMAIL{yifan.hu@epfl.ch}} 
\AUTHOR{Jie Wang}
\AFF{ISyE, Georgia Institute of Technology, \EMAIL{jwang3163@gatech.edu}}
\AUTHOR{Xin Chen}
\AFF{ISyE, Georgia Institute of Technology, \EMAIL{xin.chen@isye.gatech.edu}}
\AUTHOR{Niao He}
\AFF{Computer Science Department, ETH Z\"urich, Switzerland, \EMAIL{niao.he@inf.ethz.ch}} 
}

\ABSTRACT{

We consider stochastic optimization when one only has access to biased stochastic oracles of the objective and the gradient, and obtaining stochastic gradients with low biases comes at high costs. This setting captures various optimization paradigms,  such as conditional stochastic optimization, distributionally robust optimization, shortfall risk optimization, and machine learning paradigms, such as contrastive learning. We examine a family of multi-level Monte Carlo (MLMC) gradient methods that exploit a delicate tradeoff among bias, variance, and oracle cost. We systematically study their total sample and computational complexities for strongly convex, convex, and nonconvex objectives and demonstrate their superiority over the widely used biased stochastic gradient method. When combined with the variance reduction techniques like SPIDER, these MLMC gradient methods can further reduce the complexity in the nonconvex regime. Our results imply that a series of stochastic optimization problems with biased oracles, previously considered to be more challenging, is fundamentally no harder than the classical stochastic optimization with unbiased oracles. We also delineate the boundary conditions under which these problems become more difficult. Moreover, MLMC gradient methods significantly improve the best-known complexities in the literature for conditional stochastic optimization and shortfall risk optimization. 
Our extensive numerical experiments on distributionally robust optimization, pricing and staffing scheduling problems, and contrastive learning demonstrate the superior performance of MLMC gradient methods.
\footnote{A preliminary version of this manuscript has appeared in a conference proceeding. Please refer to \cite{hu2021bias} Y.~Hu, X.~Chen, and N.~He. On the bias-variance-cost tradeoff of stochastic optimization. {\it Advances in Neural Information Processing Systems, 34}.}

}
\maketitle

\section{Introduction}
\label{section:introduction}
\noindent
SGD and its numerous variants are commonly used optimization approaches to address large-scale, data-driven applications in machine learning and optimization.
While vanilla SGD depends crucially on unbiased gradient oracles, constructing such estimators can be prohibitively expensive or even infeasible for many emerging machine learning and optimization problems. 
Examples include 
distributionally robust optimization~\citep{levy2020large,ghosh2020unbiased},
conditional stochastic optimization~\citep{hu2020biased},
meta-learning~\citep{rajeswaran2019meta, fallah2019on}, contextual optimization~\citep{bertsimas2020from,bertsimas2021data}, variantional inference~\citep{blei2017variational,fujisawa2021multilevel}, and etc.

As an alternative, one often resorts to some naturally biased gradient estimators. 
However, existing literature~\citep{hu2020analysis, ajalloeian2020analysis,hu2016bandit, chen2018stochastic,karimi2019non,demidovich2024guide,liu2023adaptive} usually focus on the iteration complexity given biased stochastic oracles while ignoring the fact that constructing stochastic oracles with lower bias and variances incurs higher costs, such as more samples or computation resources. 

In this paper, we build a unified framework to study the fundamental tradeoff among bias, variance, and cost for stochastic optimization with biased gradient oracles. We investigate 
efficient gradient-based methods to solve optimization problems of the general form
\begin{equation}
\label{eq:problem}
\min_{x\in\RR^d} F(x),
\vspace{-0.5em}
\end{equation}
where we assume that one does not have access to an unbiased gradient estimator of $F(x)$ via simple Monte Carlo sampling.
Instead, assume there exists a sequence of approximations of $F(x)$, denoted as $\{F^l(x)\}_{l=0}^\infty$ such that the error $|F^l(x)-F(x)|$ or $\|\nabla F^l(x) - \nabla F(x)\|$ vanishes as $l\to\infty$. 
Unbiased gradient estimators of $F^l$ are attainable through stochastic oracles $\mathcal{SO}^l$, but the cost of querying $\mathcal{SO}^l$ increases with $l$, as sampling from more accurate approximations generally entails higher costs.
For simplicity of discussion, we assume the continuous differentiability of $F$ and $F^l$.

\vspace{-0.5em}
\subsection{Motivating Examples}\label{sub:motivating}
\vspace{-0.5em}

\subsubsection*{Conditional Stochastic Optimization (CSO).}
The objective of CSO~\citep{hu2020sample} involves the expectation of compositions of nonlinear function with another conditional expectation, i.e., $\min_x \EE_\xi [f_\xi (\EE_{\eta\mid \xi} g_\eta(x;\xi))]$.
It is challenging to obtain the unbiased gradient estimator due to the composition of nonlinear transformation with conditional expectation.
To tackle this difficulty, one can construct a sequence of approximation functions for the CSO objective, where the level $l$ implies that one uses $2^l$ samples to estimate the conditional expectation.
It has been shown in \citep[Theorem~4.1]{hu2020sample} that the bias of the approximation sequence decreases.
As we will discuss in Section~\ref{section:CSO}, two widely studied applications, distributionally robust optimization~(DRO)~\citep{wang2021sinkhorn,levy2020large} and contrastive learning~(CL)~\citep{chen2020simple}, are special cases of CSO. For more applications of CSO in reinforcement learning, variational inference, and others, interested readers may refer to \citep{hu2020biased}.

\subsubsection*{Pricing and Staffing in Stochastic Systems.}
The pricing and staffing task in a stochastic system seeks to optimize the long-run expected profit~\citep{lee2014optimal}, whose objective involves the expectation of steady-state distributions, which can be simulated by taking the time limit of a queuing process.
Since the unbiased gradient estimator is difficult to obtain, we resort to an approximate gradient method, where we build the $l$-th level approximated objective using the first $2^l$ samples of the trajectory from the queuing process.
As we will discuss in Section~\ref{Sec:pricing:staffing}, the bias of the approximate gradient decreases at a linear rate.

\subsubsection*{Utility-Based Shortfall Risk~(SR) Optimization.} SR optimization~\citep{bertsimas2004shortfall} appears widely in mathematical finance and power system.
The gradient of the objective involves (i) a fraction of two expectations with respect to continuous distributions and (ii) a descriptive statistic that requires a large number of samples to obtain its accurate estimation, which makes it challenging to obtain an unbiased gradient estimator.
Fortunately, it is feasible to construct an approximate gradient estimate with bias $\cO(2^{-l})$ and $\cO(1)$ variance using $\cO(2^l)$ samples, which will be discussed in Section~\ref{Sec:UBSR}.

\subsection{Multilevel Monte Carlo (MLMC) Gradient Estimation}
A classical approach to solve \eqref{eq:problem} is to perform SGD on the approximation function $F^L(x)$ with a level $L$ such that the bias is small and within a target accuracy $\epsilon$. We call such a framework 
$L$-SGD. It encapsulates several existing algorithms under different contexts, such as BSGD~\citep{hu2020biased} in conditional stochastic optimization, biased gradient method~\citep{levy2020large} in DRO,  a variant of LEON for contextual SO~\citep{diao2020distribution}, and various double-loop methods for meta-learning and bilevel optimization~\citep{franceschi2018bilevel,hu2023contextual}. Despite its simplicity, $L$-SGD generally fails to achieve the smallest total cost since it requires expensive stochastic oracles to ensure a small bias. 

On the other hand, the \emph{multilevel Monte Carlo} (MLMC) sampling technique, initially designed for stochastic simulation~\citep{giles2008multilevel}, is amenable to obtain better gradient estimators by exploiting the bias-cost tradeoff. 
Consider the following representation of the objective and gradient of $F^L$ for a sufficiently large $L$:
$$
 F^L(x) =  \sum_{l=0}^L  [ F^l(x) - F^{l-1}(x)]; \quad \nabla F^L(x) =  \sum_{l=0}^L  [\nabla F^l(x) -\nabla F^{l-1}(x)],
$$
where we denote $F^{-1}(x):=0$ and $\nabla F^{-1}(x):=0$ by convention. Let $H^l$ be an estimator of $\nabla F^l -\nabla F^{l-1}$. Such telescoping sum representation suggests that we can construct estimators of $\nabla F^L$ by summing up $H^l$ from $l=0$ to $l=L$. As a result, we can use different batches for  $H^l$ with different $l$ or assign different weights in a randomized construction.
The idea is to query fewer oracles with higher costs and smaller biases (oracles with large $l$) and query more oracles with lower costs and larger biases (oracles with small $l$); thus, the cost can be effectively reduced. Note that the construction of the estimator $H^l$ requires that the estimators of $\nabla F^l$ and $\nabla F^{l-1}$ to be correlated so that one can well control the variance of the MLMC estimator, which we will discuss for various applications in Section \ref{section:applications}. In this work, we consider four MLMC gradient constructions. The vanilla MLMC (denoted as V-MLMC) utilizes sample averages and implements varying mini-batch sizes for the estimator $H^l$  with different levels $l$. 
 There exist several variations of MLMC with randomization, including the unbiased MLMC estimator with importance sampling (denoted RU-MLMC), 
 the MLMC estimator with randomized truncation and importance sampling (denoted as RT-MLMC)~\citep{blanchet2015unbiased}, and the 
 Russian roulette estimator (denoted as RR-MLMC) that uses randomized telescope sum~\citep{kahn1955use}. 

Although the technique of MLMC could also be used to construct function value estimators or solution estimators~\citep{blanchet2019unbiased,frikha2016multi}, these two methods have their disadvantages in optimization, respectively. An MLMC function value estimator, unlike sample average approximation, is usually nonconvex even if the original objective possesses (strong) convexity since it uses estimators of $F^l(x) - F^{l-1}(x)$ as  building blocks. Thus, it would be difficult to find the global optimal solution. On the other hand, using MLMC to build solution estimators requires the empirical objective to be strongly or strictly convex so that the optimal empirical solution is unique. Otherwise, there is no guarantee of the quality of the MLMC solution estimator. In addition, it requires solving multiple empirical problems. We discuss further details in Section \ref{Sec:statistical:inference}. Note that MLMC gradient methods overcome these disadvantages.

\subsection{Our Contributions }

In this paper, we systematically compare these MLMC techniques when combined with stochastic gradient descent and variance reduction methods. Our primary focus is on the analysis of their total sampling and computation cost for achieving an $\epsilon$-optimal solution for (strongly) convex problems and an $\epsilon$-stationary point for nonconvex problems, which is largely missing in the literature. We also highlight the importance of using MLMC for gradient estimation rather than function value or solution estimation.

    \begin{table}[t]

    {
    \footnotesize
    \caption{Summary of (expected) total costs for finding an $\eps$-optimal solution for (strongly) convex $F$ or an $\eps$-stationary point for nonconvex $F$. Here $a$ is the decrease rate in the bias, i.e., $|F^l(x)- F(x)|=\cO(2^{-al})$ for (strongly) convex case, or the gradient approximation error $\|\nabla F^l(x)- \nabla F(x)\|_2^2=\cO(2^{-al})$ for nonconvex case. $b$ is the decrease rate of the variance of $H^l$, the estimator of $\nabla F^l(x) - \nabla F^{l-1}(x)$, i.e., $\mathrm{Var}(H^l)=\cO(2^{-bl})$.    $c$ is the increase rate of oracle cost, i.e., the cost to generate $H^l$ is $\cO(2^{cl})$. See more details in Assumptions \ref{assumption:general} and \ref{assumption:gradient_bias}. $\tilde \cO(\cdot)$ represents the order hiding logarithmic factors. N.A. stands for not applicable, VR stands for variance reduction. Let $n_1 = \max\{1, c/a\}$, $n_2 = \max\{1+(c-b)/a, c/a\}$.
    }
    \renewcommand\arraystretch{1.5}
    \begin{center}
     \makebox[\textwidth]{
    \begin{tabular}{c|ccc|c|c}
        \toprule[1.5pt]
        {}
        & \multirow{2}{*}{\textbf{V-MLMC}}
        & \multirow{2}{*}{\textbf{RT-MLMC}}
        & \textbf{RU-MLMC}
        & \textbf{$L$-SGD}
        & \multirow{2}{*}{Conditions of $F$, $F^l$}
        \\
        
        {}
        & 
        & 
        & \textbf{RR-MLMC}
        & \textit{(existing biased methods)}
        & 
        \\
        \hline
        
        {Unbiasedness}
        & Biased
        & Biased
        & Unbiased
        & Biased
        & 
        \\
        \cline{1-5}
        
        {Requires Mini-batch}
        & Yes
        & No
        & No
        & No
        &
        \\
        \cline{1-6}

        \multirow{3}{*}{Total Cost  if $c<b$}
        & $\tilde \cO(\eps^{-n_1})$
        & $\cO(\eps^{-1})$
        & $\cO(\eps^{-1})$
        & $\cO(\eps^{-1-c/a})$
        & Strongly Convex
        \\
        \cline{2-6}
        
        {}
        & $\cO(\eps^{-1-n_1})*$
        & $\cO(\eps^{-2})$
        & $\cO(\eps^{-2})$
        & $\cO(\eps^{-2-c/a})$
        & Convex
        \\
        \cline{2-6}

        {}
        & $\cO(\eps^{-2-2n_1})$
        & $\cO(\eps^{-4})$
        & $\cO(\eps^{-4})$
        & $\cO(\eps^{-4-2c/a})$
        & Nonconvex
        \\
        \hline

        \multirow{3}{*}{Total Cost if $c\geq b$}
        & $\tilde \cO(\eps^{-n_2})$
        & $\tilde \cO(\eps^{-1-(c-b)/a})$
        & N.A.
        & $\cO(\eps^{-1-c/a})$
        & Strongly Convex
        \\
        \cline{2-6}
        
        {}
        & $\tilde \cO(\eps^{-1-n_2})*$
        & $\tilde \cO(\eps^{-2-(c-b)/a})$
        & N.A.
        & $\cO(\eps^{-2-c/a})$
        & Convex
        \\
        \cline{2-6}
        
        {}
        & $\tilde  \cO(\eps^{-2-2n_2})$
        & $\tilde  \cO(\eps^{-4-2(c-b)/a})$
        & N.A.
        & $\cO(\eps^{-4-2c/a})$
        & Nonconvex\\
       \hline
       
        Total Cost of VR when $c<b$
        & $\cO(\eps^{-1-2n_1})$
        & $ \cO(\eps^{-3})$
        & $ \cO(\eps^{-3})$
        & $\cO(\eps^{-3-2c/a})$
        &  \multirow{2}{*}{Nonconvex }
        \\
        \cline{1-5}
        
        {Total Cost of VR when $c\geq b$}
        & $\cO(\eps^{-1-2n_2}) $
        & $\tilde \cO(\eps^{-3-2(c-b)/a})$
        & N.A.
        & $\cO(\eps^{-3-2c/a})$
        & 
        \\
         \hline
         \multicolumn{6}{c}{$*$: if the approximation function $F^l$ is only Lipschitz continuous, i.e., Assumption \ref{assumption:objective}\ref{assumption:objective:II} holds, then $n_1 = n_2 =1+c/a$.}
         \\
         \multicolumn{6}{c}{Total cost of VR denotes the total costs of variance-reduced methods using the corresponding gradient estimators.}
         \\
         \bottomrule[1.5pt]
    \end{tabular} 
    }
    \label{tab:summary}
    \end{center}
    }
    \vspace{-1em}
    \end{table} 

\subsubsection*{Non-asymptotic analysis of MLMC gradient methods.} We analyze the non-asymptotic performances and computation complexities of $L$-SGD and four MLMC gradient methods, i.e., V-MLMC, RT-MLMC, RU-MLMC, and RR-MLMC, under the general stochastic optimization with biased oracle framework in various regime, e.g., different convexity assumptions. We particularly consider different combinations of setup for the estimator $H^l$, i.e., the estimator of $\nabla F^l(x) -\nabla F^{l-1}(x)$. We discuss how the decrease rate of bias (with respect to $a$), the decrease rate of variance (with respect to $b$), and the increase rate of costs (with respect to $c$) as stated in Assumption \ref{assumption:general} affect the total costs. Our main results are summarized in Table \ref{tab:summary}. In contrast, previous results either focus on specific applications~\citep{blanchet2017unbiased, levy2020large} or only asymptotic behaviors~\citep{beatson2019efficient} in restricted regimes. 

\subsubsection*{Variance reduced MLMC gradient methods.} In the nonconvex smooth case, we combine MLMC gradient estimators with variance reduction methods~\citep{wang2018spiderboost} to further reduce the (expected) total cost. When $b>c$, VR RT-MLMC, VR RU-MLMC, and VR RR-MLMC can achieve $\cO(\eps^{-3})$ expected total cost that matches the lower bounds for classical stochastic nonconvex optimization under the average smooth condition~\citep{arjevani2019lower}. When $b\leq c$, VR RU-MLMC and VR RR-MLMC are no longer applicable, while VR RT-MLMC is still better than VR $L$-SGD. As for VR V-MLMC, similar to the V-MLMC in the SGD framework, it needs large mini-batches and can only achieve reduced total cost when $b\geq \max \{a,c\}$. 
Our result provides a general variance reduction algorithmic framework for stochastic optimization with biased oracles. Note that MLMC gradient estimators can combine with other optimization techniques, such as acceleration. 

\subsubsection*{Theoretical comparisons and new insights.} Our comparative study of different MLMC techniques yields several interesting findings: (1) When $b>c$, i.e., the variance decays faster than the increase of the cost, V-MLMC (only when $a\geq c$) and all three randomized MLMC gradient methods nearly match the iteration complexity of classical unbiased SGD, implying that the problem is no harder than classical stochastic optimization, regardless of the bias. (2) When $b\leq c$,  the unbiased MLMC constructions, RU-MLMC and RR-MLMC, are no longer applicable since either the expected per-iteration cost or the variance of the gradient estimator goes to infinity. In this regime, RT-MLMC and V-MLMC are still applicable, yet RT-MLMC does not require any mini-batch. (3) In all regimes with $0< a<\infty$ (i.e., bias exists), using the naive biased gradient method $L$-SGD is strictly sub-optimal in terms of the total cost, demonstrating the importance of the exploitation of the bias-variance-cost tradeoff in stochastic optimization with biased oracles.  {\color{black}(4) Due to potential rounding issues for finding optimal mini-batches, V-MLMC can only achieve a reduced cost comparing to $L$-SGD when the following conditions holds: using large mini-batches,  $F^l$ is smooth, and $a\geq \max\{b,c\}$, implying limited applicability comparing to RT-MLMC.}

\subsubsection*{MLMC for practical applications.}
We demonstrate that MLMC gradient methods and their variance reduction counterparts significantly improve the sample and computational complexity over existing approaches for conditional stochastic optimization, including DRO and constrastive learning, pricing and staffing problem, and shortfall risk optimization in Section~\ref{section:applications}.
We further validate the efficiency of MLMC gradient methods over the $L$-SGD method using a series of experiments with real and synthetic datasets in Section~\ref{section:numerical_results}. For pricing and staffing problems, the number of samples needed for MLMC is about 10 times smaller than that of $L$-SGD. For contrastive learning, on the CIFAR-100 dataset, the MLMC method takes about 3 hours while $L$-SGD takes about 24 hours to achieve the same test accuracy.

\subsection{Related Literature}

\subsubsection*{Stochastic Gradient Descent and Variance Reduction Methods.}
SGD and its variants form one of the most important optimization methods in machine learning and deep learning~\citep{nemirovski2009robust,bottou2018optimization}. The iteration complexity of SGD for achieving an $\eps$-optimal solution is $\cO(\eps^{-1})$ for strongly convex objectives and $\cO(\eps^{-2})$ for convex objectives. \citet{agarwal2009information} demonstrated matching lower bounds, implying that SGD is optimal in the (strongly) convex regime. In the nonconvex smooth setting, the iteration complexity of SGD for achieving an $\eps$-stationary point is $\cO(\eps^{-4})$~\citep{ghadimi2013stochastic}. Recently, several variance reduction methods~\citep{wang2018spiderboost,nguyen2017sarah,STORM, li2021page} improve the iteration complexity to $\cO(\eps^{-3})$ for nonconvex stochastic optimization under an additional average smooth assumption~\citep{arjevani2019lower}. These bounds are tight as 
 \citet{arjevani2019lower} demonstrated matching lower bounds on the iteration complexity for nonconvex stochastic optimization with or without the average smooth condition. 

 The major distinctions between existing studies and our work lie in two aspects. 
First, they assume access to unbiased stochastic gradient estimators, while we consider only biased ones.
Second, they do not consider the expenses of querying stochastic oracles. In contrast, our framework accounts for the costs related to querying low-bias oracles and focuses on optimizing the overall costs for seeking an $\eps$-optimal solution or $\eps$-stationary point.

\subsubsection*{Biased Gradient Methods.}
Two primary approaches have emerged for constructing (possibly) biased gradient estimators.
The first focuses on constructing estimators with small biases at each iteration, which aligns more with our work.
However, existing references usually ignore the cost of querying the biased oracle but focus on iteration complexity only~\citep{hu2020analysis, ajalloeian2020analysis,hu2016bandit, chen2018stochastic,karimi2019non}.
Various biased gradient methods can be categorized under the proposed $L$-SGD framework for specific applications. 
Notable examples include the BSGD mentioned above for CSO~\citep{hu2020biased},
biased gradient method for $\phi$-divergence DRO~\citep{levy2020large},
and learning enabled optimization with non-parametric estimators~(LEON) for contextual SO~\citep{diao2020distribution},
multi-step Model-Agnostic Meta-Learning~(MAML)~\citep{ji2020multi},
and double-loop algorithms for min-max/bilevel optimization~\citep{jin2020local,ghadimi2018approximation}.
In this paper, we first elucidate the total costs of $L$-SGD and propose to use more efficient MLMC gradient methods.

Another approach leverages the special problem structure to design algorithms that adaptively reduce the gradient bias, 
relying largely on a case-by-case analysis.
For instance, \citet{chen2021tighter} proposed an alternating stochastic gradient method for stochastic nested/min-max/ bilevel optimization, and \citet{yang2021faster} proposed a smoothed stochastic Gradient Descent Ascent~(GDA) for stochastic min-max optimization. \citet{demidovich2024guide} and \citet{ajalloeian2020analysis} summarized biased gradient methods from zeroth-order optimization and gradient compression problems.
\subsubsection*{MLMC Techniques. }
MLMC techniques have received recent attention in various optimization contexts.
Notably, \citet{blanchet2015unbiased} combined the MLMC idea with sample average approximation to solve stochastic optimization. 
\citet{dereich2019general} integrated V-MLMC with Robbins-Monro and Polyak-Ruppert stochastic approximation schemes and provided their convergence guarantees.
RT-MLMC was applied in stochastic compositional optimization with strongly convex objectives~\citep{blanchet2017unbiased} and $\phi$-divergence DRO~\citep{levy2020large, ghosh2020unbiased}.
\citet{beatson2019efficient}  studied the asymptotic behaviors of SGD with RT-MLMC or RR-MLMC estimators. \citet{asi2021stochastic} and  \citet{carmon2022recapp} use RT-MLMC for efficient computation of proximal operator and finite-sum and max-structured minimization, respectively. \citet{dorfman2022adapting} use MLMC to reduce dependence on mixing time for stochastic optimization with Markovian data. \citet{hu2023contextual} use MLMC for a contextual bilevel optimization problem. \citet{alacaoglurevisiting} use MLMC for min-max optimization. The purpose of our paper is to provide a systematic study and comparison of MLMC gradient methods in a general stochastic optimization with biased oracles framework and demonstrate how to use them in a more principled way for different applications that can be modeled under the framework. In addition, we are the first to investigate the variance-reduced MLMC gradient methods.

\vspace{-0.5em}
\subsection{Notations}
A function $f:\mathbb{R}^d\rightarrow\mathbb{R}$ is $L$-Lipschitz continuous if $|f(x)-f(y)|\leq L\|x-y\|_2$ holds for any $x,y\in\RR^d$. 
A function $f$ is $S$-smooth if it is continuously differentiable and $\|\nabla f(x)-\nabla f(y)\|_2\leq S\|x-y\|_2$ for any $x,y$.
A function $f$ is $\mu$-strongly convex if
$
    f(x)-f(y)-\nabla f(y)^\top(x-y)\geq \frac{\mu}{2}\|x-y\|_2^2,
$ for any $x,y$.
Let  $x^*\in \argmin_x F(x)$. 
A point $x$ is an $\eps$-optimal solution if it holds that 
$F(x) - F(x^*)\le \eps$,
and $x$ is an $\eps$-stationary point of $F$ if 
$\|\nabla F(x)\|_2^2\leq \eps^2$.
Polyak-Łojasiewicz (PL) condition~\citep{karimi2016linear} is a generalization of strong convexity such that $\|\nabla F(x)\|_2^2\geq 2\mu(F(x)-F(x^*))$ for $\mu>0$. 
Throughout the paper,  we assume that the desired accuracy $\eps>0$ is small enough so that $\eps^{-1}\gg\log(\eps^{-1})\geq 1$. We use $\cO(\cdot)$ to denote the order in terms of $\eps$ and $\tilde\cO(\cdot)$ to denote the order hiding the logarithmic dependency on $\eps$. 
For $N\in\mathbb{N}_+$, denote $[N]$ as the set $\{1,2,\ldots,N\}$.

\vspace{-0.5em}
\subsection{Organizations}
In Section \ref{section:MLMC_estimators}, we discuss the construction of $L$-SGD and four MLMC gradient estimators, together with their variance reduction counterparts.
In Section \ref{section:total_cost}, we discuss the sample complexity of these methods.
In Section \ref{section:applications}, we demonstrate three applications of MLMC estimators and further discuss the advantages of MLMC gradient methods over MLMC function value estimator and solution estimator in Section~\ref{Sec:statistical:inference} . Section~\ref{section:numerical_results} demonstrates numerical results.

\section{MLMC Gradient Methods}
\label{section:MLMC_estimators}
In this section, we first formally restate the setting for the biased oracles we use throughout the paper. Next, we demonstrate constructions of different biased gradient estimators via either $L$-SGD or MLMC methods. Lastly, we discuss the bias, variance, and cost properties of these gradient estimators.
\vspace{-0.5em}
\subsection{Biased Oracle Setting}\label{Sec:MLMC:assumption}
We first define the stochastic oracle. In particular, we assume access to the following black box stochastic oracle, which independently generates desired stochastic first-order information upon each query.
\begin{assumption}
\label{assumption:general}
There exist constants $a,b,c, M_a, M_b, M_c,\sigma > 0$ such that for any $x\in\RR^d$ and $l\in\NN$, the following conditions hold. 
\begin{enumerate}
    \item\label{assumption:general:I} The function value approximation error is bounded: $|F^l(x)-F(x)|\leq B_l:=M_a 2^{-al}.$
    \item\label{assumption:general:II} There exits a stochastic oracle, $\mathcal{SO}^l$, that for given $x$ returns stochastic estimators $h^l(x,\zeta^l)$ and $H^l(x,\zeta^l)$ such that 
    \begin{align}
        &\EE  h^l(x,\zeta^l) = \nabla F^l(x), &&\mathrm{Var}( h^l(x,\zeta^l)) \leq \sigma^2,\\
        &\EE H^l(x,\zeta^l)=\nabla F^l(x)-\nabla F^{l-1}(x), &&\mathrm{Var}(H^l(x,\zeta^l))\leq V_l:=M_b 2^{-bl}.
    \end{align}
    \item\label{assumption:general:III} The cost to query $\mathcal{SO}^l$ is bounded: 
    $
        C_l \leq M_c 2^{cl}.
    $
\end{enumerate}
\end{assumption}

Assumption \ref{assumption:general} characterizes the requirement of the stochastic oracle.
Assumptions \ref{assumption:general}\ref{assumption:general:I} and \ref{assumption:general:III} imply that obtaining an unbiased gradient estimator of more accurate approximation functions $F^l$ (namely smaller bias from the true gradient) requires higher cost.  The variance decay of the difference estimator $H^l(x,\zeta^l)$ in  Assumption~\ref{assumption:general}\ref{assumption:general:II} is the key assumption for MLMC methods. Intuitively, since $\nabla F^l(x) - \nabla F^{l-1}(x)$ becomes very small for large $l$, constructing highly correlated estimators of $\nabla F^l(x)$ and $\nabla F^{l-1}(x)$ using the same samples $\zeta^l$ will likely yield small variance of their difference. This is akin to building variance reduction in  statistical estimation.  We give several examples of constructing $H^l(x)$ for various optimization and machine learning applications in Section \ref{section:applications}.

Although we format the assumptions in terms of exponential decay in the bias and exponential increase in the costs, in various applications, there exist approximations $\{F_k(x)\}_k$ with a polynomial decaying bias at a polynomial increasing cost $k$, e.g., $|F_k(x)-F(x)|\leq k^{-1}$. In such circumstances, one can pick a subsequence $\{F^l(x)\}_l$, such that $F^l(x): = F_{2^l}(x)$, to ensure  that $|F^l(x) - F(x)|\leq 2^{-l}$ and that the costs increases exponentially. We consider exponential growth in bias, variance, and costs as it simplifies the analysis and holds in various applications.

In the nonconvex setting, since we care about $\eps$-stationarity points rather than $\eps$-optimal function values, we use the following assumption on the approximation error of the gradient to replace Assumption \ref{assumption:general}\ref{assumption:general:I}.
\begin{customthm}{1$(\mathrm{I}$$^\prime)$}
\label{assumption:gradient_bias}
There exists $M_a>0$ such that for any $l\in\NN$, $\|\nabla F^l(x) - \nabla F(x)\|_2^2\leq B_l:= M_a 2^{-al}.$
\end{customthm}
\begin{remark}[{Relationship between
\label{remark:PL}
Assumptions~\ref{assumption:general}\ref{assumption:general:I} and \ref{assumption:gradient_bias}}] Our main results in the strongly convex case can extend to the PL condition~\citep{karimi2016linear} easily, under which Assumptions \ref{assumption:general}\ref{assumption:general:I} and \ref{assumption:gradient_bias} are exchangeable for biased gradient-based methods to achieve the same convergence and total cost results. 
See detailed discussions in Appendix~\ref{appendix:pl-condition}.
In general, Assumptions \ref{assumption:general}\ref{assumption:general:I} and \ref{assumption:gradient_bias} do not imply each other but reflect two ways of constructing approximations. 
Constructing uniform function approximation first and then taking the gradient is more suitable for a problem when the form of the true gradient is unknown.
On the other hand, directly constructing gradient approximations is more suitable when one knows the form of the true gradient and the approximated gradient does not have an easily computable objective function.
\end{remark}

Besides, we make the following standard assumptions on the objective function and its approximations.
\begin{assumption} Either of the following conditions holds:
\label{assumption:objective}
\begin{enumerate}
    \item\label{assumption:objective:I} The objective $F$ and its approximation $F^l$ are $S_F$-smooth for any $l\in \NN$.
    \item\label{assumption:objective:II} The objective $F$ and its approximation $F^l$ are $L_F$-Lipschitz continuous for any $l\in\NN$.
\end{enumerate}
\end{assumption}

\vspace{-0.8em}
\subsection{SGD with MLMC Gradient Estimators}
In the sequel, we construct a direct biased gradient estimator and four different MLMC gradient estimators based on the oracle $\mathcal{SO}^l$.  
They will be fed into the generic SGD framework in Algorithm~\ref{alg:SGD-framework}. These gradient estimators are as follows:
\begin{subequations}\label{alg:general:MLMC}
\noindent
\textbf{$L$-SGD estimator:}  at a query point $x$,  query oracle $\mathcal{SO}^L$ for $n_L$ times to obtain $\{ h^L(x,\zeta_i^L)\}_{i=1}^{n_L}$, and then construct 
    \begin{equation}
    \begin{array}{c}
    \label{alg:L-SGD}
    v^\textrm{$L$-SGD}(x) := \frac{1}{n_L}{\sum}_{i=1}^{n_L} h^L(x,\zeta_i^L).
    \end{array}
    \end{equation}

\noindent
\textbf{V-MLMC estimator:} at a query point $x$, query the oracle $\mathcal{SO}^l$ for $n_l$ times to obtain $\{ H^l(x,\zeta_i^l)\}_{i=1}^{n_l}$ for $l=0,...,L$, and then construct
    \begin{equation}
    \label{alg:V-MLMC}
    \begin{array}{c}
    v^\textrm{V-MLMC}(x): = \sum_{l=0}^L  \frac{1}{n_l} \sum_{i=1}^{n_l} H^l(x,\zeta_i^l).
    \end{array}
    \end{equation}

\noindent
\textbf{RT-MLMC estimator}: at a query point $x$, first sample a random level $\iota$ according to the probability distribution $Q_1= \{q_l\}_{l=0}^L$ with $P(\iota=l)=q_l$, where $\sum_{l=0}^L q_l =1$, and then query the oracle $\mathcal{SO}^\iota$ once to obtains $H^\iota(x,\zeta^\iota)$.
	\begin{equation}
    \label{alg:RT-MLMC}
     v^\textrm{RT-MLMC}(x):= q_\iota^{-1} H^\iota(x,\zeta^\iota).
    \end{equation}

\noindent
\textbf{RU-MLMC estimator}: at a query point $x$, first sample a random level $\iota$ according to the probability distribution $Q_2= \{q_l\}_{l=0}^\infty$ with $P(\iota=l)=q_l$,  where $\sum_{l=0}^\infty q_l =1$, and then query the oracle $\mathcal{SO}^\iota$ once to obtains $H^\iota(x,\zeta^\iota)$.
	\begin{equation}
    \label{alg:RU-MLMC}
    v^\textrm{RU-MLMC}(x) := q_\iota^{-1} H^\iota(x,\zeta^\iota).
    \end{equation}

\noindent
\textbf{RR-MLMC estimator}: at a query point $x$, first sample a random level $L$ according to the probability distribution $Q_2= \{q_l\}_{l=0}^\infty$ with $P(L=l)=q_l$,  where $\sum_{l=0}^\infty q_l =1$,  and then query oracle $\mathcal{SO}^l$ once for $l=0,..,L$ to obtain $\{H^l(x,\zeta^l)\}_{l=0}^L$.
	\begin{equation}
    \label{alg:RR-MLMC}
    \begin{array}{c}
    v^\textrm{RR-MLMC}(x): = \sum_{l=0}^L p_l H^l(x,\zeta^l), 
    \end{array}
    \end{equation}
where $p_l = (1-\sum_{l^\prime=0}^{l-1} q_{l^\prime})^{-1}$,  and $\sum_{l^\prime=0}^{-1} q_{l^\prime}:=0$.
\end{subequations}
\begin{algorithm}[t]
	\caption{SGD Framework}
	\label{alg:SGD-framework}
	\begin{algorithmic}[1]
		\REQUIRE  Number of iterations $ T $, stepsizes $\{\gamma_t\}_{t=1}^{T}$, initialization point $x_1$.
		\FOR{$t=1$ to $T$ \do} 
		\STATE Construct a gradient estimator $v(x_t)$ of $\nabla F(x_t)$.
		\STATE Update $ x_{t+1}=x_t-\gamma_t v(x_t)$.
		\ENDFOR
		\ENSURE 
		$x_T$  if $F(\cdot)$ is strongly convex; otherwise, $\hat x_T$ selected uniformly at random from
		$\{x_t\}_{t=1}^T$.
	\end{algorithmic}
\end{algorithm}

Both RU-MLMC and RR-MLMC are unbiased gradient estimators of $F$, whereas the other three estimators are biased. 
As a result, combining the SGD algorithm with $L$-SGD estimator,  V-MLMC estimator, and RT-MLMC estimator only leads to an approximate optimal point (approximate stationarity in the nonconvex case) of the approximation function $F^L(x)$.

For ease of notation, let $\hat x^{\texttt{A}}_T$ denote the output of SGD using estimator $\texttt{A}$ after $T$ iterations for 
$\texttt{A}\in\{
\text{$L$-SGD}, \text{V-MLMC}, \text{RT-MLMC}, \text{RU-MLMC}, \text{RR-MLMC}
\}$. Under Assumption~\ref{assumption:general}, we have the decomposition of errors in the (strongly) convex case:
\begin{equation}
\label{eq:error_decomposition}
\begin{split}
        &\EE \Big[\underbrace{F(\hat x^\texttt{A}_T) -  F(x^*)}_{\text{Error of Algorithm $\texttt{A}$ on } F}\Big] \\=& 
        \EE \Big[
            \underbrace{F(\hat x^\texttt{A}_T) - F^L(\hat x^\texttt{A}_T)}_{\text{Approximation error}~\leq ~B_L}
        + 
            \underbrace{F^L(\hat x^\texttt{A}_T)-F^L(x^L)}_{\text{Error of SGD on } F^L}\Big] 
        +
            \underbrace{F^L(x^L)-F^L(x^*)}_{\leq ~0 \text{ by optimality of } x^L}
        +
            \underbrace{F^L(x^*)-F(x^*)}_{\text{Approximation error}~\leq ~B_L},
\end{split}
\end{equation}
where $x^L$ is a minimizer of $F^L(x)$. 
The decomposition~\eqref{eq:error_decomposition} suggests that when the level $L$ is large enough, e.g., $L=\lceil a^{-1}\log(4M_a\eps^{-1})\rceil$ such that $2B_L\leq \eps/2$, and number of iterations $T$ is large enough such that expected error of SGD on $F^L$ is smaller than $\eps/2$, these methods return an $\eps$-optimal solution in the convex case.  Let $C_\textrm{iter}^{\texttt{A}}$ and $\mathrm{Var}(v^{\texttt{A}})$ denote the expected computation cost and variance of the estimator $\texttt{A}$, respectively.

\begin{remark}[Intuition on the superiority of MLMC gradient over $L$-SGD]
By the analysis of SGD (see Appendix \ref{app_sec:SGD}), the number of iterations $T
\propto \Var(v^\mathtt{A}(x))$. 
Then, the total cost of the gradient method $\texttt{A}$ is $T\cdot c_\mathrm{iter}^\mathtt{A}\propto \Var(v^\mathtt{A}(x))\cdot c_\mathrm{iter}^\mathtt{A}$, where $\propto$ denotes proportional to.
It is noteworthy from Table~\ref{tab:bias-variance-cost} that $L$-SGD has the larger value of $\Var(v^\mathtt{A}(x))\cdot c_\mathrm{iter}^\mathtt{A}$ over all four MLMC gradient methods.
Next we discuss why $L$-SGD admits a larger $\Var(v^\mathtt{A}(x))\cdot c_\mathrm{iter}^\mathtt{A}$ using RT-MLMC for illustration. The major bottleneck of $L$-SGD is that it requires querying the expensive stochastic oracle $\mathcal{SO}^L$ at each iteration. Notice that the RT-MLMC gradient estimator is an unbiased estimator of $L$-SGD. Thus, their bias levels towards $\nabla F(x)$ remain the same. However, RT-MLMC queries $\mathcal{SO}^l$ with small $l$ with a high probability and queries $\mathcal{SO}^l$ with large $l$ with a low probability. Thus, the expected costs to query the oracle are greatly reduced from $C_L$ to $\sum_{l=0}^L q_l C_l$. One may question that for large $l$, $q_l$ is small, and RT-MLMC divides $q_l$, which may introduce a large variance. In fact, since $H^l(x)$ admits an exponentially decreasing variance, which mitigates the increase of variance introduced by dividing $q_l$. 
\end{remark}

We summarize the bias, variance, and (expected) computation cost in Table \ref{tab:bias-variance-cost} when  $x\in\RR^d$ is
independent from $v^{\texttt{A}}$.
See Lemmas~\ref{lm:bias} and \ref{lm:variance_cost} for a detailed discussion on parameter selection and the corresponding bias, variance, and cost.

\begin{table}[!ht]
\vspace{-1em}
    {
    \footnotesize
    \caption{Bias, Variance, and Cost of Gradient Estimators.}
    \renewcommand\arraystretch{1.3}
    \begin{center}
    \begin{tabular}{c|c|c|c|c}
        \toprule[1pt]
        {Estimators $\texttt{A}$}
        & Expectation $\EE v^{\texttt{A}}(x)$
        & Variance $\mathrm{Var}(v^{\texttt{A}}(x))$
        & Cost $C_\textrm{iter}^{\texttt{A}}$
        & $\mathrm{Var}(v^{\texttt{A}}(x))$ $\times$ $C_\textrm{iter}^{\texttt{A}}$ when $b>c$
        \\
        \hline
        
        {$L$-SGD}
        & $\nabla  F^L(x)$
        & $n_L^{-1}\sigma^2$
        & $n_L C_L$
        & $\cO(\eps^{-c/a})$
        \\
        \hline
        
        {V-MLMC}
        & $\nabla F^L(x)$
        & ${\sum}_{l=0}^L n_l^{-1}V_l$
        & ${\sum}_{l=0}^L n_l C_l$
        & $\cO(1)$
        \\
        \hline
        
        {RT-MLMC}
        & $\nabla F^L(x)$
        & ${\sum}_{l=0}^L q_l^{-1}V_l$
        & ${\sum}_{l=0}^L q_l C_l$
        & $\cO(1)$
        \\
        \hline

        {RU-MLMC}
        & $\nabla F(x)$
        & ${\sum}_{l=0}^\infty q_l^{-1}V_l$
        & ${\sum}_{l=0}^\infty q_l C_l$
        & $\cO(1)$
        \\
        \hline
        
        {RR-MLMC}
        & $\nabla F(x)$
        & $\sum_{L=0}^\infty q_L \left(\sum_{l=0}^L p_l^2 V_l\right)$
        & $\sum_{L=0}^\infty q_L\left(\sum_{l=0}^L C_l\right)$
        & $\cO(1)$
        \\
        \bottomrule[1pt]
    \end{tabular} 
    \label{tab:bias-variance-cost}
    \end{center}
    }\vspace{-1em}
    \end{table}

\subsection{Variance Reduced Methods with MLMC Gradient Estimators}
In the nonconvex smooth regime with unbiased oracles, variance reduction methods, such as SPIDER~\citep{fang2018spider,wang2018spiderboost}, SARAH~\citep{nguyen2017sarah}, STORM~\citep{STORM}, and PAGE~\citep{li2021page}, can reduce the complexity for finding approximate stationary points.
We adopt the algorithm idea from SPIDER to develop variance-reduced algorithms using biased oracles and further demonstrate their (expected) total costs to achieve $\eps$-stationary points.
We use the prefix VR to denote variance reduction methods using gradient estimators $L$-SGD, V-MLMC, RT-MLMC, RU-MLMC, and RR-MLMC. Algorithm~\ref{alg:VR-framework} is the variance-reduced algorithmic framework. Note that the variance reduction effect uses the control variate technique: $v_k(x_t,\xi_{tk})$ and $v_k(x_{t-1},\xi_{tk})$ are built from the same samples in the gradient estimator construction.

\begin{algorithm}[ht]
	\caption{Variance Reduction Counterpart of SGD Framework}
	\label{alg:VR-framework}
	\begin{algorithmic}[1]
		\REQUIRE  Number of iterations $ T $, stepsizes $\gamma$, initialization point $x_1$, batch size $D_1$, $D_2$, epoch length $Q_E$.
		\FOR{$t=1$ to $T$ \do} 
		\STATE Construct the recursive gradient estimator
		\begin{equation}
  \label{eq:variance_reduction_alg}
        m_t(x_t) = 
        \begin{cases} \frac{1}{D_1}\sum_{k=1}^{D_1}v_k(x_t,\xi_{tk}), & \text{ if } t\equiv0(\text{mod }Q_E),\\
       m_{t-1}(x_{t-1}) - \frac{1}{D_2}\sum_{k=1}^{D_2}v_k(x_{t-1},\xi_{tk})+ \frac{1}{D_2}\sum_{k=1}^{D_2}v_k(x_t,\xi_{tk}), & \text{ otherwise,} 
        \end{cases}
        \end{equation}
        where $\{v_k(x_t,\xi_{tk})\}_{k=1}^D$ are i.i.d gradient estimators of $\nabla F(x_t)$ and $\{v_k(x_{t-1},\xi_{tk})\}_{k=1}^D$ are i.i.d gradient estimators of $\nabla F(x_{t-1})$. 
		\STATE Update $ x_{t+1}=x_t-\gamma m_t(x_t)$.
		\ENDFOR
		\ENSURE $\hat x_T$ uniformly selected from
		$\{x_t\}_{t=1}^T$.
	\end{algorithmic}
\end{algorithm}
\vspace{-1em}
We use VR RT-MLMC to demonstrate how we construct $m_t(x_t)$. Let $\{v_k(x,\xi_k)\}_{k=1}^D$ be independent RT-MLMC gradient estimators constructed via independently generating level $l_k\in\{ 0,...,L \}$ following the distribution $Q$ for $k=1,...,D$, obtaining $\{H^{l_k}(x,\zeta^{l_k})\}_{k=1}^D$ by querying $\mathcal{SO}^{l_1},\ldots,\mathcal{SO}^{l_D}$ and constructing 
$$
\frac{1}{D}\sum_{k=1}^{D}v_k(x,\xi_k) = \frac{1}{D}\sum_{k=1}^{D} q_{l_k}^{-1}H^{l_k}(x;\zeta^{l_k}),
$$
where $D$ refers to $D_1$ or $D_2$ depending on the iteration $t$.
Then $m_t(x)$ is constructed by Algorithm~\ref{alg:VR-framework}. As discussed earlier, the variance reduction effect is important to ensure exponential decaying variance of  $H^{l_k}$. Note that the variance reduction in Algorithm \ref{alg:VR-framework} is achieved via a control variate in $v_k(x_t,\xi_{tk}) - v_k(x_{t-1},\xi_{tk})$ in \eqref{eq:variance_reduction_alg}. Thus the VR RT-MLMC methods in fact leverages the control variate technique twice, one in estimation of $\nabla F^l(x) - \nabla F^{l-1}(x)$ and one in the construction of the recursive gradient estimator. The effect are also different, where the former aims to reduce the oracle costs for building low-cost low variance gradient estimators and the latter aims to reduce the number of oracle queries needed.

\section{(Expected) Total Cost Analysis}
\label{section:total_cost}

In this section, we demonstrate the (expected) total cost of gradient estimators in \eqref{alg:general:MLMC} under different convexity conditions, different smoothness conditions, and different combinations of $a,b,c$. We further discuss the applicability of MLMC methods. 
We define the following notations for the simplicity of the presentation.
Let $\hat x_T^{\texttt{A}}$ be selected uniformly from $\{x_t^{\texttt{A}}\}_{t=1}^T$ for $\texttt{A}\in\{
\text{$L$-SGD}, \text{V-MLMC}, \text{RT-MLMC}, \text{RU-MLMC}, \text{RR-MLMC}
\}$. 
Recall that $C_\textrm{iter}^{\texttt{A}}$  denotes the expected computation cost associated with the estimator $\texttt{A}$. Let $T^{\texttt{A}}$ denote the iteration complexity of ${\texttt{A}}$ for achieving $\eps$-optimality or $\eps$-stationarity,
  and  $C:=T^{\texttt{A}}C_\mathrm{iter}^{\texttt{A}}$ denote the (expected) total cost.  Denote the upper bounds on $\mathrm{Var}(v^{\texttt{A}}(x))$ as $\bV(v^{\texttt{A}})$.

\vspace{-0.5em}
\subsection{Total Cost of \texorpdfstring{$L$-SGD}{}}
In this subsection, we derive the total cost of the $L$-SGD method, which covers the biased gradient methods in various applications. 
We first show the variance and the cost of the $L$-SGD gradient estimator.

\begin{lemma}[Variance and Per-Iteration Cost of $L$-SGD]
\label{lm:variance_cost_LSGD}
Under Assumption \ref{assumption:general}, for any $x\in\RR^d$, the variance and per-iteration cost of $L$-SGD estimator $v^{\LSGD}(x)$ with batch size $n_L$ satisfy
\begin{equation*}
        \mathrm{Var}(v^{\LSGD}(x)) 
    \leq 
        \frac{\sigma^2}{n_L};
    \quad
        C_\mathrm{iter}^{\LSGD} 
    \leq
        n_L M_c 2^{cL}.     
\end{equation*}
\end{lemma}
In the (strongly) convex case, we set $L=\lceil a^{-1}\log(4M_a\eps^{-1})\rceil$ so that $2B_L\leq \eps/2$. Without loss of generality, consider the mini-batch size $n_L=1$. Then  the cost of $L$-SGD estimator is $M_c2^{cL}=\cO(\eps^{-c/a})$ and the variance is of order $\cO(1)$.
The total cost $C=T\cdot C_\mathrm{iter}^{\LSGD} $. From  the analysis of SGD (See Theorem \ref{lm:SGD} in Appendix \ref{app_sec:SGD}), to ensure $\eps$-optimality, the iteration $T=\cO(\epsilon^{-1})$ when $F^L$ is strongly convex, and $T=\cO(\epsilon^{-2})$ when $F^L$ is convex.  Therefore,  the total costs for $L$-SGD to achieve an $\eps$-optimal solution are $\cO(\eps^{-1-c/a})$  and $\cO(\eps^{-2-c/a})$ in the strongly convex and convex settings, respectively.  The following theorem formally describes the complexity bound, whose proof is in Appendix~\ref{appendix:L-SGD}.
\begin{theorem}[Total cost of $L$-SGD]
\label{thm:L-SGD}
For $L$-SGD with batch size $n_L=1$,%
with properly chosen hyper-parameters as in Table~\ref{tab:assumptions:LSGD}, we have the following results.
\begin{enumerate}
    \item 
Suppose $F^L$ is $\mu$-strongly convex and $S_F$-smooth and Assumption~\ref{assumption:general} holds, 
the total cost of $L$-SGD for finding an $\eps$-optimal solution of $F$ is $\cO(\eps^{-1-c/a}).$
    \item
 Suppose $F^L$ is convex and Assumptions \ref{assumption:general} holds, 
 under either Assumption~\ref{assumption:objective}\ref{assumption:objective:I} or \ref{assumption:objective}\ref{assumption:objective:II}, 
 the total cost of $L$-SGD for finding an $\eps$-optimal solution of $F$ is $\cO( \eps^{-2-c/a}).$
    \item
Suppose $F^L$ is $S_F$-smooth, and Assumptions \ref{assumption:general}\ref{assumption:general:II}\ref{assumption:general:III} and \ref{assumption:gradient_bias} hold, the total cost of $L$-SGD for finding an $\eps$-stationary point of $F$ is $\cO(\eps^{-4-2c/a}).$
\end{enumerate}
\vspace{-1em}
\begin{table}[!ht]
\caption{Hyper-parameters of $L$-SGD in Theorem~\ref{thm:L-SGD}}
\centering
\small{
\begin{tabular}{c|c|c|c}
	\toprule%
\bf{Assumptions on $F^L$} & \bf{Level $L$} & \bf{Step Size $\gamma_t$} & \bf{Iteration $T$} \\[8pt]
\hline
$\begin{array}{c}
    \mbox{$\mu$-strongly convex}\\
    \mbox{ and $S_F$-smooth}
\end{array}$
& $\lceil \frac{\log(4M_a\eps^{-1})}{a}  \rceil$&  $
\Big[t+S_F^2/\mu^2\Big]^{-1}
$
& $\cO(\eps^{-1})$
\\[8pt]
\hline
convex &$\lceil \frac{\log(4M_a\eps^{-1})}{a}  \rceil$ &  $\left\{ \begin{aligned}
    \Big[T\sigma^2\Big]^{-1/2}
    ,&\quad\text{under Assumption~\ref{assumption:objective}\ref{assumption:objective:I}}\\
   \Big[T(\sigma^2+L_f^2)\Big]^{-1/2}
   , &\quad\text{under Assumption~\ref{assumption:objective}\ref{assumption:objective:II}}
\end{aligned}\right.$
&
$\cO(\eps^{-2})$
\\[8pt]
\hline
$S_F$-smooth &$\lceil \frac{\log(4M_a\eps^{-2})}{a}  \rceil$ &  $\Big[T\sigma^2\Big]^{-1/2}$
&
$\cO(\eps^{-4})$
\\[8pt]
\bottomrule
\end{tabular}}
\label{tab:assumptions:LSGD}
\vspace{-2em}
\end{table}
\end{theorem}

\vspace{-0.5em}
\subsection{Expected Total Cost of RT-MLMC}
RT-MLMC is an unbiased estimator of $L$-SGD using randomization, and thus, RT-MLMC is unbiased for $\nabla F^L(x)$. We obtain the bias of RT-MLMC directly through Assumption~\ref{assumption:general}\ref{assumption:general:I} or \ref{assumption:gradient_bias}.
For $\alpha\in\mathbb{R}$,
denote  
\begin{equation}
R^L(\alpha):=\sum_{l=0}^L 2^{\alpha l}  
    = 
        \frac{1-2^{\alpha(L+1)}}{1-2^{\alpha}} \text{ for } \alpha \not = 0; \quad R^\infty(\alpha):=\lim_{L\rightarrow \infty} R^L(\alpha)=(1-2^{\alpha})^{-1} \text{ for } \alpha < 0.
\label{Eq:RL:infty}
\end{equation}
\vspace{-1em}
\begin{lemma}[Variance and Per-Iteration Cost of RT-MLMC]
\label{lm:variance_cost_RT-MLMC}
Under Assumption \ref{assumption:general}, construct a distribution $Q=\{q_l\}_{l=0}^l$ with the probability mass value $q_l = 2^{-(b+c)l/2}R^L(-\frac{b+c}{2})^{-1}$ and $\sum_{l=0}^L q_l =1$. The variance and the expected per-iteration cost of RT-MLMC satisfy:
\[
\begin{array}{lll}
\text{If $c\ne b$:}&~~
 \displaystyle \mathrm{Var}(v^{\RTMLMC}(x) ) 
    \leq M_bR^L\big(\frac{c-b}{2}\big)R^L\big(-\frac{b+c}{2}\big),   
&~~ \displaystyle C_\mathrm{iter}^{\RTMLMC} \le M_cR^L\big(\frac{c-b}{2}\big)R^L\big(-\frac{b+c}{2}\big)^{-1}, \\
\text{If $c= b$:}&~~
 \displaystyle   \mathrm{Var}(v^{\RTMLMC}(x) ) 
    \leq M_b(L+1) R^L\big(-\frac{b+c}{2}\big),  & ~~ \displaystyle C_\mathrm{iter}^{\RTMLMC} \le M_c (L+1) R^L\big(-\frac{b+c}{2}\big)^{-1}.
\end{array}
\]
\end{lemma}
To ensure that the bias is of order $\cO(\eps)$, the lemma implies that $\mathrm{Var}(v^{\RTMLMC}(x))$ and $C_\mathrm{iter}^\mathrm{RT-MLMC}$ are of order $\cO(1)$ when $b>c$,  and $\cO(\log(\eps^{-1}))$ when $b=c$. When $L$-SGD picks the mini-batch $n_L=1$, the variance of RT-MLMC and $L$-SGD is of the same scale, yet the cost to build RT-MLMC gradient estimator reduces from $\cO(\eps^{-c/a})$ for $L$-SGD to $\cO(1)$.

\begin{remark}[{Construction of the  distribution $Q := \{q_l\}_{l=0}^L$}]
We use {RT-MLMC} for convex objectives as an example to explain how $Q$ is constructed.
To minimize the (expected) total cost for achieving $\eps$-optimality, conceptually, we solve the following optimization problem:
\begin{equation}
\label{problem:minimize_cost}
\min_{Q}\Big\{ 
 T^{\RTMLMC} C_\mathrm{iter}^{\RTMLMC}:\quad 
 \sum_{l=0}^L q_l = 1,
 q_l\geq 0,  \forall l=0,\ldots,L
\Big\},
\end{equation}
where the constraints are to make sure $Q$ is a well-defined probability mass vector. 
Although Problem~\eqref{problem:minimize_cost} is nonconvex and intractable to solve, we can show that under properly selected stepsizes, %
\begin{equation}
\label{eq:selection_proportation}
    T^{\RTMLMC} C_\mathrm{iter}^{\RTMLMC}\propto \mathrm{Var}(v^{\RTMLMC})C_\mathrm{iter}^{\RTMLMC} =\left(\sum_{l=0}^L q_l C_l\right)\left(\sum_{l=0}^L  V_l/q_l\right).
\end{equation}
To ensure small bias, we take $L=\cO(\log(\eps^{-1}))$.
Optimizing the right-hand-side of \eqref{eq:selection_proportation} over $Q$ gives $q_l \propto 2^{-(b+c)l/2}$. Although it might not be the optimal solution of \eqref{problem:minimize_cost}, we show that such a choice of $q_l$ can significantly reduce the expected total cost compared to $L$-SGD. As we will show later, when $b>c$, the expected total cost of RU-MLMC and RT-MLMC reaches the lower bounds of unbiased stochastic optimization, implying the distribution $Q$ is an optimal choice in terms of the dependence on the accuracy $\eps$.
\end{remark}

\begin{theorem}[Expected Total cost of RT-MLMC]
\label{thm:RT-MLMC}
For RT-MLMC with a distribution $Q=\{q_l\}_{l=0}^L$ such that $q_l \propto 2^{-(b+c)l/2}$, 
with properly chosen hyper-parameters as in Table~\ref{tab:assumptions:RTMLMC},%
\begin{enumerate}
    \item 
Suppose $F^L$ is $\mu$-strongly convex and $S_F$-smooth and Assumption~\ref{assumption:general} holds, the total cost of RT-MLMC for finding an $\eps$-optimal solution of $F$ is 
\[%
\begin{array}{llllll}
 \cO(\eps^{-1}),    & \quad\text{if $c<b$;} & \qquad \cO(\eps^{-1}\log^2(1/\eps)), & \quad\text{if $c=b$;}
 & \qquad\cO( \eps^{-1-(c-b)/a} ), &\quad \text{if $c>b$.}
\end{array}
\]
  \item
Suppose $F^L$ is convex and Assumption~\ref{assumption:general} holds, 
under either Assumption~\ref{assumption:objective}\ref{assumption:objective:I} or \ref{assumption:objective}\ref{assumption:objective:II},
the total cost of RT-MLMC for finding an $\eps$-optimal solution of $F$ is 
\[%
\begin{array}{llllll}
 \cO(\eps^{-2}),    & \quad\text{if $c<b$;} & \qquad \cO(\eps^{-2}\log^2(1/\eps)), & \quad\text{if $c=b$;}
 & \qquad\cO( \eps^{-2-(c-b)/a} ), &\quad \text{if $c>b$.}
\end{array}
\]
   \item
Suppose $F^L$ is $S_F$-smooth and Assumptions \ref{assumption:general}\ref{assumption:general:II}\ref{assumption:general:III} and \ref{assumption:gradient_bias} hold,
the total cost of RT-MLMC for finding an $\eps$-stationarity point of $F$ is 
\[
\begin{array}{llllll}
 \cO(\eps^{-4}),    & \quad\text{if $c<b$;} & \qquad \cO(\eps^{-4}\log^2(1/\eps)), & \quad\text{if $c=b$;}
 & \qquad\cO( \eps^{-4-2(c-b)/a} ), &\quad \text{if $c>b$.}
\end{array}
\] 
\end{enumerate}
\vspace{-1.8em}
\begin{table}[!ht]
\caption{Hyper-parameters of RT-MLMC in Theorem~\ref{thm:RT-MLMC}}
\centering
\small{
\begin{tabular}{c|c|c|c}
	\toprule%
\bf{Assumptions on $F^L$} & \bf{Level $L$} & \bf{Step Size $\gamma_t$} & \bf{Iteration $T$} \\[8pt]
\hline
$\begin{array}{c}
    \mbox{$\mu$-strongly convex}\\
    \mbox{ and $S_F$-smooth}
\end{array}$
& $\lceil \frac{\log(4M_a\eps^{-1})}{a}  \rceil$&  $
\Big[t+S_F^2/\mu^2\Big]^{-1}
$
& $\cO(\eps^{-1})$
\\[8pt]
\hline
convex &$\lceil \frac{\log(4M_a\eps^{-1})}{a}  \rceil$ &  $\left\{ \begin{aligned}
    \Big[\bV(v^\mathrm{RT-MLMC})T\Big]^{-1/2}
    ,&\quad\text{under Assumption~\ref{assumption:objective}\ref{assumption:objective:I}}\\
   \Big[(\bV(v^\mathrm{RT-MLMC})+L_F^2)T\Big]^{-1/2}
   , &\quad\text{under Assumption~\ref{assumption:objective}\ref{assumption:objective:II}}
\end{aligned}\right.$
&
$\cO(\eps^{-2})$
\\[8pt]
\hline
$S_F$-smooth &$\lceil \frac{\log(4M_a\eps^{-2})}{a}  \rceil$ &  $\Big[\bV(v^\mathrm{RT-MLMC})T\Big]^{-1/2}$
&
$\cO(\eps^{-4})$
\\[8pt]
\bottomrule
\end{tabular}}
\label{tab:assumptions:RTMLMC}%
\vspace{-1.8em}
\end{table}
\end{theorem}
The proof of Theorem~\ref{thm:RT-MLMC} is in Appendix~\ref{appendix:RT-MLMC}.
When $b>c$, RT-MLMC achieves an $\cO(1)$ expected per-iteration cost and $\cO(1)$ variance, and thus the total cost of RT-MLMC is the same as the costs for classical SGD with unbiased oracles, i.e., RT-MLMC reaches the lower bounds for first-order methods for stochastic optimization with unbiased oracles.
Compared to $L$-SGD, RT-MLMC is always  $\tilde \cO(\eps^{-c/a})$ better in terms of complexity bounds in the (strongly) convex case no matter what is the relationship of $a$, $b$, and $c$.
\vspace{-0.5em}
\subsection{Total Cost of V-MLMC}
We show that V-MLMC achieves a reduced total cost compared to $L$-SGD  under some additional assumptions on the relationship between $a,b,c$ and large mini-batch size. In the main context, we only demonstrate the case when $b\geq c$~(See Theorem~\ref{thm:V-MLMC}). The full version and its formal proof are in Appendix~\ref{appendix:V_MLMC}.
As usual, we first demonstrate the variance and the per-iteration cost of V-MLMC. 
\begin{lemma}[Variance and Per-Iteration Cost of V-MLMC]
\label{lm:variance_cost_V-MLMC}
Under Assumption \ref{assumption:general} and setting the batch size $n_l = \lceil  2^{-(b+c)l/2} N \rceil, l=0,\ldots,L$ for a constant $N>0$, the variance and  per-iteration cost satisfy%
\[
\begin{array}{lll}
\text{If $c\ne b$: }&~~
 \displaystyle \mathrm{Var}(v^{\VMLMC}(x)) 
    \leq M_b R^L\big(\frac{c-b}{2}\big) N^{-1},   
&~~ \displaystyle C_\mathrm{iter}^{\VMLMC} \le M_cR^L\big(\frac{c-b}{2}\big)N + M_cR^L(c), \\
\text{If $c= b$: }&~~
 \displaystyle  \mathrm{Var}(v^{\VMLMC}(x)) 
    \leq M_b(L+1) N^{-1},  & ~~ \displaystyle C_\mathrm{iter}^{\VMLMC} \le M_c (L+1)N+M_c R^L(c).
\end{array}
\]
\end{lemma}
Similar to RT-MLMC, the optimal selection of $\{n_l\}_{l=0}^L$ can be obtained via solving the V-MLMC counterpart of \eqref{problem:minimize_cost}. However, it would become an integer program since one needs to ensure that the batch sizes $\{n_l\}_{l=0}^L$ are integers. Here we consider a relaxation by setting $n_l= q_l\cdot N$, where $\sum_{l=0}^L q_l=1$ and $q_l\geq 0$.  As a result, we have $n_l\propto 2^{-(b+c)l/2} N$. Given that $n_l$ is the mini-batch size, we set $n_l = \lceil 2^{-(b+c)l/2} N\rceil$. Unfortunately, rounding the optimal selection $2^{-(b+c)l/2} N$ to integers would incur extra costs, i.e., the term $M_cR^L(c)$ in Lemma \ref{lm:variance_cost_V-MLMC}. We show later that the extra cost would lead to extra requirements on $a,b,c$, and $N$ for V-MLMC to achieve a reduced total cost as other MLMC methods compared to the $L$-SGD method.

\begin{theorem}[Total cost for V-MLMC when $b\geq c$]
\label{thm:V-MLMC}
For V-MLMC with batch sizes $n_l = \lceil 2^{-(b+c)l/2}N\rceil$ for some $N>0$, 
with properly chosen hyper-parameters as in Table~\ref{tab:assumptions:VMLMC}, %
\begin{enumerate}
    \item 
Suppose $F^L$ is $\mu$-strongly convex and $S_F$-smooth and Assumption~\ref{assumption:general} holds, 
the total cost of V-MLMC for finding an $\eps$-optimal solution of $F$ is $\tilde \cO( \eps^{-\max\{1,c/a\}}).$
    \item
Suppose $F^L$ is convex and $S_F$-smooth, and Assumptions \ref{assumption:general} holds, 
 the total cost of V-MLMC for finding an $\eps$-optimal solution of $F$ is $\tilde \cO( \eps^{-1-\max\{1,c/a\}}).$
    \item
Suppose $F^L$ is $S_F$-smooth, and Assumptions \ref{assumption:general}\ref{assumption:general:II}\ref{assumption:general:III} and \ref{assumption:gradient_bias} hold, the total cost of V-MLMC for finding an $\eps$-stationarity point of $F$ is $\tilde \cO( \eps^{-2-2\max\{1,c/a\}}).$
\end{enumerate}

\begin{table}[!ht]
\vspace{-1em}
\caption{Hyper-parameters of V-MLMC in Theorem~\ref{thm:V-MLMC}}
\centering
\small{
\begin{tabular}{c|c|c|c|c}
	\toprule%
\bf{Assumptions on $F^L$} & \bf{Level $L$} & \bf{Step Size $\gamma_t$} & $\begin{array}{c}
    \textbf{Multiplicative Factor}\\
    \textbf{$N$ of Batch Sizes}
\end{array}$& \bf{Iteration $T$} \\[8pt]
\hline
$\begin{array}{c}
    \mbox{$\mu$-strongly convex}\\
    \mbox{ and $S_F$-smooth}
\end{array}$
& $\lceil \frac{\log(4M_a\eps^{-1})}{a}  \rceil$&  $
1/S_F
$
& 
$\tilde{\cO}(\eps^{-1})$
&
$\cO(\log\frac{1}{\eps})$
\\[8pt]
\hline
convex &$\lceil \frac{\log(4M_a\eps^{-1})}{a}  \rceil$ &  $1/(2S_F)$& 
$\tilde{\cO}(\eps^{-1})$
&
$\cO(\eps^{-1})$
\\[8pt]
\hline
$S_F$-smooth &$\lceil \frac{\log(4M_a\eps^{-2})}{a}  \rceil$ &  $1/S_F$& 
$\tilde{\cO}(\eps^{-2})$
&
$\cO(\eps^{-2})$
\\[8pt]
\bottomrule
\end{tabular}}
\label{tab:assumptions:VMLMC}%
\vspace{-2em}
\end{table}
\end{theorem}
\begin{remark}[Why V-MLMC requires a large mini-batch size]
\label{rem:V-MLMC}
According to Lemma \ref{lm:variance_cost_V-MLMC} and the expression of $n_l$, the per-iteration cost of V-MLMC satisfies
\begin{equation}
\label{eq:cost_iteration_VMLMC}
    C_\text{iter}^\text{V-MLMC} \leq \sum_{l=0}^L 2^{(c-b)l/2}N +\sum_{l=0}^L 2^{cl}\overset{(\text{if } b>c)}{=}\cO(N)+\cO(\eps^{-c/a}) .
\end{equation} 
where $L=\lceil 1/a\log(4M_a\eps^{-1})\rceil$ as specified in Theorem~\ref{thm:V-MLMC}. 
The first term on the right-hand side~(RHS) in \eqref{eq:cost_iteration_VMLMC} denotes the desired balanced per-iteration cost of MLMC, while the second term denotes the cost incurred by rounding the batch size to an integer. When $c< b$ and $N=\cO(1)$, the second term dominates the per-iteration cost. At the same time, since the variance of V-MLMC is $\cO(N^{-1})=\cO(1)$, the iteration complexity and per-iteration cost of V-MLMC are the same as $L$-SGD. Therefore, the total costs of V-MLMC are the same as the total cost of $L$-SGD when $N=\cO(1)$.

To handle the high complexity bound caused by the rounding issue, we use a large  $N=\cO(\eps^{-1})$ to reduce the variance of V-MLMC to $\cO(\eps)$. Thus, V-MLMC gradient method behaves like deterministic gradient descent. The iteration complexity of V-MLMC reduces from $\cO(\eps^{-1})$ to logarithmic level in the strongly convex case and from $\cO(\eps^{-2})$ to $\cO(\eps^{-1})$ in the convex case, and the total cost becomes $\tilde\cO(\eps^{-\max\{1,c/a\}})$ and $\cO(\eps^{-1-\max\{1,c/a\}})$, respectively. In such a case, V-MLMC is always better than $L$-SGD in 
 terms of the total cost by $\cO(\eps^{-\min\{1,c/a\}})$. The drawback is that V-MLMC has to use very large $N$ and thus large mini-batches $n_l$ for small $l$. In contrast, $L$-SGD cannot use large batch to reduce the total cost as $\bV^{L\mathrm{-SGD}}C_\mathrm{iter}^{L\mathrm{-SGD}}=\cO(\eps^{-c/a})$ for any $n_L>0$.
For $b<c$, please refer to the proof of Theorem \ref{thm:V-MLMC} in  Appendix \ref{appendix:V_MLMC}.

\end{remark}

\begin{remark}[{V-MLMC cannot reduce total cost under Lipschitz continuity Condition}]
\label{rem:vmlmc_vs_lsgd}
Note that smoothness is the key for V-MLMC to outperform $L$-SGD. As the total cost is proportional to  $ \mathrm{Var}(v^A)C_\mathrm{iter}^A$, the large mini batches ensures that $\mathrm{Var}(v^\mathrm{V-MLMC})$ is of order $\cO(\eps)$ and $ \mathrm{Var}(v^A)C_\mathrm{iter}^A=\cO(1)$.
However, in the Lipschitz continuous case, the total cost is proportional to $
 (\mathrm{Var}(v^A)+L_F^2)C_\mathrm{iter}^A
$, and even if using a large mini-batch size $N=\cO(\eps^{-1})$, the term $\mathrm{Var}(v^A)+L_F^2$ cannot be reduced. As a result, V-MLMC cannot outperform $L$-SGD in the Lipschitz continuous setting. On the contrary, RT-MLMC does not require a mini-batch and can still outperform $L$-SGD even if the objective is only Lipschitz continuous.
\end{remark}

\vspace{-1em}
\subsection{Expected Total Cost of RU-MLMC and RR-MLMC}
For two unbiased MLMC methods, i.e., RU-MLMC and RR-MLMC, we derive their variance and per-iteration cost following similar steps as in Lemma~\ref{lm:variance_cost_RT-MLMC} and let $L\rightarrow\infty$. However, when $b\leq c$, either the variance or the expected per-iteration cost of them is unbounded for any selection of the distribution $Q$. Hence, the total expected costs of these two unbiased MLMC estimators are unbounded for any selection of $Q$ where $b\leq c$, implying that they are not theoretically applicable unless with stronger assumptions.

Given these, we focus on the case when $b>c$, and we show that the variance and the per-iteration cost of RU-MLMC and RR-MLMC are both $\cO(1)$. The expected total costs of these two methods follow from the analysis of unbiased SGD  (see, e.g.,~\citet{nemirovski2009robust,ghadimi2013stochastic}). We summarize the result in Theorem~\ref{thm:unbiased_MLMC} with the proof in Appendix~\ref{appendix:unbiased}.

\begin{theorem}[Expected Total cost of RU-MLMC and RR-MLMC]
\label{thm:unbiased_MLMC}
Let Assumption~\ref{assumption:general}\ref{assumption:general:II}\ref{assumption:general:III} hold with $b>c$ and Assumption~\ref{assumption:gradient_bias} hold, and define $q_l = 2^{-(b+c)l/2}(1-2^{-(b+c)/2})$ when using RU-MLMC or RR-MLMC estimator.
With properly chosen hyper-parameters as in Table~\ref{tab:assumptions:RUMLMC:RRMLMC},
\begin{enumerate}%
    \item 
    Suppose $F$ is $\mu$-strongly convex and $S_F$-smooth, 
the expected total cost of RU-MLMC or RR-MLMC for finding an $\eps$-optimal solution of $F$ is $\cO(\eps^{-1}).$
    \item
   Suppose $F$ is convex and either Assumption~\ref{assumption:objective}\ref{assumption:objective:I} or \ref{assumption:objective}\ref{assumption:objective:II} holds,
the expected total cost of RU-MLMC or RR-MLMC for finding an $\eps$-optimal solution of $F$ is $\cO(\eps^{-2}).$
    \item Suppose $F$ is $S_F$-smooth~(i.e., Assumption \ref{assumption:objective}\ref{assumption:objective:I} holds), the expected total cost of RU-MLMC or RR-MLMC for finding an $\eps$-stationarity point of $F$ is $\cO(\eps^{-4}).$
\end{enumerate}
\begin{table}[!ht]
\caption{Hyper-parameters of RU-MLMC and RR-MLMC in Theorem~\ref{thm:V-MLMC}.
The notation $\texttt{A}$ refers to RU-MLMC or RR-MLMC.
}
\centering
\small{
\begin{tabular}{c|c|c|c|c}
	\toprule%
\bf{Assumptions on $F$} & \bf{Step Size $\gamma_t$} & \bf{Iteration $T$} \\[8pt]
\hline
$\begin{array}{c}
    \mbox{$\mu$-strongly convex}\\
    \mbox{ and $S_F$-smooth}
\end{array}$
&  $\mu^{-1}(t+2S_F/\mu^2)^{-1}$
&
$\cO(\eps^{-1})$
\\[8pt]
\hline
convex &  $\left\{ \begin{aligned}
    \Big[\bV(v^\texttt{A})T\Big]^{-1/2}
    ,&\quad\text{under Assumption~\ref{assumption:objective}\ref{assumption:objective:I}}\\
   \Big[(\bV(v^\texttt{A})+L_F^2)T\Big]^{-1/2}
   , &\quad\text{under Assumption~\ref{assumption:objective}\ref{assumption:objective:II}}
\end{aligned}\right.$
&
$\cO(\eps^{-2})$
\\[8pt]
\hline
$S_F$-smooth&  $\Big[\bV(v^\texttt{A})T\Big]^{-1/2}$
&
$\cO(\eps^{-4})$
\\[8pt]
\bottomrule
\end{tabular}}
\label{tab:assumptions:RUMLMC:RRMLMC}\vspace{-1em}
\end{table}
\end{theorem}

\subsection{Total Cost of Variance Reduction Methods}
\label{section:variance_reduction}

To demonstrate the convergence of the variance reduction gradient MLMC methods, we require an additional assumption on the stochastic oracle $\mathcal{SO}^l$, called average smooth condition~\citep{arjevani2019lower}. 
\begin{assumption}
\label{assumption:average_smoothness}
For any given $x$, the output of the stochastic oracle $\mathcal{SO}^l$ satisfies the average smooth condition, i.e., there exist constants $S_h,S_H>0$ such that 
$\EE_{\zeta^l} \|h^l(x,\zeta^l)-h^l(y,\zeta^l)\|_2\leq S_h\|x-y\|_2$ and $\EE_{\zeta^l} \|H^l(x,\zeta^l)-H^l(y,\zeta^l)\|_2\leq S_H\|x-y\|_2$.
\end{assumption}
With this additional assumption, Theorem~\ref{thm:VR} shows the (expected) total costs of variance reduction counterparts of MLMC gradient methods, with the proof provided in Appendix~\ref{app_sec:variance_reduction}.

\begin{theorem}[(Expected) Total Costs of Variance Reduction Methods]
\label{thm:VR}
Under Assumptions \ref{assumption:general}\ref{assumption:general:II}\ref{assumption:general:III}, \ref{assumption:gradient_bias}, \ref{assumption:objective}\ref{assumption:objective:I} and \ref{assumption:average_smoothness}, setting $D_1= \cO(\eps^{-2})$, $D_2 = \cO(\eps^{-1})$, $Q_E=\cO(\eps^{-1})$, $\gamma\leq 1/(3S_F)$, to find an $\eps$-stationary point for nonconvex smooth objective $F$, with properly chosen hyper-parameters as in Table~\ref{tab:summary:VR},%
\begin{enumerate}
    \item 
The total cost of VR $L$-SGD is $\cO(\eps^{-3-2c/a}).$
    \item
When $b>c$, the (expected) total costs of VR RT-MLMC and VR V-MLMC are $ \cO(\eps^{-3})$ and $\cO(\eps^{-1-2\max\{1+(c-b)/a, c/a}\})$, respectively; 
when $b=c$, the (expected) total costs are $\tilde \cO(\eps^{-3})$ and $\tilde \cO(\eps^{-1-2\max\{1,c/a\}})$;
and when $b<c$, the (expected) total costs are $\cO(\eps^{-3-2(c-b)/a})$ and $\cO(\eps^{-1-2\max\{1,c/a}\})$.
    \item
When $b>c$, the expected total cost of VR RU-MLMC or VR RR-MLMC is $\cO(\eps^{-3})$.
\end{enumerate}
\begin{table}[!ht]
	\vspace{-1em}
	\caption{
	Hyper-parameters of various gradient estimators in Theorem~\ref{thm:VR}.
	}
	\label{tab:summary:VR}
	\small
	\begin{center}%
{
\setlength\extrarowheight{4pt}
			\begin{tabular}{c|m{10cm}}
			\toprule%
Gradient Estimators & Hyper-parameters \\[4pt]
\hline
\textbf{VR $L$-SGD} & 
 $L = \lceil 1/a \log(4M_a\eps^{-2})  \rceil,\quad 
 n_L=1
 $
\\[4pt]
\hline
\textbf{VR RT-MLMC} & $L = \lceil 1/a \log(4M_a\eps^{-2})  \rceil,\quad q_l = 2^{-(b+c)l/2}(\frac{1-2^{-(b+c)(L+1)/2}}{1-2^{-(b+c)/2}})^{-1}$
\\
\hline
\textbf{VR V-MLMC} & $L = \lceil 1/a \log(4M_a\eps^{-2})  \rceil,\quad n_l = \lceil 2^{-(b+c)l/2}N\rceil,\quad N=\tilde{\cO}(\eps^{-2})$
\\
\hline
$\begin{array}{l}
\textbf{VR RU-MLMC}
\text{ or }\textbf{VR RR-MLMC}
\end{array}
$& $q_l = 2^{-(b+c)l/2}(1-2^{-(b+c)/2})^{-1}$
\\
   \bottomrule
\end{tabular}
}
	\end{center}
 \vspace{-1.2em}
\end{table}
\end{theorem}
It is worth noting that in Theorem~\ref{thm:VR}, we specify $D_1$, $D_2$, $Q_E$, and $\gamma$ following the same setup as in SPIDERBoost~\citep{wang2018spiderboost}.
One may combine our proposed gradient estimators with other variants of variance reduction methods to obtain similar sample complexity results in the nonconvex smooth case.

When $b>c$, intuitively, one may interpret RU-MLMC and RR-MLMC as a way to generate unbiased gradient estimators of $F$ at $\cO(1)$ costs. Thus, VR RU-MLMC and VR RR-MLMC should achieve a total cost of $\cO(\eps^{-3})$, since the total number of unbiased gradient estimators needed by variance reduced method is $\cO(\eps^{-3})$. RT-MLMC is more versatile regarding $a,b,c$, but it introduces bias. As long as the bias can be well controlled, VR RT-MLMC could further reduce the complexity bounds from $\cO(\eps^{-4})$ to $\cO(\eps^{-3})$ when compared to RT-MLMC. 
These observations further imply that MLMC gradient estimators can be easily integrated with other optimization techniques previously developed for stochastic optimization with unbiased oracles. It would greatly enhance the applicability of MLMC gradients in the optimization and machine learning problems. For example, it is an interesting open question to consider Nesterov's accelerated MLMC for certain applications with biased oracles, which we leave for future investigation.

\vspace{-0.5em}
\subsection{Comparison of Different MLMC Gradient Methods}
    When $b>c$, RU-MLMC and RR-MLMC are the most favorable among the four MLMC methods since they have unbiased gradient estimators and do not need to specify $L$ in advance. 
When $b\leq c$, in theory, only RT-MLMC and V-MLMC are applicable as RU-MLMC and RR-MLMC either admit unbounded variance or unbounded expected per-iteration costs. Since RT-MLMC can be treated as imposing a deterministic truncation on RU-MLMC, RT-MLMC, in fact, introduces bias to avoid the high computation cost. It further reflects the importance of balancing bias, variance, and cost tradeoffs.

RT-MLMC is the most versatile algorithm among the four MLMC methods since it has no restrictions on $a,b,c$ and does not need any mini-batch. It suits situations when the per-iteration budget is limited, e.g., when one only has limited samples or computation power. 
V-MLMC requires large mini-batches that lead to a very small variance, and thus a constant stepsize $\cO(1)$ is sufficient to guarantee convergence.

\section{Applications}
\label{section:applications}
In this section, we demonstrate how to perform MLMC gradient methods on CSO for DRO and constrastive learning problems, pricing and staffing in stochastic systems, and UBSR optimization. We further establish the complexity bounds of the MLMC gradient methods by analyzing the $a,b,c$ parameters for these problems. 

\vspace{-0.5em}
\subsection{Conditional Stochastic Optimization}
\label{section:CSO}
Conditional stochastic optimization appears widely in applications from machine learning, including the optimal control in linearly-solvable Markov decision process~\citep{pmlr-v54-dai17a}, policy evaluation and control in reinforcement learning~\citep{pmlr-v54-dai17a,pmlr-v80-dai18c,nachum2020}, meta-learning~\citep{hu2020biased}, instrumental variable regression~\citep{muandet2019dual}. Previously, \citet{hu2020biased} considered biased SGD and its variance reduction counterpart using SPIDER~\citep{fang2018spider}. 
These methods are special cases of $L$-SGD or VR $L$-SGD proposed in this manuscript.

In the following, we show that MLMC methods can significantly reduce the sample complexity and achieve better results than previous biased variance reduction methods proposed in \citep{hu2020biased}. The CSO problem and its approximation function $F^l(x)$ admit the following form:
\begin{equation}
\label{pro:CSO}
\min_{x\in\RR^d} F(x):=\EE_\xi f_\xi(\EE_{\eta|\xi} g_\eta(x,\xi)), \quad\ F^l(x)= \EE_{\xi^l} \EE_{\{\eta_j^l\}_{j=1}^{2^l}|\xi^l}\Big[ f_{\xi^l}\Big(\frac{1}{2^l}\sum_{j=1}^{2^l} g_{\eta_j^l}(x,\xi^l)\Big)\Big],
\end{equation}
where $\xi^l\sim \PP(\xi)$ and $\{\eta_{j}^l\}_{j=1}^{2^l}\sim\PP(\eta|\xi_i^l)$.
Denote the collection of random samples as $\zeta^l = (\xi^l,\{\eta_{j}^l\}_{j=1}^{2^l})$ and take $\hat g_{n_1:n_2}(x,\zeta^l) = (n_2-n_1+1)^{-1} \sum_{j=n_1}^{n_2} g_{\eta_{j}^l}(x,\xi^l)$ for some $1\le n_1\le n_2$.
For each query point $x$, $\mathcal{SO}^l$ returns $(h^l(x,\zeta^l), H^l(x,\zeta^l))$ such that $h^l(x,\zeta^l) 
     = 
     \nabla_xf_{\xi^l}(\hat g_{1:{2^l}}(x,\zeta^l))$, and
\begin{equation}
\label{eq:cso_oracle_return}
\begin{split}
     H^l(x,\zeta^l) 
     = 
     \nabla_x\Big[ f_{\xi^l}(\hat g_{1:{2^l}}(x,\zeta^l))- 
     \frac{1}{2}[f_{\xi^l}(\hat g_{1:{2^{l-1}}}(x,\zeta^l)) + f_{\xi^l}(\hat g_{1+{2^{l-1}}:{2^l}}(x,\zeta^l))]
     \Big].
\end{split}
\end{equation}
Under appropriate conditions, it can be shown that Assumptions~\ref{assumption:general} and \ref{assumption:gradient_bias} hold with the parameter $a=b=c=1$ (see Proposition~\ref{prop:cso} in Appendix~\ref{Appendix:applications:1}).
Since $b=c$, RU-MLMC and RR-MLMC are not applicable. 
Leveraging previous results on MLMC gradient methods (i..e, Theorems~\ref{thm:RT-MLMC}, \ref{thm:V-MLMC}, and \ref{thm:VR}), 
we derive the sample complexity of V-MLMC and RT-MLMC for solving CSO and summarize it in Theorem~\ref{cor:cso}.

\begin{theorem}[Sample Complexity for Solving CSO]
\label{cor:cso}
Under mild assumptions~(see Assumption~\ref{assumption:cso} in Appendix~\ref{Appendix:applications:1}), for strongly convex $F(x)$, the sample complexity of V-MLMC and RT-MLMC for finding $\eps$-optimal solution is $\tilde \cO(\eps^{-1})$; for convex $F(x)$, the sample complexity of V-MLMC and RT-MLMC is $\tilde \cO(\eps^{-2})$; for nonconvex smooth $F(x)$, the sample complexity of V-MLMC and RT-MLMC for finding $\eps$-stationary point is $\tilde \cO(\eps^{-4})$, which reduces to $\tilde \cO(\eps^{-3})$ when using VR V-MLMC and VR RT-MLMC.
\end{theorem}
The BSGD in \citet{hu2020biased} achieved $\tilde \cO(\eps^{-2})$, $\cO(\eps^{-3})$, and $\cO(\eps^{-6})$ sample complexity in the strongly convex, convex and nonconvex smooth setting, respectively. Since BSGD is a special case of the $L$-SGD framework in our paper, we directly recover the corresponding sample complexity when plugging $a=b=c=1$ into Theorem \ref{thm:L-SGD}. \citet{hu2020biased} further investigated variance reduction techniques in the nonconvex smooth setting and achieved $\cO(\eps^{-5})$ sample complexity, 
albeit inferior to the sample complexity of  V-MLMC and RT-MLMC methods and the corresponding variance-reduced counterparts. In Appendix \ref{Appendix:CSO:special_case}, we further show that when $g_\eta$ is linear in $x$, we have $a=c=1$ and $b=2$. In such cases, RU-MLMC and RR-MLMC are applicable and have theoretical guarantees. Table \ref{tab:CSO} summarizes the complexity bounds of CSO problems.
\begin{remark}[Comparison with Lower Bounds for CSO in \citet{hu2020biased}]
Notably, the lower bounds on sample complexities in~\citet{hu2020biased} do not apply to our results.
The reason is that RT-MLMC and V-MLMC use a low-cost white-box oracle rather than a high-cost black-box oracle defined in the aforementioned reference.
It further highlights the importance of considering the oracle cost when studying the complexity of biased oracle models. 
\end{remark}

\begin{table}[!ht]
    \vspace{-1em}
    {
    \footnotesize
    \caption{Comparison of Algorithms on CSO}
    \renewcommand\arraystretch{1}
    \begin{center}
    \begin{tabular}{c|ccc|c}
        \toprule[1pt]
        {Condition on $F$}
        & {V-MLMC}
        & RT-MLMC
        & RU-MLMC/RR-MLMC *
        & BSGD~($L$-SGD)~\citep{hu2020biased}
        \\
        \hline

        {Strongly Convex}
        & $\tilde \cO(\eps^{-1})$
        & $\tilde  \cO(\eps^{-1})$
        & $ \cO(\eps^{-1})$
        & $\tilde  \cO(\eps^{-2})$
        \\
        \hline

        {Convex}
        & $\tilde \cO(\eps^{-2})$
        & $\tilde \cO(\eps^{-2})$
        & $\cO(\eps^{-2})$
        & $\cO(\eps^{-3})$
        \\
        \hline

        {Nonconvex}
        & $\tilde \cO(\eps^{-4})$
        & $\tilde \cO(\eps^{-4})$
        & $\cO(\eps^{-4})$
        & $\cO(\eps^{-6})$
        \\
        \hline
        
        {Individual Smoothness}
        & $\tilde \cO(\eps^{-3})$
        & $\tilde \cO(\eps^{-3})$
        & $\cO(\eps^{-3})$
        & $\cO(\eps^{-5})$
        \\
        \hline
        
        \multicolumn{5}{c}{* RU-MLMC and RR-MLMC are only applicable under the special case when $g_\eta(x,\xi)$ is linear (Assumption \ref{assumption:cso_additional}).}
        \\

    \bottomrule[1pt]
    \end{tabular}  
    \label{tab:CSO}
    \end{center}
    }
    \vspace{-1em}
    \end{table}

\subsubsection{Applications of CSO in DRO.}
DRO has shown promise to address the challenge of decision-making under uncertainty, especially in high-stake environments.
Given a reference distribution $\mathbb{P}^0$ of $\xi$ that usually takes the form of an empirical distribution, DRO proposes to solve
\begin{equation}
    \label{Eq:DRO}\min_{x}\max_{\mathbb{P}:~\mathcal{D}(\mathbb{P}, \mathbb{P}^0)\le \rho}~\Big\{ 
\mathbb{E}_{\xi\sim\mathbb{P}}[\ell(x;\xi)]
    \Big\},
\end{equation}
where $l(x;\xi)$ is the loss function dependent on the decision $x$ and the random variable $\xi$, and the constraint $\mathcal{D}(\mathbb{P}, \mathbb{P}^0)\le \rho$ allows to take into account all distributions near $\mathbb{P}^0$ up to $\rho$-radius ball using a pre-specified divergence $\mathcal{D}(\cdot,\cdot)$ (such as $\phi$-divergence and regularized Wasserstein distance).
\citet{levy2020large} utilized MLMC gradient methods to solve $\phi$-divergence DRO.
When specifying the $\phi$-divergence as the Kullback-Leibler~(KL) divergence, the corresponding DRO is a special CSO~\eqref{pro:CSO} with $f_\xi(x)$ as a log function and the conditional distribution $\PP(\eta\mid\xi)$ as the reference distribution~\citep{qi2022stochastic}.
In the context of the KL-divergence regularized $2$-Wasserstein DRO model (also known as the Sinkhorn DRO model, see references ~\cite{wang2021sinkhorn, blanchet2020semi, azizian2023regularization}), by taking the strong dual reformulation of \eqref{Eq:DRO} and fixing the Lagrangian multiplier associated with the divergence ball constraint, the goal is to solve the problem
\begin{equation}
\min_{x}~\mathbb{E}_{\xi}\Big[ 
\lambda\log\mathbb{E}_{\eta\mid\xi}\left[ 
e^{\ell(x; \eta)/\lambda}
\right]
\Big],\label{Eq:SDRO}
\end{equation}
where $\xi\in\mathbb{R}^{d_{\xi}}$ follows the reference distribution $\PP^0(\xi)$, $\eta\mid\xi$ follows the conditional distribution $\mathcal{N}(\xi, \tau^2\mathbf{I}_{d_{\xi}})$, and $\lambda, \tau^2>0$ are problem coefficients.
Here, $\ell(x; \eta)$ denotes the loss function dependent on the decision variable $x$ and random data $\eta$.
This formulation can also be viewed as a special CSO problem~\eqref{pro:CSO} with $f_{\xi}(\cdot)=\lambda\log(\cdot)$ and $g_{\eta}(\cdot, \xi)=e^{\ell(\cdot;\eta) / \lambda}$.
Following the assumption in \cite[Assumption~2]{wang2021sinkhorn} that $\ell$ is $B$-uniformly bounded, 
$L_{\ell}$-Lipschitz continuous, 
and $S_{\ell}$-smooth,
one can see that
Assumption~\ref{assumption:cso} holds with $\sigma^2=B^2, S_f=L_{\ell}=\lambda, L_g=e^{B/\lambda}L_{\ell}/\lambda, S_g=e^{B/\lambda}(S_{\ell}+L_{\ell}^2)/\lambda$.
Therefore, the sample complexity results of MLMC gradient methods in Theorem~\ref{cor:cso} can directly apply to solving KL regularized 2-Wasserstein DRO, i.e., the Sinkhorn DRO.
We present a comprehensive numerical study in Section~\ref{Sec:syn} to validate the performance of various gradient methods for solving this problem.

\subsubsection{Applications of CSO in CL.}\label{Sec:CSO:con:learn}
CL is a popular unsupervised learning technique, aiming to learn representations such that similar instances are close while dissimilar instances are far apart in the representation space~\citep{chopra2005learning,oord2018representation,wang2020understanding}.
In detail, let $\cv x$ be an input batch of $N$ samples. We augment these samples and treat two augmented samples from the same sample in $\cv x$ as similar and otherwise dissimilar.
For $k=1,\ldots,N$, let $\tilde{\cv x}[2k-1]$ and $\tilde{\cv x}[2k]$ denote two augmented samples from the $k$-th sample from $\cv x$.
Let $\cv z[i], i=1,\ldots,2N$ denote the neural network outputs of $\tilde{\cv x}[i], i=1,\ldots,2N$.
CL aims to train the neural network model such that 
(i) the (cosine) similarity between distinct samples $i,j$, denoted as $s_{i,j}\triangleq \frac{\cv z[i]^{\top}\cv z[j]}{\|\cv z[i]\|\|\cv z[j]\|}$, is minimized; and 
(ii) the similarity between different data augmentations $\tilde{\cv x}[2k-1], \tilde{\cv x}[2k]$ of the same sample $k$, denoted as $s_{2k-1,2k}\triangleq \frac{\cv z[2k-1]^{\top}\cv z[2k]}{\|\cv z[2k-1]\|\|\cv z[2k]\|}$, is maximized.
The well-known simple framework for contrastive learning of visual representations~(SimCLR)~\citep{chen2020simple} solves the CL task based on
\[
\begin{aligned}
\min_{\theta}&~\mathcal{L}_{\mathrm{SimCLR}}(\theta)\triangleq\mathbb{E}_{\cv x}\left[ 
\frac{\sum_{k\in[N]}\ell_{2k-1,2k} + \ell_{2k,2k-1}}{2N}
\right],\quad
\mbox{ where }\ell_{i,j}=-\log\frac{\exp(s_{i,j}/\tau)}{
\sum_{k\in[2N]}\textsf{1}_{k\ne i}\exp(s_{i,k}/\tau)
}.
\end{aligned}
\]
Here the variable $\theta$ corresponds to the weight of the neural network for the data-to-representation mapping.
It is worth noting that the SimCLR objective is a special case of CSO~\citep{Qiu23}:
\begin{align*}
\mathcal{L}_{\mathrm{SimCLR}}(\theta)
&=
\mathbb{E}_{\cv x, i\sim\mathrm{Uni}([N])}\Big[ 
f\Big( 
\mathbb{E}_kg_{k}(\theta; i)
\Big)
-
s_{2i,2i-1}/\tau
\Big] 
+ \mathrm{Constant},
\end{align*}
where the index $i$ follows uniform distribution over $[N]$ (denoted as $i\sim \mathrm{Uni}([N])$), 
the outer loss function $f:~\mathbb{R}_+^2\to\mathbb{R}$ is defined as $f(y_1,y_2) = \frac{1}{2}\log(y_1y_2)$, 
the inner loss function 
$g_k(\theta,i) = (e^{s_{2i-1,k_1}}/\tau, e^{s_{2i, k_2}}/\tau)$,
and the random vector $k=(k_1,k_2)$ satisfies $k_1\sim \mathrm{Uni}([2N]\setminus\{2i-1\}), k_2\sim \mathrm{Uni}([2N]\setminus\{2i\})$.

\citet{chen2020simple} has highlighted the effectiveness of employing a large batch size $N$ for achieving satisfactory performance in CL.
Nevertheless, a notable drawback of this approach lies in the escalated computational and memory costs. 
In this context, the MLMC gradient methods can substantially alleviate these computational and memory burdens, whereas the $L$-SGD corresponds to the vanilla training procedure in the existing literature.
We will demonstrate the empirical superior performance of MLMC methods over the baseline in Section~\ref{Sec:con:learning} for a moderately large batch size.

\vspace{-0.5em}
\subsection{Pricing and Staffing in Stochastic Systems}\label{Sec:pricing:staffing}
Consider the joint pricing and staffing task for a stochastic system~\citep{lee2014optimal}, where the service provider's goal is to 
seek the optimal service fee and service capacity such that the long-run expected profit (the service revenue minus the staffing cost and customer delay penalty) is maximized, leading to the following stochastic optimization problem:
\begin{equation}
\begin{aligned}
\min_{\mu, p}&~\Big\{F(\mu,p):= h_0\EE[Q_{\infty}(\mu,p)] + c(\mu) - p\lambda(p)\Big\}.\\
\end{aligned}
\label{Eq:pricing:capacity:formula}
\end{equation}
In the problem above, $Q_{\infty}(\mu,p)$ denotes the steady-state queue length of customers under the price $p$ and capacity $\mu$, and $h_0$ denotes the holding cost coefficient; $c(\mu)$ denotes the cost incurred from providing service capacity $\mu$; $\lambda(p)$ denotes the rate of demand when the customers are charged with service fee $p$.

We consider the stochastic system to be a GI/GI/1 queue, i.e., 
the interarrival time of customers and the service time both follow a general i.i.d. distribution.
Unfortunately, the gradient of the objective function cannot be numerically simulated exactly, particularly for the gradient estimator of $\EE[Q_{\infty}(\mu,p)]$.
By~\cite[Lemma~5]{chen2023online}, the gradient of the objective is
\begin{align*}
\nabla F(\mu,p)=\begin{pmatrix}
\frac{\partial F(\mu,p)}{\partial p}\\
\frac{\partial F(\mu,p)}{\partial \mu}
\end{pmatrix}=
\begin{pmatrix}
    -\lambda(p)-p\lambda'(p)+
h_0\lambda'(p)\Big( 
\EE[W_{\infty}(\mu,p)+X_{\infty}(\mu,p)]
+\frac{1}{\mu}
\Big)\\
c'(\mu)
-h_0\frac{\lambda(p)}{\mu}\Big( 
\EE[W_{\infty}(\mu,p)+X_{\infty}(\mu,p)]
+\frac{1}{\mu}
\Big)
\end{pmatrix},
\end{align*}
where $W_{\infty}(\mu,p)$ and $X_{\infty}(\mu,p)$ denotes the steady-state customer waiting time and service's busy time under price $p$ and capacity $\mu$, respectively.
Note that $W_{\infty}(\mu,p)$ and $X_{\infty}(\mu,p)$ are the limits of processes $\{W_n\}_n$ and $\{X_n\}_n$, which can be simulated using the recursive relation~\citep{lindley1952theory}
\[
 W_n(\mu,p)=\Big(W_{n-1}(\mu,p) + \frac{V_n}{\mu} - \frac{U_n}{\lambda(p)}\Big)^+,\quad 
 X_n(\mu,p)=\textbf{1}_{\{W_n(\mu,p)>0\}}\cdot \Big( 
X_{n-1}(\mu,p) + \frac{U_n}{\lambda(p)}
\Big),
\]
with $W_0(\mu,p)=X_0(\mu,p)=0$, and $\{U_n\}_n, \{V_n\}_n$ being i.i.d. random variables of (normalized) customer interarrival time and service time, respectively. They are normalized to ensure $\mathbb{E}[U_n]=\mathbb{E}[V_n]=1$.

We simulate the gradient in the following way.
Denote $x=(\mu,p)$, the collection of random parameters $\zeta^l_x = \{W_n(x), X_n(x)\}_{n=1}^{2^l}$ for $l\in\mathbb{N}$ and the function $\hat g_{-n_1:n_2}(x,\zeta^l) = \frac{1}{n_1} \sum_{j=n_2-n_1+1}^{n_2}(W_j(x) + X_j(x))$ for some $1\leq n_1\leq n_2$. 
We construct the approximate gradient with negligible bias as
\begin{equation}\label{Eq:approx:Fl:ss}
\nabla F^l(x)=
\begin{bmatrix}
    -\lambda(p)-p\lambda'(p)+
h_0\lambda'(p)\Big( 
\EE[
\hat g_{-m:2^l}(x,\zeta^l)
]
+\frac{1}{\mu}
\Big)
&~~
c'(\mu)
-h_0\frac{\lambda(p)}{\mu}\Big( 
\EE[\hat g_{-m:2^l}(x,\zeta^l)]
+\frac{1}{\mu}
\Big)
\end{bmatrix}^{\mathrm{T}},
\end{equation}
where $m\in\mathbb{N}_+$ denotes the mini-batch size of samples. Note that this gradient does not have a corresponding function $F^l$.
For each query on a query point $x$, $\mathcal{SO}^l$ returns $(h^l(x,\zeta^l), H^l(x,\zeta^l))$ such that
\begin{equation}
\label{eq:pricing}
\begin{aligned}
    &  h^l(x,\zeta^l) 
     = 
    \begin{bmatrix}
    -\lambda(p)-p\lambda'(p)+
h_0\lambda'(p)\Big( 
\hat{g}_{-m:2^l}(x, \zeta^l)
+\frac{1}{\mu}
\Big)\qquad&\quad
c'(\mu)
-h_0\frac{\lambda(p)}{\mu}\Big( 
\hat{g}_{-m:2^l}(x, \zeta^l)
+\frac{1}{\mu}
\Big)
\end{bmatrix}^{\mathrm{T}},\\
 &    H^l(x,\zeta^l) 
     = 
     \begin{bmatrix}
h_0\lambda'(p)
\big(\hat{g}_{-m:2^l}(x, \zeta^l) - \hat{g}_{-m:2^{l-1}}(x, \zeta^l)\big)
\qquad&\quad
-h_0\frac{\lambda(p)}{\mu}\big(\hat{g}_{-m:2^l}(x, \zeta^l) - \hat{g}_{-m:2^{l-1}}(x, \zeta^l)\big)
\end{bmatrix}^{\mathrm{T}}.
\end{aligned}
\end{equation}
Using this oracle, one can efficiently simulate the unbiased estimator of \eqref{Eq:approx:Fl:ss}.
According to the technical assumption in \citep[Section~3.1]{chen2023online} and Proposition~\ref{Prop:pricing}, the bias and variance of $H^l(x,\zeta^l)$ decrease in exponential of exponential rate. In this case, the sample complexity of $L$-SGD and RT-MLMC methods will be the same, i.e., they achieve complexity $\tilde{\cO}(\eps^{-1})$ to find $\eps$-optimal solution for PL objective and $\tilde{\cO}(\eps^{-4})$ to find $\eps$-stationary point for Lipschitz continuous and smooth objective. 
However, the exponential of exponential decay rate in the bias and variance, $\cO(e^{-0.25 \mathfrak{c} 2^l})$, has an exponent $\mathfrak{c}$ that is very small in practice. In such a case, $e^{-0.25 \mathfrak{c} 2^l}\approx\cO(2^{-l})$ for a moderate $l$. It motivates us to consider the more conservative polynomial decay rate of the bias and variance for practical applications, under which this oracle satisfies Assumptions~\ref{assumption:general}\ref{assumption:general:II}\ref{assumption:general:III} and \ref{assumption:gradient_bias} with $a=b=c=1$.
From the numerical results in Section~\ref{Sec:pricing:capacity},
we observe that the algorithms can successfully find the global optimal solution despite nonconvexity.
This fact motivates us to present the complexity analysis for the case where the smooth $F(\cdot)$ satisfies PL condition in the following Theorem~\ref{result:pricing:capacity:formula}. See Remark \ref{remark:PL} and \citet{karimi2016linear} for further details about PL condition.
\begin{theorem}[Sample Complexity for Pricing and Staffing]
\label{result:pricing:capacity:formula}
Under appropriate conditions (see Assumption~\ref{Assumption:pricing} in Appendix~\ref{Appendix:applications:2}), for smooth $F(\cdot)$ satisfying PL condition, the sample complexity of RT-MLMC for finding $\eps$-optimal solution is $\tilde{\cO}(\eps^{-1})$, whereas that of $L$-SGD is ${\cO}(\eps^{-2})$; 
for smooth and Lipschitz continuous $F(\cdot)$, 
the sample complexity of RT-MLMC for finding $\eps$-stationary point is $\tilde{\cO}(\eps^{-4})$, 
whereas that of $L$-SGD is ${\cO}(\eps^{-6})$.
\end{theorem}
It is noteworthy that our previous analysis on MLMC methods (i.e., Theorems~\ref{thm:L-SGD} and \ref{thm:RT-MLMC}) assumes that the approximate objective $F^l$ satisfies the smoothness condition, which does not hold in this application. 
In contrast, Theorem~\ref{result:pricing:capacity:formula} imposes technical assumptions on the ground truth objective $F$ only.
By adopting a variant of Theorems~\ref{thm:L-SGD} and \ref{thm:RT-MLMC}, we still obtain the sample complexity of $L$-SGD and RT-MLMC methods.
We numerically validate their performance in Section~\ref{Sec:pricing:capacity} and show that MLMC gradient method indeed achieves a significant acceleration. We leave verifying the PL condition for future studies as it goes beyond the scope of the  current work.

\subsection{Utility-Based Shortfall Risk~(UBSR) Optimization}
\label{Sec:UBSR}
UBSR has gained wide applications in financial engineering~\citep{delage2022shortfall, hu2016convex, giesecke2008measuring, follmer2002convex}, aiming to balance the trade-off of controlling the risk and maximizing the return.
\citet{hegde2021online} proposed to find the optimal decision in terms of UBSR:%
\begin{equation}
\tag*{(UBSR-O)}
\begin{aligned}
\min_{\theta}&\quad \Big\{\mathcal{SR}_{\lambda}(\theta):=
\inf_{t:~\mathbb{E}[\ell(-X(\theta)-t)]\le\lambda}~t\Big\},
\end{aligned}
\label{Eq:expression}
\end{equation}
where $\lambda$ is a pre-specified risk level, $\ell(\cdot)$ is a given convex loss, and $X(\theta)$ is a random variable that is dependent on the decision variable $\theta$.
By~\citep[Section~5.1]{hegde2021online}, the gradient of UBSR equals
\begin{equation}
\nabla \mathcal{SR}_{\lambda}(\theta)
=
-\frac{\mathbb{E}_{X(\theta)}\big[\ell'\big(-X(\theta)
-\mathcal{SR}_{\lambda}(\theta)
\big)X'(\theta)\big]}{\mathbb{E}_{X(\theta)}\big[\ell'\big(-X(\theta)
-\mathcal{SR}_{\lambda}(\theta)
\big)\big]}.
\label{Eq:grad:UBSR}
\end{equation}
In practice, the distribution of $X(\theta)$ is not given in explicit form, whereas i.i.d. samples of $X(\theta)$ (and correspondingly $X'(\theta)$) are available.
Although obtaining the unbiased stochastic estimator of \eqref{Eq:grad:UBSR} is challenging, one can apply MLMC methods to simulate the gradient estimator with negligible bias efficiently.

Denote the random parameter $\zeta^l(\theta)=(\{X_i(\theta)\}_{i=1}^{2^l},
\{\tilde{X}_i(\theta)\}_{i=1}^{2^l})$, where $X_i(\theta), \tilde{X}_i(\theta)$ are i.i.d. copies of $X(\theta)$.
Let $t_{l}(\{\tilde{X}_i(\theta)\}_{i=1}^{2^l})$ be a stochastic estimator of $\mathcal{SR}_{\lambda}(\theta)$ adopted from \citep[Section~4]{hegde2021online} using samples $\{\tilde{X}_i(\theta)\}_{i=1}^{2^l}$.
As have shown in \citep[Theorem~1]{hegde2021online}, this estimator has low bias $\mathcal{O}(2^{-l})$.
Next, define the approximation gradient %
\begin{equation}
\nabla\mathcal{SR}_{\lambda}^l(\theta) = -\frac{\EE_{\zeta^l(\theta)}\big[
\frac{1}{2^l}\sum_{i=1}^{2^l}
\ell'\big(-X_i(\theta)
-t_{l}(\{\tilde{X}_i(\theta)\}_{i=1}^{2^l})
\big)X_i'(\theta)
\big]}{\EE_{\zeta^l(\theta)}\big[
\frac{1}{2^l}\sum_{i=1}^{2^l}
\ell'\big(-X_i(\theta)
-t_{l}(\{\tilde{X}_i(\theta)\}_{i=1}^{2^l})
\big)
\big]}.
\label{Eq:exp}
\end{equation}
One can apply MLMC methods to obtain sample estimates of $\nabla\mathcal{SR}_{\lambda}^l(\theta)$. %
Define 
\[
\hat g_{n_1:n_2}(\theta,\zeta^l) = -\frac{
(n_2-n_1+1)^{-1} \sum_{j=n_1}^{n_2} \ell'\big(-X_j(\theta)
-t_{l}(\{\tilde{X}_i(\theta)\}_{i=n_1}^{n_2})
\big)X_j'(\theta)
}{
(n_2-n_1+1)^{-1} \sum_{j=n_1}^{n_2} \ell'\big(-X_j(\theta)
-t_{l}(\{\tilde{X}_i(\theta)\}_{i=n_1}^{n_2})
\big)
}.
\]
For each query on point $\theta$, we construct the oracle $\mathcal{SO}^l$ to return 
$(h^l(\theta,\zeta^l), H^l(\theta,\zeta^l))$ with
\begin{equation}
\label{eq:UBSR}
\begin{split}
    &  h^l(\theta,\zeta^l) 
     = 
    \hat g_{1:2^l}(\theta,\zeta^l),
     \quad 
 H^l(\theta,\zeta^l) 
     = 
      \hat g_{1:2^l}(\theta,\zeta^l) - \frac{1}{2}[\hat g_{1:2^{l-1}}(\theta,\zeta^l) +  \hat g_{2^{l-1}+1:2^l}(\theta,\zeta^l)].
\end{split}
\end{equation}
By Proposition~\ref{prop:UBSR} in Appendix~\ref{Appendix:applications:3}, this oracle satisfies Assumptions~\ref{assumption:general}\ref{assumption:general:II}\ref{assumption:general:III} and \ref{assumption:gradient_bias} with $a=b=c=1$.
Following a similar proof procedure as in Theorem~\ref{result:pricing:capacity:formula}, we obtain the sample complexity results of V-MLMC and RT-MLMC for solving \ref{Eq:expression}.
\begin{theorem}[Sample Complexity for UBSR-O]
\label{cor:UBSR_O}
Under technical conditions in \cite[Assumptions~1 to 9]{hegde2021online}, for smooth $\mathcal{SR}_{\lambda}(\theta)$ satisfying PL condition, the sample complexity of V-MLMC and RT-MLMC for finding $\eps$-optimal solution is $\tilde{\cO}(\eps^{-1})$; for smooth and Lipschitz continuous $\mathcal{SR}_{\lambda}(\theta)$, the sample complexity of V-MLMC and RT-MLMC for finding $\eps$-stationary point is $\tilde{\cO}(\eps^{-4})$.%
\end{theorem}
The state-of-the-art algorithm in \cite{hegde2021online} is a special case of $L$-SGD in our framework, which has a complexity $\tilde{\cO}(\eps^{-2})$ in the strongly convex case (a special case of the PL case). 
We improve the rate for this case to $\tilde{\cO}(\eps^{-1})$.
For the nonconvex case, we present the first sample complexity results, where nonconvex case is common when adapting a complicated loss function.

\vspace{-0.5em}
\subsection{Comparison of MLMC Gradient Methods against MLMC on Function Value or Solution
}\label{Sec:statistical:inference}

In this subsection, we discuss the benefits of MLMC gradient methods over MLMC function value methods and solution methods in optimization.  We use the  CSO problem as an example but the implications hold for general stochastic optimization with biased oracles. 

\noindent\textbf{MLMC Function Value Estimator.}
Consider the CSO problem in the form of \eqref{pro:CSO}. 
A natural idea to use the nested \emph{sample average approximation}~(SAA) and solve the empirical counterpart
\[
\hat{F}_{n,m}(x) := \frac{1}{n}\sum_{i=1}^nf_{\xi_i}\Big( 
\frac{1}{m}\sum_{j=1}^m g_{\eta_{i,j}}(x, \xi_i)
\Big).
\]
Similar to the $L$-SGD methods, such a nested SAA methods suffers from high sample complexity~\citep{hu2020sample}.
We show how to use MLMC to construct a function value estimator with a better sample complexity.

For a maximum level $L\in\mathbb{N}_+$ and mini-batch sizes $m_l = 2^{l-1}$, 
at each level $l=0,\ldots,L$, generate $n_l$ i.i.d. samples $\{\xi_i^l\}_{i=1}^{n_l}$ from $\PP(\xi)$, and conditioned on each $\xi_i^l$, generate $m_l$ i.i.d. samples $\{\eta_{ij}^l\}_{j=1}^{m_l}$ from $\PP(\eta|\xi_i^l)$. 
Denote the random parameters as $\zeta^l_i = (\xi_i^l,\{\eta_{ij}^l\}_{j=1}^{m_l})$.
Finally, construct MLMC function value estimator
\begin{equation*}
\hat F(x) = \sum_{l=0}^L \frac{1}{n_l}\sum_{i=1}^{n_l}V^l(x, \zeta^l_i),
\text{ }V^l(x, \zeta^l_i) := 
f_{\xi_i^l}(\hat{g}_{1:m_l}(x, \zeta^l_i) )- \frac{
f_{\xi_i^l}(\hat{g}_{1:m_{l-1}}(x, \zeta^l_i))
+
f_{\xi_i^l}(\hat{g}_{m_{l-1}+1:m_l}(x, \zeta^l_i))
}{2}.\label{Eq:MLMC:estimator}
\end{equation*}

Unfortunately, even if the population objective and the empirical objective in \eqref{pro:CSO} are convex, the MLMC function value estimator is not necessarily convex since the estimator $V^l$ involves the difference of convex functions.
As a result, the empirical MLMC function value estimator cannot readily be solved to global optimality.  
In fact, the MLMC function value estimator is more suitable when estimating the true objective for a given $x$, i.e., in statistical inference settings, instead of in optimization. For completeness, we establish the improved sample complexity to ensure the uniform convergence between $\hat{F}$ and $F$ under the same assumptions in \citet{hu2020sample} for CSO problem. 
\begin{theorem}\label{Theorem:uniform:converge}
Suppose Assumption~\ref{Assume:CSO:special} holds, 
then for any $\bar{\epsilon}>0$ and $\alpha\in(0,1)$, there exists $\eps_1>0$ such that for any $\eps\in(0, \eps_1)$, the total cost to ensure $\PP
\left( 
\sup_{x\in\cX}~|F(x) - \hat{F}(x)|\ge \eps
\right)\le \alpha$ is $\tilde{\cO}(d\eps^{-2})$.
\end{theorem}
One may also refer to \citep{blanchet2019unbiased} for discussions on RU-MLMC function value estimator.

\subsubsection*{MLMC Solution Estimator.} To overcome the nonconvexity issue brought in by the MLMC function value estimator, one could use V-MLMC or RU-MLMC to build solution estimators~\citep{frikha2016multi,blanchet2019unbiased}. We demonstrate the idea via a novel RT-MLMC solution estimator.  Assuming that the empirical problem $\hat F_{n,m}(x)$ is strongly convex, we repeat the following procedure for $K$ independent trials: for $k=1,...,K$, solve the following  to obtain $\hat x(k)$. 
\begin{align*}
    & \hat x^l := \underset{x\in\mathcal{X}}{\argmin} \frac{1}{n}\sum_{i=1}^{n} f_{\xi_i}\big(\hat{g}_{1:m_l}(x, \zeta^l_i)\big),\quad
     \hat x^l_a:=\underset{x\in\mathcal{X}}{\argmin} \frac{1}{n}\sum_{i=1}^{n} f_{\xi_i}\big(\hat{g}_{1:{m_{l-1}}}(x, \zeta^l_i)\big),\\
    &\hat x^l_b:=\underset{x\in\mathcal{X}}{\argmin} \frac{1}{n}\sum_{i=1}^{n} f_{\xi_i}\big(\hat{g}_{1+{m_{l-1}}:{m_{l}}}(x, \zeta^l_i)\big),\quad
     \hat x(k) :=\frac{1}{p_l} \Big(\hat x^l -\frac{1}{2} (\hat x^l_a + \hat x^l_b)\Big).
\end{align*}
The MLMC solution estimator is $\hat x := \frac{1}{K}\sum_{k=1}^K \hat x(k)$. 
Such an estimator requires solving multiple empirical problems to optimality, while MLMC gradient methods only solve once. Secondly, it requires these nested SAA problems to be strongly/strictly convex to admit a unique solution. Recall that the essence of the MLMC method is to ensure a variance reduction effect, i.e., $\Var(\hat x^l -\frac{1}{2} (\hat x^l_a + \hat x^l_b))$ should be small. Unfortunately, for convex problems that admit multiple optimal solutions, the variance reduction effect is not easily achievable. For nonconvex problems, it is hard to find optimal solutions and multiple stationary points could exist.

In summary, both MLMC function value and solution estimator are not suitable for optimization due to some limitation. Instead, the studied MLMC gradient method do not have such limitations and can be easily combined with existing first-order algorithms and are effective under various structural assumptions.

\section{Numerical Experiments}
\label{section:numerical_results}

In this section, we test our MLMC gradient methods in three examples: DRO, CL, and pricing and capacity sizing for stochastic systems.
Experiments are conducted on Google Colab with an Nvidia Tesla V100-SXM2-16GB GPU.
\vspace{-0.5em}
\subsection{Synthetic Problem with Biased Oracles}\label{Sec:syn}
Consider the Sinkhorn DRO problem in \eqref{Eq:SDRO} with parameters $(\sigma^2,\lambda)=(0.1,20)$.  %
For the formulation therein, we specify $\xi:=(a,b)$ to follow an empirical distribution from feature-label data $\{(a_i,b_i)\}_{i=1}^n$; 
the feature of the data $\eta$~(denoted as $z$) to follow $\mathbf{N}(a, \sigma^2\mathbb{I}_d)$ whereas its label equal $b$;
and %
$\ell(x;\eta)=(f_x(z) - b)^2$.

We evaluate the performance of our proposed optimization algorithms for this problem in both convex and nonconvex setups.
In the convex setup, we take data points $\{(a_i,b_i)\}_{i=1}^n$ from three real-world LIBSVM~\citep{libsvmregression} datasets~(housing, mg, and mpg), and $f_x(z)$ as a linear predictor;
in the nonconvex setup, we take data points $\{(a_i,b_i)\}_{i=1}^n$ from two additional real-world UCI~\citep{dua2017} datasets that are more intricate in terms of sample size and data dimension~(Gas and bike), and $f_x(z)$ as a neural network.

\begin{figure}[!ht]
\vspace{-1em}
    \centering
        \includegraphics[width=0.8\textwidth]{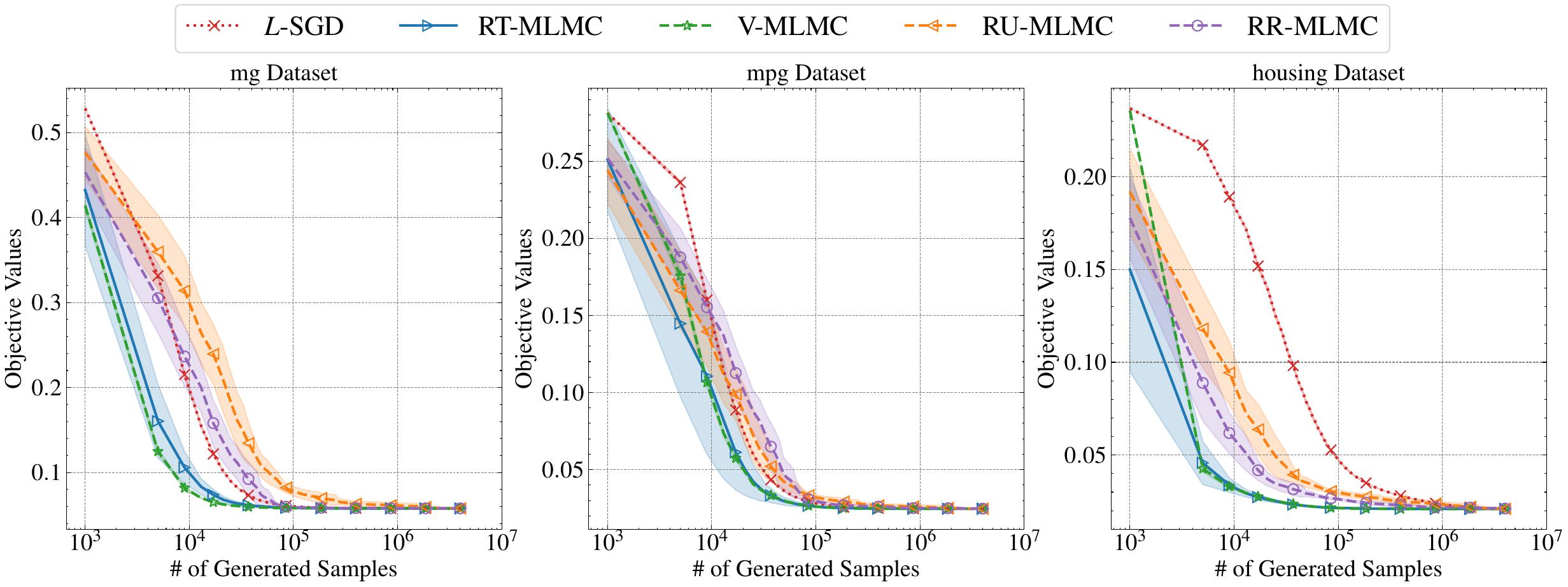}
    \caption{Comparison results of $L$-SGD, V-MLMC, RT-MLMC, RU-MLMC, and RR-MLMC methods when $f_x(z)$ is convex. The $x$-axes represent the number of generated samples, and $y$-axes represent objective values. The results are averaged with error bars based on $10$ independent runs.}
    \label{fig:new:revision}
    \vspace{-2em}
\end{figure}

Figure~\ref{fig:new:revision} reports the performance of various gradient algorithms in the convex setting. 
The $x$-axis corresponds to the number of generated samples, and the $y$-axis corresponds to the objective value.
The plot clearly shows that V-MLMC and RT-MLMC algorithms consistently have faster convergence rates than other algorithms, which aligns with the fact that their theoretical sample complexity rate is optimal.
In contrast, the RR-MLMC and RU-MLMC algorithms have large variances throughout the training process, as suggested by the large shaded error bars depicted on the plot.
Although those two types of algorithms do not have convergence guarantees for this CSO problem as the theoretical variance of these two gradient estimators is unbounded, in practice, we observe in some scenarios that they also have relatively good performance: for the housing dataset, their performance is worse than RT-MLMC/V-MLMC method but better than $L$-SGD; in other cases, they have slightly worse performance than $L$-SGD.
We do not recommend using RR-MLMC/RU-MLMC algorithms for CSO problems due to their unstable performance.
\begin{figure}[!ht]
    \centering
   \subfigure[Results for Gas Dataset]{\includegraphics[width=0.9\textwidth]{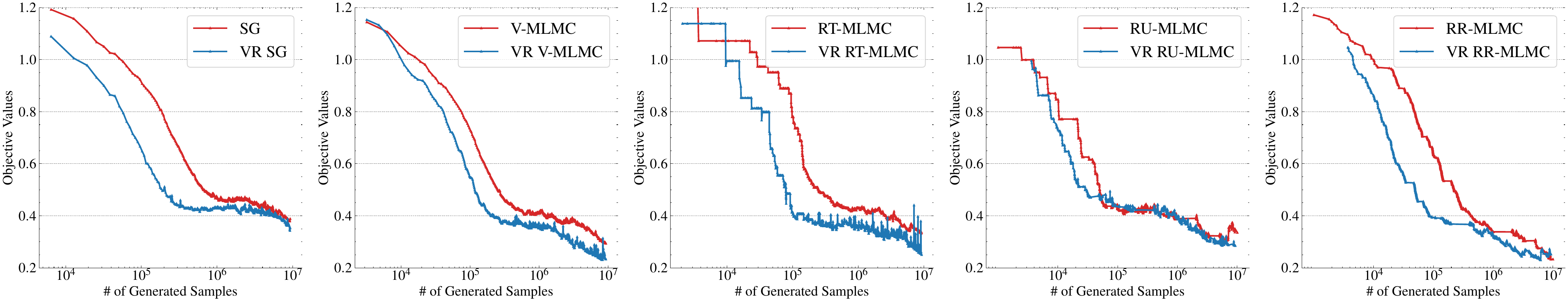}\label{fig:sub1}}
  \hfill
  \subfigure[Results for bike Dataset]{\includegraphics[width=0.9\textwidth]{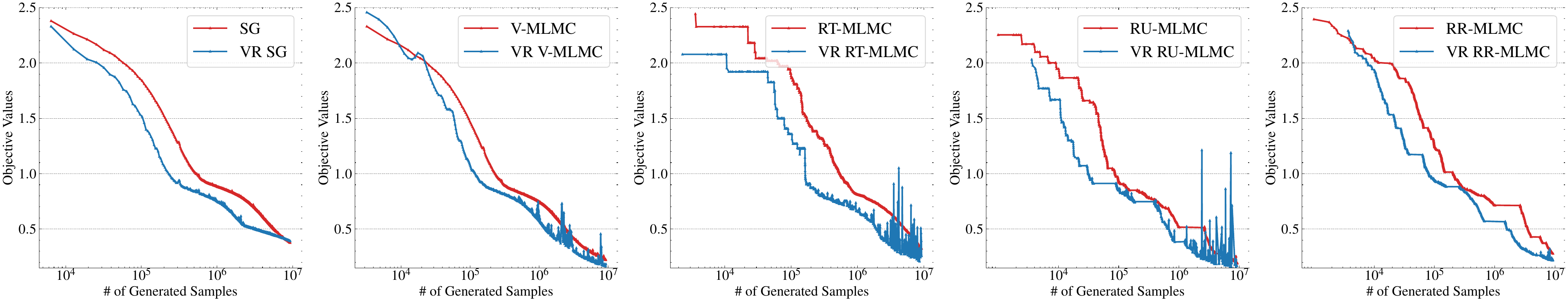}\label{fig:sub2}}

    \caption{Comparison results of various gradient methods versus their variance reduction counterparts when $f_x(z)$ is nonconvex.%
    }
    \label{fig:new:revision:noncvx}
    \vspace{-2em}
\end{figure}
\begin{figure}[!ht]
  \begin{center}
  \subfigure[Results for Gas Dataset]{
    \includegraphics[width=0.2\textwidth]{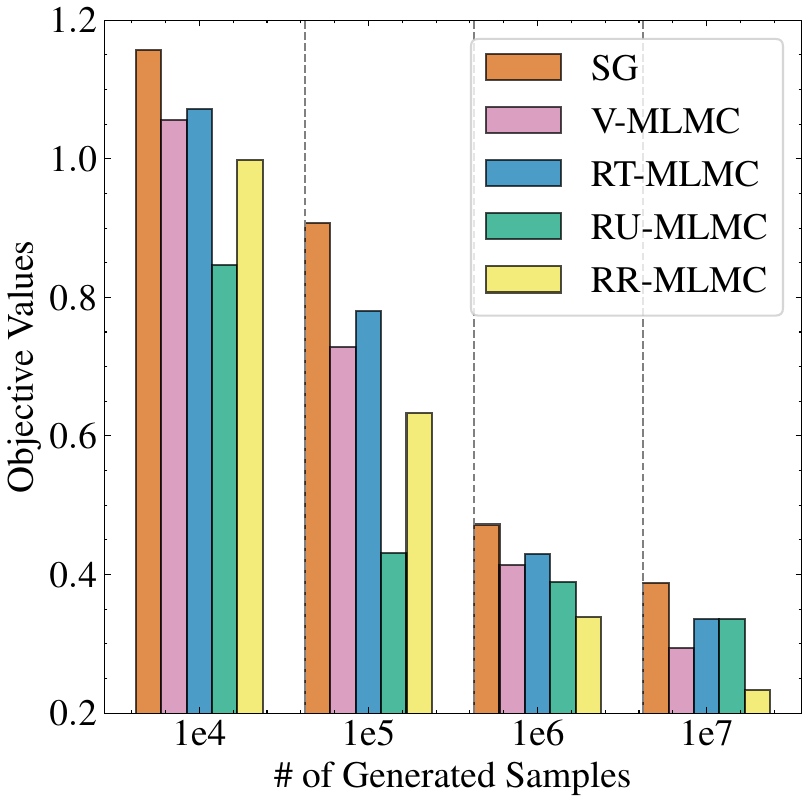}
     \includegraphics[width=0.2\textwidth]{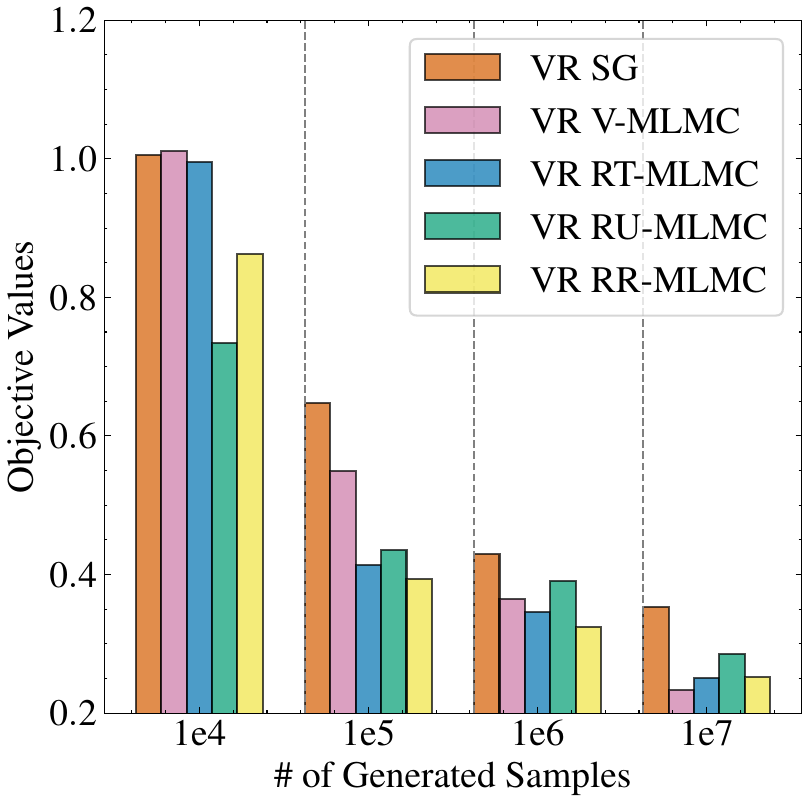}}
     \subfigure[Results for bike Dataset]{
      \includegraphics[width=0.2\textwidth]{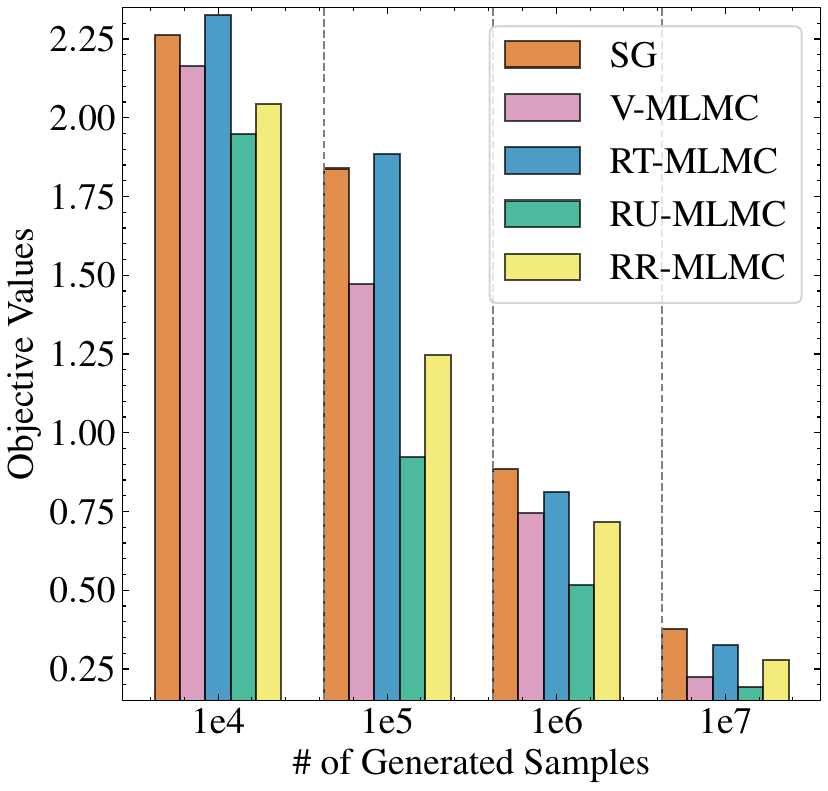}
     \includegraphics[width=0.2\textwidth]{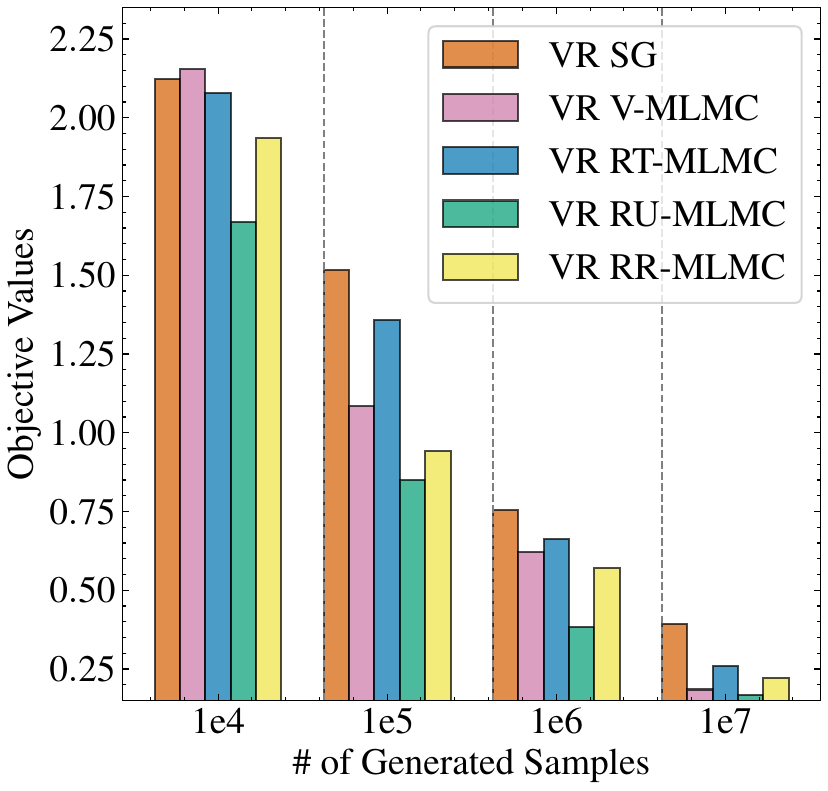}}
  \end{center}%
  \caption{%
    Left: comparison results of $L$-SGD, V-MLMC, RT-MLMC, RU-MLMC, and RR-MLMC methods when $f_x(z)$ is nonconvex; 
    Right: comparison results of variance reduction counterparts of MLMC gradient estimators.}
    \label{fig:new:revision:noncvx:summary}
    \vspace{-1.5em}
\end{figure}

In each subplot of Figure~\ref{fig:new:revision:noncvx}, we compare the performance of specific gradient algorithms and their corresponding variance reduction counterparts in the nonconvex setting. 
These subplots indicate that the variance reduction techniques further accelerate the convergence of MLMC.
Note that while we are unable to provide theoretical convergence guarantees for RU-MLMC/RR-MLMC and their variance reduction counterparts for solving the CSO problem, empirical observation demonstrates their convergence in this numerical example, with enhanced performance through variance reduction methods.

In each subplot of Figure~\ref{fig:new:revision:noncvx:summary},
we show the objective values of various gradient algorithms when their generated number of samples has reached $\texttt{1e4},\texttt{1e5},\texttt{1e6},\texttt{1e7}$, respectively.
The left plots correspond to five MLMC gradient methods, whereas the right plots correspond to their variance reduction counterparts. 
From these plots, we find, in general, V-MLMC/RT-MLMC (and their corresponding variance reduction adaptations) outperform the $L$-SGD (or VR $L$-SGD) estimator.
Besides, RU-MLMC/RR-MLMC and their variance reduction counterparts also exhibit satisfactory performance. 
Providing theoretical guarantees for these unbiased gradient estimators for solving the CSO problem remains an open challenge. It is plausible that, under more stringent theoretical assumptions, these gradient estimators lead to sharper sample complexities.

\vspace{-0.5em}
\subsection{Contrastive Learning}
\label{Sec:con:learning}

\begin{figure}[!ht]
\vspace{-2em}
    \centering
    \includegraphics[width=0.7\textwidth]{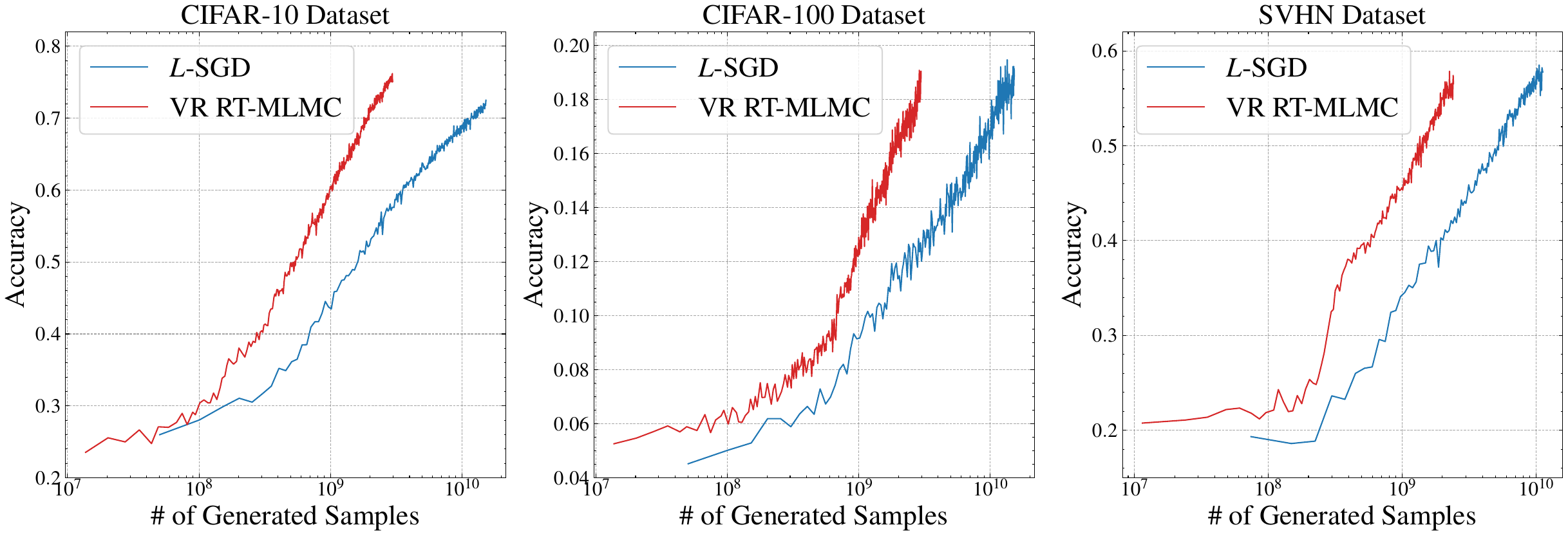}
    \caption{Comparison results of $L$-SGD and VR RT-MLMC methods on contrastive learning with CIFAR-10, CIFAR-100, and SVHN datasets.}
    \label{fig:contrastive}
   \vspace{-2em}
\end{figure}

Next, we examine the performance of $L$-SGD and VR RT-MLMC gradient estimators in contrastive learning~(See Section~\ref{Sec:CSO:con:learn}) using three large-scale datasets: CIFAR10~\citep{krizhevsky2009learning}, CIFAR100~\citep{krizhevsky2009learning}, and SVHN~\citep{goodfellow2013multi}.
We train at batch size $N=512$ (corresponding to $L=9$) for 100 epochs, where each epoch means we pass through all the training data points once.
We examine the performance at each $0.5$ epoch interval.
Experiment results are in Figure~\ref{fig:contrastive}, and other details regarding datasets, training, and testing procedure are in Table~\ref{Tab:achievable:rate:summary}.
Implementation details are provided in Appendix~\ref{Appendix:con:learn}.

Figure~\ref{fig:contrastive} shows that the VR MLMC method accomplishes training with a substantially lower number of generated samples within the same 100-epoch timeframe while maintaining comparable classification accuracy to the $L$-SGD baseline. 
Notably, the $L$-SGD method typically demands a significantly larger number of samples to attain an equivalent level of accuracy as the VR RT-MLMC method. 
The plots also reveal that optimal performance is not attained within 100 epochs.
From  Table~\ref{Tab:achievable:rate:summary}, we note that the $L$-SGD method takes an average of 25 hours to complete 100-epoch training, whereas the VR RT-MLMC method accomplishes the same in only 3 hours. Despite this time discrepancy, both methods yield reasonably good classification accuracy, enabling fair and efficient comparisons between those two methods.%

\begin{table}[!ht]
\vspace{-1em}
\caption{
Details regarding datasets, training, and testing procedure for CL} 
\label{Tab:achievable:rate:summary}
    \centering
        \begin{tabular}{c|cc|cc|c|c}
        \hline
        
        \multirow{2}{*}{Dataset} & \multicolumn{2}{c|}{$L$-SGD} & \multicolumn{2}{c|}{VR RT-MLMC} & \multirow{2}{*}{$\begin{array}{c}
           \text{Training/Testing}\\\text{Sample Size}
        \end{array}$} & \multirow{2}{*}{$\begin{array}{c}
           \text{Number of}\\\text{Labels}
        \end{array}$} \\[4pt] 
        \cline{2-5} 
        & $\begin{array}{c}
           \text{Testing}\\\text{Accuracy}
        \end{array}$ & $\begin{array}{c}
           \text{Training}\\\text{Time}
        \end{array}$ & $\begin{array}{c}
           \text{Testing}\\\text{Accuracy}
        \end{array}$ & $\begin{array}{c}
           \text{Training}\\\text{Time}
        \end{array}$ & & \\[4pt]
        \hline
       CIFAR-10  & 72.50\% & 23h50m &  75.12\% & 2h51m& 50000/10000&10\\[4pt]
        \hline
       CIFAR-100  & 18.88\% & 23h50m& 19.90\% & 2h51m& 50000/10000&100\\[4pt]
        \hline
       SVHN  & 57.72\% & 27h34m& 56.52\%& 3h22m& 50257/49032&10\\[4pt]
       \hline
    \end{tabular}
  \vspace{-1em}
    \label{tab:mytable}
\end{table}

\subsection{Pricing and Capacity
Sizing for Stochastic Systems}
\label{Sec:pricing:capacity}

Finally, we validate the performance of $L$-SGD and VR RT-MLMC estimator for the joint pricing and staffing problem \eqref{Eq:pricing:capacity:formula} based on the approximate gradient in \eqref{Eq:approx:Fl:ss} with maximum level $L=18$.
In this case, the $L$-SGD estimator requires $\cO(2^L)\approx \cO(10^{5})$ operations to construct a single gradient, which is a costly choice.
We follow the setup in \cite{chen2023online} to perform the numerical study.
We model (normalized) inter-arrival time $U_n$ to follow the exponential distribution and the (normalized) service time $V_n$ to follow the Erlang, exponential, or hyper-exponential distribution.

\begin{figure}[!ht]
\vspace{-1em}
  \begin{center}
  \subfigure[Erlang distribution]{
    \includegraphics[width=0.23\textwidth]{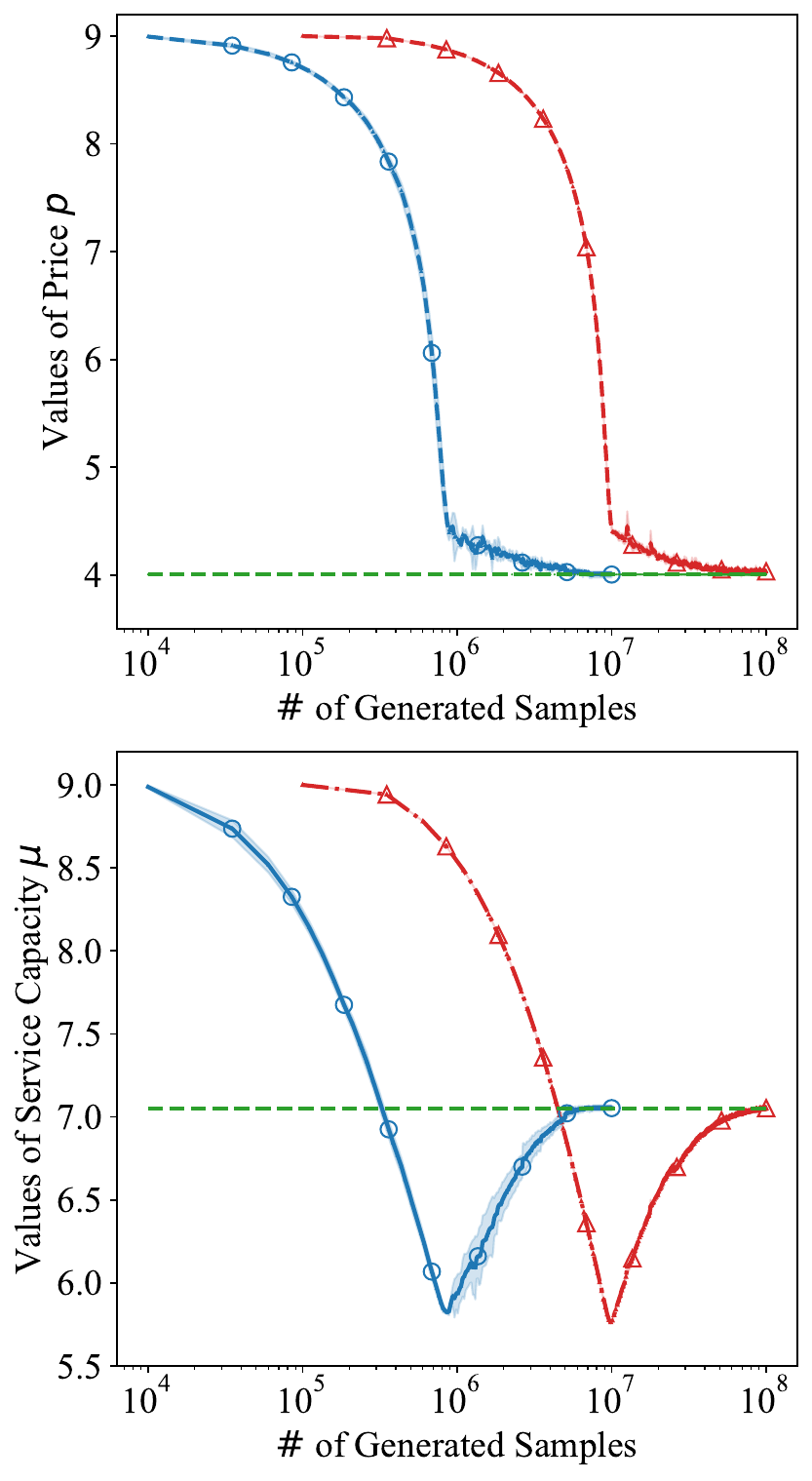} %
     }
      \subfigure[Exponential distribution]{
    \includegraphics[width=0.23\textwidth]{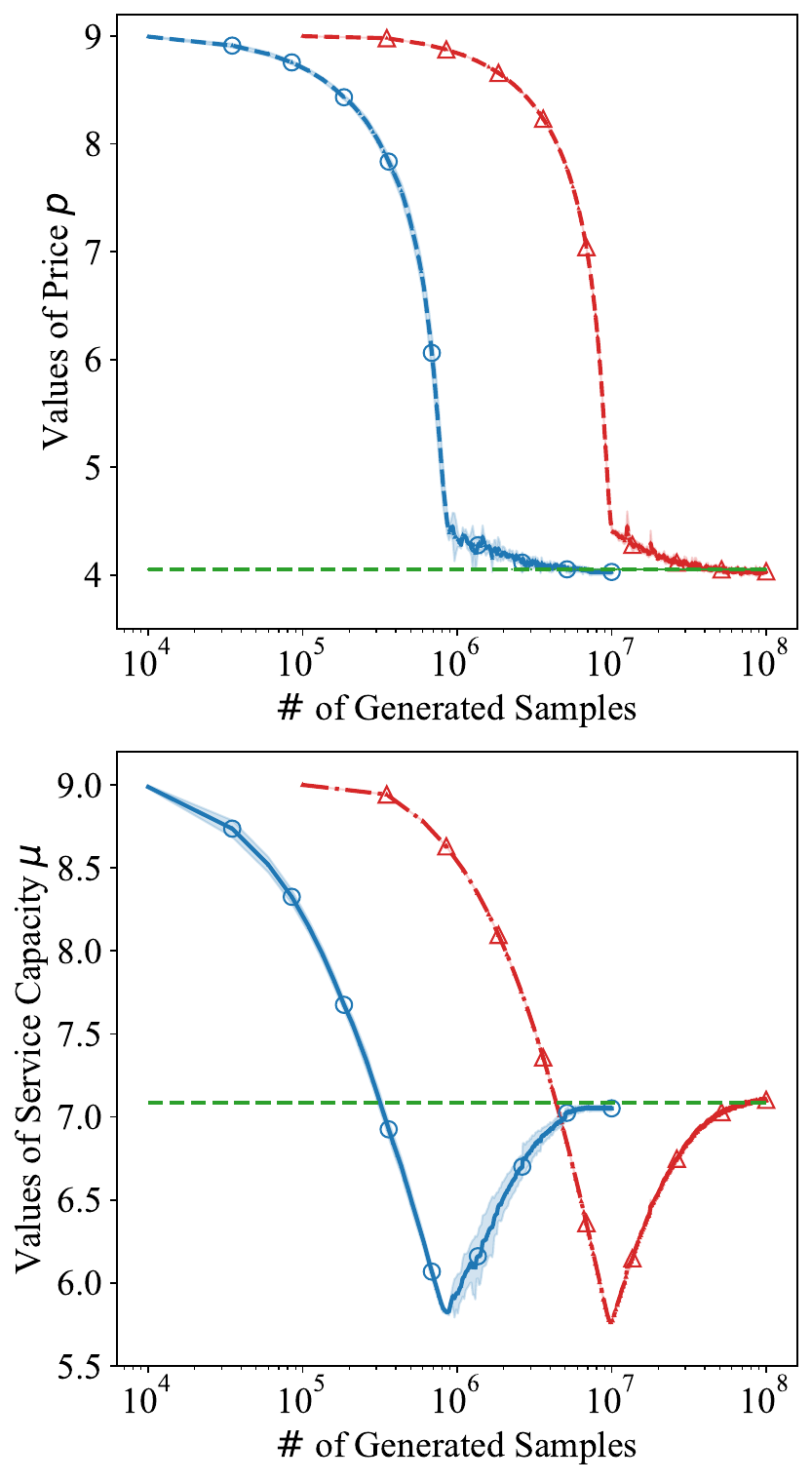}%
     }
     \subfigure[Hyperexponential distribution]{
    \includegraphics[width=0.348\textwidth]{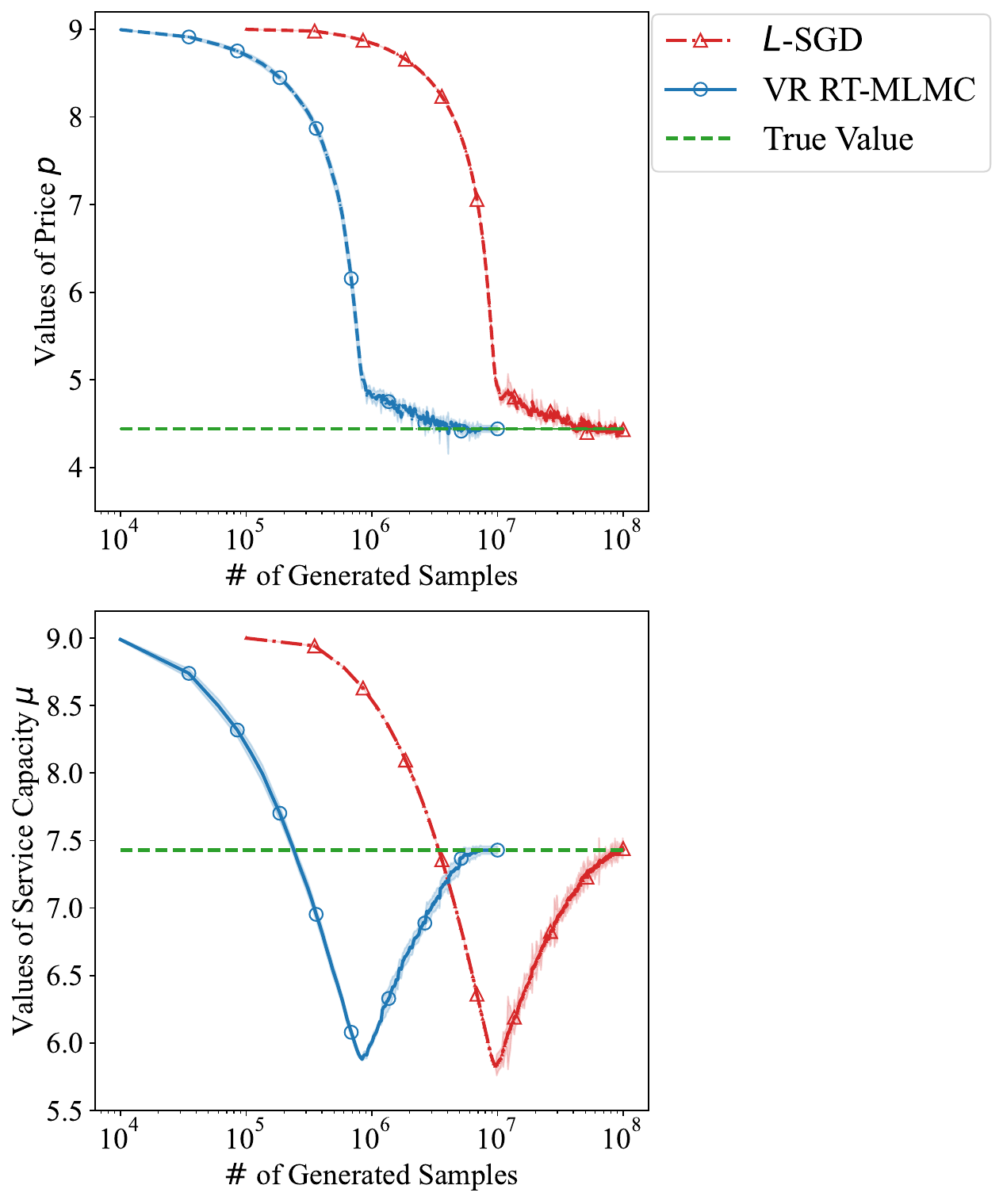}
     }
  \end{center}%
  \caption{%
    Comparison results of $L$-SGD and VR RT-MLMC estimators for the joint pricing and staffing task in a stochastic system with various service time distributions.
    The results are averaged with error bars based on $10$ independent runs.}
    \label{fig:queue}
    \vspace{-1.5em}
\end{figure}

Figure~\ref{fig:queue} reports the performance of $L$-SGD and VR RT-MLMC gradient estimators under three different service time distributions.
The $x$-axis corresponds to the number of generated samples, the $y$-axis corresponds to either the estimated price or service capacity and the shaded areas correspond to the error bars generated using $10$ independent trials.
For those three special instances, we have a closed-form reformulation of the objective function. We use an exhaustive search to obtain the ground-truth optimal choices of price and capacity parameters, which are the green dashed lines in those plots.

Despite the non-convexity of the objective functions, we find these gradient methods converge to the global optimum solution in practice, which aligns with the observation in \citet{chen2023online}, and its global convergence analysis remains an interesting theoretical question for future investigation.
From those plots, we also find that the VR RT-MLMC method obtains the optimal solution with a significantly smaller number of generated samples, which indicates it is a powerful choice for the joint pricing and staffing task.

\section{Conclusion}
\label{section:conclusion}
This paper systematically studies the bias-variance-cost tradeoff of several MLMC gradient methods for stochastic optimization under a generic biased oracle model, shedding light on their superiority and limitations under different situations. 
We further incorporate variance reduction methods to reduce the (expected) total cost in the nonconvex smooth case. 
We demonstrate the significant benefits of MLMC gradient methods for various applications in theory and numerical study.  
Several open questions remain for future work. First, it is promising to design gradient estimators using MLMC that adaptively reduce the bias when the iteration increases.
Second, it remains interesting to understand if $b<c$, i.e., the increase rate of the costs is higher than the decrease rate of the variance, what are the fundamental limits on the total costs, and if it is possible to develop more efficient algorithms. 

\clearpage
\ACKNOWLEDGMENT{Yifan Hu is supported by NCCR Automation of Swiss NSF.}
{
\bibliographystyle{plainnat}
\bibliography{main.bib}
}

\begin{APPENDICES}

\newpage
\renewcommand{\theHsection}{A\arabic{section}}
\begin{APPENDICES}
\ECHead{Online Appendices}

\section{Preliminaries on SGD}
\label{app_sec:SGD}
In this section, we show the convergence property of the SGD framework when the gradient estimator $v(x)$ is unbiased and has bounded variance. 
Similar results appear in \citet{nemirovski2009robust,rakhlin2012making,bottou2018optimization}. 
Lemma \ref{lm:SGD} gives a summary of the results that 
we in this paper.

Note that $L$-SGD, V-MLMC, and RT-MLMC are biased gradient estimators for optimizing $F$ but unbiased estimators for optimizing  $F^L$. We use Lemma \ref{lm:SGD} together with Assumption \ref{assumption:general}\ref{assumption:general:I} and relation~\eqref{eq:error_decomposition} to derive the global convergence of biased methods on $F$ in the PL condition or convex case; with Assumption~\ref{assumption:gradient_bias} and relation \eqref{eq:truncate_nonconvex_error} to derive the stationary convergence of biased methods on $F$ in the nonconvex smooth case.
Notably, $\mu$-strong convexity is a sufficient condition of PL condition with parameter $\mu$~\citep{karimi2016linear}.

\begin{lemma}[Convergence of SGD]
\label{lm:SGD}
Let $\{x_t\}_{t=1}^T$ be the updating trajectory of Algorithm \ref{alg:SGD-framework}, $x^*$ be a minimizer of $F(x)$ on $\RR^d$, and  $ \hat x_T $ be selected uniformly randomly from $\{x_t\}_{t=1}^T$.
Suppose that  there exists a constant $V>0$ such that $\EE v(x) = \nabla F(x), \mathrm{Var}(v(x))\leq \mathbf{V},$
where expectations are taken with respect to~(w.r.t.) the randomness in $v$. 
If $F(x)$ is $S_F$-smooth, the following results hold.
\begin{enumerate}
    \item[\mylabel{Lemma:EC1:a1}{$\mathrm{(I}$-$\mathrm{a)}$}] If $F(x)$ satisfies PL condition with parameter $\mu$, for fixed stepsizes $\gamma_t = \gamma \in(0, \frac{1}{S_F}]$, we have
    \begin{align*}
        \EE[F(x_T)-F(x^*)] 
    \leq 
        (1-\gamma \mu)^{T-1}[F(x_1)-F(x^*)]+\frac{S_F\gamma \mathbf{V}}{2\mu}.
    \end{align*}
    \item[\mylabel{Lemma:EC1:a2}{$\mathrm{(I}$-$\mathrm{b)}$}]  If $F(x)$ satisfies PL condition with parameter $\mu$, for time-varying stepsizes $\gamma_t = \frac{2}{\mu(t+2S_F/\mu-1)}$, we have
    \begin{align*}
        \EE [F(x_{T})-F(x^*)]
    \leq
        \frac{2\max\{S_F\mathbf{V}, 2\mu S_F(F(x_1)-F(x^*))\}}{\mu^2 (T+2S_F/\mu-1)}.
    \end{align*}
    \item[\mylabel{Lemma:EC1:b}{$\mathrm{(II)}$}] If $F(x)$ is convex,  for stepsizes $\gamma_t=\gamma\in(0, \frac{1}{2S_F}]$, we have
    \begin{align*}
        \EE [F(\hat x_T)- F(x^*)] 
    \leq  
        \gamma \mathbf{V}+\frac{\|x_1-x^*\|_2^2}{\gamma T}.
    \end{align*}
    \item[\mylabel{Lemma:EC1:c}{$\mathrm{(III)}$}] If $F(x)$ is nonconvex, for fixed stepsizes $\gamma_t =\gamma\in(0, \frac{1}{S_F}]$, we have 
    \begin{align*}
        \EE \|\nabla F(\hat x_T)\|_2^2 
    \leq 
        \frac{2(F(x_1)- F(x^*))}{\gamma T} + S_F\gamma\mathbf{V}.
    \end{align*}
\end{enumerate}

\begin{enumerate}
  \item[\mylabel{Lemma:EC1:d}{$\mathrm{(IV)}$}]
    If $F(x)$ is convex and $L_F$-Lipschitz continuous, it holds that
    \begin{align*}
        \EE [F(\hat x_T)- F(x^*)] 
    \leq  
        \frac{\|x_1-x^*\|_2^2}{2\gamma T} + \frac{\gamma(L_F^2+\mathbf{V})}{2}.
    \end{align*}
\end{enumerate}
\end{lemma}
Note that we do not specify the stepsizes in cases \ref{Lemma:EC1:a1}, \ref{Lemma:EC1:a2}, \ref{Lemma:EC1:c}, and \ref{Lemma:EC1:d}. When the variance of $v(x)$ is of order $\cO(\eps)$, one can use stepsizes that are independent of $\eps$ to guarantee $\eps$-optimality or $\eps$-stationarity. The algorithm would behave similarly to gradient descent. On the other hand, when the variance of $v(x)$ is of order $\cO(1)$ or even larger, one should use stepsizes $\gamma_t= \cO(1/t)$ in the strongly convex setting and $\gamma_t = \cO(1/\sqrt{T})$ in the convex or nonconvex smooth setting. 
\begin{proof}{Proof of Lemma~\ref{lm:SGD}.}
We divide the proof into the smooth and the Lipschitz continuous cases.\\
\noindent \textbf{Part I: Smooth Case.}
Since $x^*$ is a minimizer of $F$ over $\RR^d$, we have $\nabla F(x^*) =0$.
\subsubsection*{PL condition and smooth case~\ref{Lemma:EC1:a1}.}
By smoothness, it holds that
\begin{align*}
        F(x_{t+1}) - F(x_t) 
    \leq & 
        \nabla F(x_t)^\top (x_{t+1}-x_t) + \frac{S_F}{2}\|x_{t+1} - x_{t}\|_2^2 
    =   
        -\gamma_t \nabla F(x_t)^\top v(x_t) +\frac{S_F\gamma_t^2}{2}\|v(x_t)\|_2^2.
\end{align*}
Taking expectation conditioned on $x_t$ on both side, since $\EE [v(x_t)|x_t]=\nabla F(x_t)$ it holds that
\begin{equation}
\begin{split}
\label{eq:stronglyconvex_1}
        \EE [F(x_{t+1})-F(x_t)|x_t]
    \leq &
        -\gamma_t \|\nabla F(x_t)\|_2^2 +\frac{S_F\gamma_t^2}{2}\EE[\|v(x_t)\|_2^2|x_t]\\
    = & 
        -\gamma_t \|\nabla F(x_t)\|_2^2 +\frac{S_F\gamma_t^2}{2}\|\nabla F(x_t)\|_2^2 +\frac{S_F\gamma_t^2}{2}\mathrm{Var}(\|v(x_t)\|_2^2|x_t)\\
    = & 
        -\Big(\gamma_t-\frac{S_F\gamma_t^2}{2}\Big)\|\nabla F(x_t)\|_2^2+\frac{S_F\gamma_t^2}{2}\mathrm{Var}(\|v(x_t)\|_2^2|x_t) \\
    = &
        -\Big(1-\frac{S_F\gamma_t}{2}\Big)\gamma_t\|\nabla F(x_t)\|_2^2+\frac{S_F\gamma_t^2}{2}\mathrm{Var}(\|v(x_t)\|_2^2|x_t) \\
    \leq & 
        -\frac{\gamma_t}{2}\|\nabla F(x_t)\|_2^2+\frac{S_F\gamma_t^2}{2}\mathrm{Var}(\|v(x_t)\|_2^2|x_t) \\
    \leq &
        -\gamma_t \mu(F(x_t)-F(x^*))+ \frac{S_F\gamma_t^2}{2}\mathrm{Var}(\|v(x_t)\|_2^2|x_t),
\end{split}
\end{equation}
where the second to last inequality uses the assumption that $\gamma_t \leq \frac{1}{S_F}$, the last inequality uses PL condition.

Subtracting $F(x^*)$ on both sides of \eqref{eq:stronglyconvex_1} and taking full expectation, we have
\begin{equation}
\label{eq:sc-recursion}
\begin{aligned}
        \EE [F(x_{t+1})-F(x^*)] 
    \leq & 
        (1-\gamma_t \mu)\EE[F(x_t)-F(x^*)] +\frac{S_F\gamma_t^2}{2}\mathrm{Var}(\|v(x_t)\|_2^2)\\
    \leq &
        (1-\gamma_t \mu)\EE[F(x_t)-F(x^*)] +\frac{S_F\gamma_t^2}{2}\mathbf{V}.
\end{aligned}   
\end{equation}
For fixed $\gamma_t=\gamma \leq \frac{1}{S_F}$, by induction, it is easy to see
\begin{align*}
    \EE[F(x_T)-F(x^*)] 
    \leq (1-\gamma \mu)^{T-1}[F(x_1)-F(x^*)]+\frac{S_F\gamma\mathbf{V}}{2\mu}.
\end{align*}

\subsubsection*{PL condition and smooth case~\ref{Lemma:EC1:a2}.}
Starting from \eqref{eq:sc-recursion} and substituting $\gamma_t = \frac{2}{\mu (t+\mathfrak{a})}$ with $\mathfrak{a}:=2S_F/\mu-1$ such that $\gamma_t\in(0,1/S_F]$, by induction, it is easy to see
\begin{align*}
    \EE [F(x_{T})-F(x^*)]
\leq
    \frac{2\max\{S_F\mathbf{V}, \mu^2(1+\mathfrak{a})(F(x_1)-F(x^*))\}}{\mu^2 (T+\mathfrak{a})}.
\end{align*}

\subsubsection*{Convex smooth Case~\ref{Lemma:EC1:b}.}

Denote $a_t = \frac{1}{2}\EE\norm{x_t-x^*}_2^2$, then it holds that
\begin{equation}
\label{eq:recursion}
\begin{split}
        a_{t+1} 
    & 
    =  
        \frac{1}{2}\EE\norm{x_{t}-\gamma_t v(x_t)-x^*}_2^2 
    = 
        a_t + \frac{1}{2}\gamma_t^2\EE \norm{v(x_t)}_2^2 - \gamma_t \EE v(x_t)^\top(x_{t}-x^*).
\end{split}
\end{equation}
Dividing $\gamma_t$ on both sides of \eqref{eq:recursion}, it holds that
\begin{equation*}
    \EE v(x_t)^\top(x_{t}-x^*)\leq \frac{a_t-a_{t+1}}{\gamma_t}  + \frac{1}{2}\gamma_t\EE \norm{v(x_t)}_2^2.
\end{equation*}
Since $F(x)$ is convex and $v(x_t)$ is an unbiased gradient estimator of $F(x_t)$ conditioned on $x_t$, we have
\begin{equation*}
\begin{split}
    -\EE v(x_t)^\top (x_t-x^*) = - \nabla F(x_t)^\top (x_t-x^*)
    \leq F(x^*)-F(x_t).
\end{split}
\end{equation*}
Summing up the two inequalities above and rearranging, it holds  that conditioned on $x_t$, 
\begin{equation*}
\begin{split}
    & 
        F(x_t)-F(x^*)
    \leq 
        \frac{a_t-a_{t+1}}{\gamma_t}  + \frac{1}{2}\gamma_t\EE \norm{v(x_t)}_2^2
    \leq 
        \frac{a_t-a_{t+1}}{\gamma_t}  + \frac{1}{2}\gamma_t \mathrm{Var}(v(x_t))+\frac{1}{2}\gamma_t\|\nabla F(x_t)\|_2^2\\
    = & 
        \frac{a_t-a_{t+1}}{\gamma_t}  + \frac{1}{2}\gamma_t \mathrm{Var}(v(x_t))+\frac{1}{2}\gamma_t\|\nabla F(x_t)-\nabla F(x^*)\|_2^2\\
    \leq & 
        \frac{a_t-a_{t+1}}{\gamma_t}  + \frac{1}{2}\gamma_t \mathrm{Var}(v(x_t))+\gamma_t S_F[F(x_t)-F(x^*)],
\end{split}
\end{equation*}
where the last inequality holds as $F$ is convex and $S_F$-smooth. It further implies
\begin{align*}
        \EE [F(x_t)-F(x^*)]
    \leq 
        \frac{1}{1-\gamma_t S_F}\Big[\frac{a_t-a_{t+1}}{\gamma_t}  + \frac{1}{2}\gamma_t \mathrm{Var}(v(x_t))\Big].
\end{align*}
For $\gamma_t\leq \frac{1}{2S_F}$, it holds that $\frac{1}{1-\gamma_t S_F}\leq 2$.
By definition of $\hat x_T$, Jensen's inequality, and convexity of $F$,  %
\begin{equation}
\begin{split}
    & 
        \EE[F(\hat{x}_T)-F(x^*)]  
    \leq
        \frac{1}{T} \sum_{t=1}^T\EE[F(x_t)-F(x^*)]\\
    \leq &
        \frac{1}{T}\sum_{t=1}^T\gamma_t \mathrm{Var}(v(x_t)) + \frac{1}{T}\sum_{t=2}^T \|x_t-x^*\|_2^2\left(\frac{1}{\gamma_{t}} -\frac{1}{\gamma_{t-1}}\right) +\frac{1}{\gamma_{1}T}\|x_1-x^*\|_2^2
    \leq 
        \gamma\mathbf{V}  +\frac{\|x_1-x^*\|_2^2}{\gamma T}. 
\end{split}
\end{equation}
where the last inequality holds as $\gamma_t = \gamma$.

\subsubsection*{Nonconvex smooth case~\ref{Lemma:EC1:c}.}
Since $F(x)$ is $S_F$-smooth, it holds that 
\begin{equation*}
\begin{split}
    &
        F(x_{t+1})-F(x_t) 
    \leq 
        \nabla F(x_t)^\top (x_{t+1}-x_t) + \frac{S_F}{2}\|x_{t+1}-x_t\|_2^2
    =
        -\nabla F(x_t)^\top \gamma_t v(x_t) + \frac{S_F\gamma_t^2}{2}\|v(x_t)\|_2^2.
\end{split}
\end{equation*}
Taking expectations on both sides, it holds that
\begin{align*}
 \EE [F(x_{t+1})-F(x_t)] 
    &\leq 
        -\gamma_t \EE \|\nabla F(x_t)\|_2^2 +\frac{S_F\gamma_t^2}{2}\EE \|v(x_t)\|_2^2 \\
    &= 
        \Big(-{\gamma_t} +\frac{S_F\gamma_t^2}{2}\Big)\EE \|\nabla F(x_t)\|_2^2 +\frac{S_F\gamma_t^2}{2}\mathrm{Var}(v(x_t)).
\end{align*}
Summing up from $t=1$ to $t=T$,  taking full expectation, and setting $\gamma_t = \gamma \leq 1/S_F$, we have
\begin{equation*}
\begin{split}
    &
        \EE [F(x_{T+1})-F(x_1)] 
    \leq 
        \Big(-\gamma+\frac{S_F\gamma^2}{2}\Big)\sum_{t=1}^{T}  \EE\|\nabla F(x_t)\|_2^2 +\sum_{t=0}^{T-1}\frac{S_F\gamma^2}{2}\mathrm{Var}(v(x_t))\\
    \leq & 
        -\frac{\gamma}{2}\sum_{t=1}^{T}  \EE\|\nabla F(x_t)\|_2^2 +\sum_{t=1}^{T}\frac{S_F\gamma^2}{2}\mathrm{Var}(v(x_t)).
\end{split}
\end{equation*}
As a result,%
\begin{equation*}
\begin{split}
        \EE \|\nabla F(\hat x_T)\|_2^2
    \leq &
        \frac{1}{T}\sum_{t=1}^{T}  \EE\|\nabla F(x_t)\|_2^2
    \leq 
        \frac{2\EE(F(x_1)-F(x_{T+1}))}{\gamma T}+ S_F\gamma \mathbf{V}
    \leq 
        \frac{2(F(x_1)-F(x^*))}{\gamma T}+ S_F\gamma\mathbf{V}.
\end{split}
\end{equation*}

\noindent \textbf{Part II: Lipschitz Continuous Case.}
If $F$ is $L_F$-Lipschitz continuous, the second moment of the gradient estimator $v$ is upper bounded:
$$
\EE \|v(x)\|_2^2 =\mathrm{Var}(v(x)) + \|\EE v(x)\|_2^2 =  \mathrm{Var}(v(x)) + \|\nabla F(x)\|_2^2 \leq L_F^2+ \mathbf{V}.
$$

\subsubsection*{Convex Lipschitz continuous case~\ref{Lemma:EC1:d}.}
Dividing $\gamma_t$ on both sides of \eqref{eq:recursion}, it holds that
\begin{equation*}
    \EE v(x_t)^\top(x_{t}-x^*)\leq \frac{a_t-a_{t+1}}{\gamma_t}  + \frac{1}{2}\gamma_t\EE \norm{v(x_t)}_2^2.
\end{equation*}
By convexity and upper bounds on the second moment, we have 
$$
\EE [F(x_t)-F(x^*)]\leq \frac{a_t-a_{t+1}}{\gamma_t}  + \frac{1}{2}\gamma_t(L_F^2+\mathbf{V}).
$$
Summing up from $t=1$ to $t=T$ and setting $\gamma_t = \gamma$, we have
\begin{equation}
\EE [F(\hat x_T)-F(x^*)]\leq \EE\frac{1}{T}\sum_{t=1}^T[F(x_t)-F(x^*)]\leq \frac{\|x_1-x^*\|_2^2}{2\gamma T} + \frac{\gamma(L_F^2+ \mathbf{V})}{2}.
\tag*{\QED}
\end{equation}
\end{proof}

\section{Discussion on Assumptions}
Since the PL condition is a generalization of strong convexity, the (expected) total cost of $L$-SGD, V-MLMC, RT-MLMC, RU-MLMC, and RR-MLMC under PL condition provides an upper bound on the (expected) total cost achieved under the strong convexity assumption. For a more detailed discussion on PL condition and strong convexity, see \citet{karimi2016linear}.

\subsection{Substituteable Bias Assumption under PL Condition}
\label{appendix:pl-condition}
Recall that  Assumption \ref{assumption:general}\ref{assumption:general:I} assumes that the bias of the function value estimator is uniformly upper-bounded. In contrast, Assumption \ref{assumption:gradient_bias} assumes that the bias of the gradient estimator is uniformly upper-bounded.
In this part, we show Assumptions~\ref{assumption:general}\ref{assumption:general:I} and \ref{assumption:gradient_bias} are substitutable under PL condition (strong convexity) setting in the sense that the (expected) total cost of $L$-SGD, V-MLMC, and RT-MLMC would remain the same under either assumption. 
However, Assumption \ref{assumption:gradient_bias} cannot replace Assumption \ref{assumption:general}\ref{assumption:general:I} for convex case. See Remark~\ref{Remark:counter} for a counter-example.

Let $x_T^{\texttt{A}}$ denote the $T$-th iteration of the algorithm using $\texttt{A}\in\{\text{$L$-SGD, V-MLMC, RT-MLMC}\}$.
The key step of such replacement is to show that under Assumption \ref{assumption:gradient_bias}, the expected error of the algorithm using $\texttt{A}$, i.e.,
$
 \EE [F(x_T^\texttt{A})-F(x^*)]    
$,
has a similar error decomposition like \eqref{eq:error_decomposition}. 

\begin{proposition}\label{Proposition:PL:SGD}
Suppose that $F$ is $S_F$-smooth and satisfies the PL condition with parameter $\mu$, under Assumptions \ref{assumption:general}\ref{assumption:general:II}\ref{assumption:general:III} and Assumption \ref{assumption:gradient_bias} and take $\mathfrak{a}=2S_F/\mu-1$, then the following results hold.
\begin{enumerate}%
    \item When using fixed stepsizes $\gamma_t = \gamma\in(0,1/S_F]$, %
    \begin{align}
    \label{eq:compare_gradient}
       \EE [F(x_T^\texttt{A})-F(x^*)]    
    \leq 
        (1-\gamma \mu)^{T-1} [F(x_1)-F(x^*)] + \frac{1}{2\mu}M_a2^{-aL} + \frac{S_F \gamma}{2\mu}\mathbf{V}(v^\texttt{A}).
    \end{align}
    \item When using time-varying stepsizes $\gamma_t = \frac{2}{\mu(t+\mathfrak{a})}$, we have
    \begin{equation}
    \label{eq:compare_gradient_2}
   \EE [F(x_T^\texttt{A})-F(x^*)]    
\leq 
    \frac{2\max\{S_F\mathbf{V}(v^\texttt{A}), \mu^2(1+\mathfrak{a})(F(x_1)-F(x^*))\}}{\mu^2 (T+\mathfrak{a})} + \frac{1}{2\mu}M_a 2^{-aL}.
    \end{equation}
\end{enumerate}
\end{proposition}

\begin{proof}{Proof of Proposition~\ref{Proposition:PL:SGD}.} 
By smoothness of $F$ and taking expectation conditioned on $x_t$, we have 
\begin{align*}
    \EE F(x_{t+1}^A)
\leq & 
    F(x_t^A)+\nabla F(x_t^A)^\top \EE (x_{t+1}^A-x_t^A)+\frac{S_F}{2}\EE\|x_{t+1}^A-x_t^A\|_2^2 \\
= &
    F(x_t^A)+\nabla F(x_t^A)^\top \EE (x_{t+1}^A-x_t^A)+\frac{S_F \gamma_t^2}{2}\EE\|v(x_t^A)\|_2^2  \\
= &
    F(x_t^A)-\frac{\gamma_t}{2}(2\nabla F(x_t^A)^\top\EE v(x_t^A)) +\frac{S_F \gamma_t^2}{2}\|\EE v(x_t^A)\|_2^2+ \frac{S_F \gamma_t^2}{2}\mathrm{Var}( v(x_t^A)) \\
\leq &
    F(x_t^A)+\frac{\gamma_t}{2}(-2\nabla F(x_t^A)^\top\EE v(x_t^A)+\|\EE v(x_t^A)\|_2^2)+ \frac{S_F \gamma_t^2}{2}\mathrm{Var}( v(x_t^A)) \\
= &
    F(x_t^A)+\frac{\gamma_t}{2}(-\|\nabla F(x_t^A)\|_2^2+\|\EE v(x_t^A)-\nabla F(x_t^A)\|_2^2)+ \frac{S_F \gamma_t^2}{2}\mathrm{Var}( v(x_t^A))\\
\leq & 
    F(x_t^A)-\frac{\gamma_t}{2}\|\nabla F(x_t^A)\|_2^2 + \frac{\gamma_t}{2} B_L+\frac{S_F \gamma_t^2}{2}\mathbf{V}(v^\texttt{A}),
\end{align*}
where the first inequality uses smoothness, the second inequality uses $\gamma_t \leq 1/S_F$, the third inequality uses Assumption \ref{assumption:gradient_bias}, the last equality uses $\|a-b\|_2^2 = \|a\|_2^2 -2a^\top b +\|b\|_2^2$.
Using PL condition and taking full expectation, we have 
\begin{equation}
\label{eq:gradient_recursion}
\begin{split}
&
    \EE [F(x_{t+1}^\texttt{A})-F(x^*)] 
\leq 
    \EE [F(x_t^\texttt{A})-F(x^*)]-\frac{\gamma_t}{2}2\mu\EE (F(x_t^\texttt{A})-F(x^*)) + \frac{\gamma_t}{2} B_L+\frac{S_F \gamma_t^2}{2}\mathbf{V}(v^\texttt{A})\\
= &
    (1-\gamma_t\mu)\EE (F(x_t^\texttt{A})-F(x^*))+\frac{\gamma_t}{2} B_L+\frac{S_F \gamma_t^2}{2}\mathbf{V}(v^\texttt{A}).
\end{split}
\end{equation}
\begin{itemize}[leftmargin=2em]
    \item If using fixed stepsizes, plugging in $\gamma_t = \gamma\in(0,1/S_F]$, by induction, we have
\begin{align*}
    & \EE [F(x_{T}^\texttt{A})-F(x^*)] 
\leq (1-\gamma \mu)^{T-1} [F(x_1)-F(x^*)] + \frac{1}{2\mu}B_L + \frac{S_F \gamma}{2\mu}\mathbf{V}(v^\texttt{A}).%
\end{align*}

\item If using time-varying stepsizes, plugging $\gamma_t = \frac{2}{\mu (t+\mathfrak{a})}\in(0,1/S_F]$ into \eqref{eq:gradient_recursion}, by induction, we have 
\begin{equation*}
\begin{split}
&   \EE [F(x_{T}^\texttt{A})-F(x^*) ]
\leq
    \frac{2\max\{S_F\mathbf{V}(v^\texttt{A}), \mu^2(1+\mathfrak{a})(F(x_1)-F(x^*))\}}{\mu^2 (T+\mathfrak{a})} + \frac{1}{2\mu}B_L.
\end{split}
\end{equation*}
\end{itemize}
Substituting $B_L=M_a2^{-aL}$ gives the desired results.
\QED
\endproof
\end{proof}

In the previous analysis of $L$-SGD and RT-MLMC, we used time-varying stepsizes and Assumption~\ref{assumption:general}\ref{assumption:general:I}.
The key is the following relation:
\begin{align}
\label{eq:compare_function_2}
    &\EE [F(x_{T}^\texttt{A})-F(x^*)]
\leq %
    2M_a 2^{-aL} + \frac{2\max\{S_F\mathbf{V}(v^\texttt{A}), 2\mu S_F(F(x_1)-F(x^*))\}}{\mu^2 (T+2S_F/\mu-1)}.
\end{align}
Comparing \eqref{eq:compare_gradient_2} and \eqref{eq:compare_function_2}, the only differences are in the constants, which do not affect the rate of the total cost. 

In the previous analysis of V-MLMC, we used fixed stepsizes and Assumption~\ref{assumption:general}\ref{assumption:general:I}.
The key is via 
\begin{align}
\label{eq:compare_function}
    &\EE [F(x_{T}^\texttt{A})-F(x^*)]
\leq %
    2M_a 2^{-aL} + (1-\gamma \mu)^{T-1}[F^L(x_1)-F^L(x^L)]+\frac{S_F\gamma \mathbf{V}(v^\texttt{A})}{2\mu}.
\end{align}
Also, we find \eqref{eq:compare_gradient} and \eqref{eq:compare_function} lead to the same error bound but with different constants, which do not affect the rate of the total cost.

\begin{remark}\label{Remark:counter}
Unlike Assumption \ref{assumption:general}\ref{assumption:general:I}, Assumption \ref{assumption:gradient_bias} is not sufficient for obtaining a global optimality gap when solving unconstrained optimization with convex objective $F$ using biased gradient methods. Suppose that Assumption \ref{assumption:gradient_bias} holds and that one finds an $\eps/4$-stationarity point of $F^L$ via some biased methods. The point is an $\eps$-stationarity point of $F$ by Assumption \ref{assumption:gradient_bias} for certain $L$. However, there is no link between the gradient norm of the point and the optimality gap for unconstrained optimization with convex objectives. In fact, one can show that for any $\eps\in(0,1)$, there exists a convex smooth function $F:\RR^d\rightarrow\RR$ and  a point $x_0$ such that $\|\nabla F(x_0)\|_2^2=\eps^2$ and $F(x_0)-\inf_x F(x) > 1$. One example is the Huber function defined in the following, and $x_0\in\RR^d$ is such that $\|x_0\|_2=2/\eps>\eps$.
\begin{align*}
    F^\mathrm{Huber}(x) = 
    \begin{cases}
        \frac{1}{2}\|x\|_2^2 & \text{ if } \|x\|_2< \eps, \\[5pt]
        \eps (\|x\|_2 - \frac{\eps}{2}) & \text{ if } \|x\|_2\geq \eps.
    \end{cases}
\end{align*}
One can easily see $\|\nabla F(x_0)\|_2^2=\eps^2$ but $F(x_0)-\inf_x F(x) = 2-\frac{\eps^2}{2}>1$.
When encountering such functions, the biased gradient methods do not guarantee to converge to an $\eps$-optimal solution under Assumption \ref{assumption:gradient_bias}.

\end{remark}

\section{Bias, Variance, and Cost of Gradient Estimators}
\label{app_sec:bias_variance_cost}
This section demonstrates the bias, variance, and (expected) per-iteration cost of the $L$-SGD and the MLMC-based gradient estimators. The main manuscript in Table~\ref{tab:bias-variance-cost} has simplified results. %
The following lemma formally characterizes the expectation of $L$-SGD and MLMC-based gradient estimators, demonstrating that RU-MLMC and RR-MLMC are unbiased gradient estimators while the others are biased.

\begin{lemma}[Bias]
\label{lm:bias}
For any $x\in\RR^d$, it holds that%
\begin{align*}
    \EE v^{\LSGD}(x) = \EE  v^{\VMLMC}(x) = \EE v^{\RTMLMC}(x)=  \nabla F^L(x). 
\end{align*}
If additionally Assumption \ref{assumption:gradient_bias} holds, %
\begin{align*}
    \EE v^{\RUMLMC}(x) = \EE v^{\RRMLMC}(x) = \nabla F(x).
\end{align*}
\end{lemma}

\begin{proof}{Proof of Lemma~\ref{lm:bias}.}
The proof for $L$-SGD is straightforward. 
For V-MLMC, it holds that
\begin{equation*}
    \EE v^{\VMLMC}(x)  = \EE \bigg[\sum_{l=0}^L \frac{1}{n_l}\sum_{i=1}^{n_l} H^l(x,\zeta_i^l)\bigg] =  \sum_{l=0}^L  [\nabla F^l(x)-\nabla F^{l-1}(x)] = \nabla F^L(x).
\end{equation*}
For RT-MLMC, we have
\begin{align*}
        \EE v^{\RTMLMC}(x) 
    = 
        \EE \bigg[\frac{H^l(x,\zeta^l)}{q_l}\bigg]
    = 
        \sum_{l=0}^L q_l \frac{\nabla F^l(x)-\nabla F^{l-1}(x)}{q_l} = \nabla F^L(x).
\end{align*}
For RU-MLMC, letting $L\rightarrow\infty$, we have
\begin{align*}
        \EE v^{\RUMLMC}(x) 
    = 
        \EE \bigg[ \frac{H^l(x,\zeta^l)}{q_l} \bigg] 
    = 
        \sum_{l=0}^\infty q_l \frac{\nabla F^l(x)-\nabla F^{l-1}(x)}{q_l} 
    = 
        \lim_{L\rightarrow\infty}\nabla F^L(x) 
    =
        \nabla F(x),
\end{align*}
where the last equality uses Assumption \ref{assumption:gradient_bias}.
As for the RR-MLMC estimator, we have
\begin{align*}
    & 
        \EE v^{\RRMLMC}
    =  
        \EE   \sum_{l=0}^L p_l H^l(x,\zeta^l) 
    =
        \EE_L   \sum_{l=0}^L p_l [\nabla F^l(x) -\nabla F^{l-1}(x)]\\
    = & 
        \sum_{L=0}^\infty q_L\bigg( \sum_{l=0}^L p_l[\nabla F^l(x) -\nabla F^{l-1}(x)] \bigg)
    =  
        \sum_{L=0}^\infty \sum_{l=0}^L q_L p_l[\nabla F^l(x) -\nabla F^{l-1}(x)] \mathbb{I}\{l\leq L\}\\
    = & 
        \sum_{l=0}^\infty [\nabla F^l(x) -\nabla F^{l-1}(x)]\sum_{L=l}^\infty q_L p_l 
    =  
        \sum_{l=0}^\infty [\nabla F^l(x) -\nabla F^{l-1}(x)] \frac{\sum_{L=l}^\infty q_L}{1-\sum_{l^\prime=0}^{l-1} q_{l^\prime}} \\
    = & 
        \sum_{l=0}^\infty [\nabla F^l(x) -\nabla F^{l-1}(x)] 
    =  
        \lim_{L\rightarrow\infty}\nabla F^L(x) 
    = 
        \nabla F(x).
\end{align*}
By convention, we let $\sum_{l^\prime=0}^{-1} q_{l^\prime}=0$ and $\nabla F^{-1}(x) =0$. 
Similarly, the last equality is by Assumption~\ref{assumption:gradient_bias}.
\QED%
\end{proof}

In the following, we demonstrate the variance and per-iteration cost of these estimators. 
Recall that $\texttt{A}$ refers to either of $L$-SGD, V-MLMC, RT-MLMC, RU-MLMC, and RR-MLMC; 
$C_\mathrm{iter}^\texttt{A}$ denotes the (expected) per-iteration cost of $\texttt{A}$;
$\mathrm{Var}(v^\texttt{A})$ denotes the variance of the gradient estimator using $\texttt{A}$;
 $T^\texttt{A}$  denotes the iteration complexity of $\texttt{A}$ for achieving $\eps$-optimality or $\eps$-stationarity;
 and $C=T^\texttt{A}C_\mathrm{iter}^\texttt{A}$ denotes the (expected) total cost. 

Note that the (expected) per-iteration cost $C_\mathrm{iter}^\texttt{A}$ depends on different gradient constructions. By Lemma~\ref{lm:SGD}, the iteration complexity $T^\texttt{A}$ depends on the desired accuracy $\eps$ and $\mathrm{Var}(v^\texttt{A})$. 
To upper bound the (expected) total cost $C$, we first provide upper bounds on variance $\mathrm{Var}(v^\texttt{A})$ and (expected) per-iteration cost $C_\mathrm{iter}^\texttt{A}$.
The following summarizes Lemmas \ref{lm:variance_cost_LSGD}, \ref{lm:variance_cost_RT-MLMC}, \ref{lm:variance_cost_V-MLMC} and includes results for RU-MLMC/RR-MLMC estimators.

\begin{lemma}[Variance and (Expected) Per-Iteration Cost]
\label{lm:variance_cost}
Under Assumption \ref{assumption:general}, for any $x\in\RR^d$, the following results hold.
\begin{itemize}[leftmargin=2em]
\item For $v^{\LSGD}$ with batch size $n_L$, 
\begin{equation*}
        \mathrm{Var}\big(v^{\LSGD}(x)\big) 
    \leq 
        \frac{\sigma^2}{n_L},\qquad 
        C_\mathrm{iter}^{\LSGD}
    \leq
        n_L M_c 2^{cL}.     
\end{equation*}
\item For $v^{\RTMLMC}(x)$, suppose we take $q_l = 2^{-(b+c)l/2}R^L(-\frac{b+c}{2})^{-1}$, then
\[
\begin{array}{lll}
\text{If $c\ne b$: }&~~
 \displaystyle \mathrm{Var}(v^{\RTMLMC}(x) ) 
    \leq M_bR^L\big(\frac{c-b}{2}\big)R^L\big(-\frac{b+c}{2}\big),   
&~~ \displaystyle C_\mathrm{iter}^{\RTMLMC} \le M_cR^L\big(\frac{c-b}{2}\big)R^L\big(-\frac{b+c}{2}\big)^{-1}, \\
\text{If $c= b$: }&~~
 \displaystyle   \mathrm{Var}(v^{\RTMLMC}(x) ) 
    \leq M_b(L+1) R^L\big(-\frac{b+c}{2}\big),  & ~~ \displaystyle C_\mathrm{iter}^{\RTMLMC} \le M_c (L+1) R^L\big(-\frac{b+c}{2}\big)^{-1}.
\end{array}
\]
\item For $v^{\VMLMC}(x)$, suppose we take $n_l = \lceil  2^{-(b+c)l/2} N \rceil$ for a constant $N>0$, then
\[
\begin{array}{lll}
\text{If $c\ne b$: }&~~
 \displaystyle \mathrm{Var}(v^{\VMLMC}(x)) 
    \leq M_b R^L\big(\frac{c-b}{2}\big) N^{-1},   
&~~ \displaystyle C_\mathrm{iter}^{\VMLMC} \le M_cR^L\big(\frac{c-b}{2}\big)N + 
M_c R^L(c),\\ 
\text{If $c= b$: }&~~
 \displaystyle  \mathrm{Var}(v^{\VMLMC}(x)) 
    \leq M_b(L+1) N^{-1},  & ~~ \displaystyle C_\mathrm{iter}^{\VMLMC} \le M_c (L+1)N+M_c R^L(c).
\end{array}
\]
\item 
For $v^{\RUMLMC}(x)$ and $v^{\RRMLMC}(x)$, when $c\geq b$, either its expected per-iteration cost or its variance is unbounded for any distribution $Q=\{q_l\}_{l=0}^\infty$. When $c<b$, setting $q_l = 2^{-(b+c)l/2}R^\infty(-\frac{b+c}{2})^{-1}$, then
\[
\begin{array}{ll}
   \displaystyle\mathrm{Var}(v^{\RUMLMC}(x) ) 
    \leq
        M_bR^\infty\big(\frac{c-b}{2}\big)R^\infty\big(-\frac{b+c}{2}\big),  &\quad \displaystyle C_\mathrm{iter}^{\RUMLMC}\leq 
        M_c R^\infty\big(\frac{c-b}{2}\big)R^\infty\big(-\frac{c+b}{2}\big)^{-1},  \\
     \displaystyle\mathrm{Var}(v^{\RRMLMC}(x)) 
    \leq
        M_b R^\infty\big(\frac{c-b}{2}\big), & \quad \displaystyle C_\mathrm{iter}^{\RRMLMC}
    \leq 
        \frac{M_c2^c}{2^c-1}R^\infty\big(\frac{c-b}{2}\big)R^\infty\big(-\frac{c+b}{2}\big)^{-1}.  
\end{array}
\]
\end{itemize}
\end{lemma}

\begin{proof}{Proof of Lemma~\ref{lm:variance_cost}.}
Recall the definition of $R^L$ and $R^{\infty}$ in \eqref{Eq:RL:infty}, then it holds that
\begin{equation}
\label{eq:series_limit_1}
        R^L(\alpha)=\sum_{l=0}^L 2^{\alpha l}  
    = 
        \frac{1-2^{\alpha(L+1)}}{1-2^{\alpha}} \text{ for } \alpha \not = 0; \quad R^\infty(\alpha)=(1-2^{\alpha})^{-1} \text{ for } \alpha < 0.
\end{equation}
In the following, we show the variance and per-iteration cost for all gradient estimators.
\subsubsection*{$L$-SGD.} The result follows directly by Assumption~\ref{assumption:general} and the fact that $\{\nabla h^L(x,\zeta_i^L)\}_{i=1}^{n_L}$ are independent for any given $L$. 

\subsubsection*{V-MLMC.}
Notice that $\{H^l(x,\zeta_i^l)\}_{i=1}^{n_l}$ are independent for any given $l$ and $H^l(x,\zeta_i^l)$ with different $l$ are independently generated.
Using \eqref{eq:series_limit_1} and the fact that $n_l$ is an integer between $2^{-(b+c)l/2}N$ and $2^{-(b+c)l/2} N +1$, we derive
\begin{equation*}
    \mathrm{Var}(v^{\VMLMC}(x)) 
    \leq 
        \sum_{l=0}^L \frac{V_l}{n_l} 
    \leq  
        M_b\left[\sum_{l=0}^L \frac{2^{-bl}}{2^{-(b+c)l/2}}\right]N^{-1}
\end{equation*}
and
\begin{equation*}
\begin{split}
        C_\mathrm{iter}^{\VMLMC}
    = &
        \sum_{l=0}^L n_l C_l 
    \leq  
        M_c\left[\sum_{l=0}^L2^{cl}2^{-(b+c)l/2}\right]N+M_c\sum_{l=0}^L2^{cl}.
\end{split}
\end{equation*}
Moreover, one can simplify these expressions by considering cases $c=b$ or $c\ne b$.
\subsubsection*{RT-MLMC.} 
The result follows by simplifying
\[
\mathrm{Var}(v^{\RTMLMC}(x))\le \sum_{l=0}^L\frac{V_l}{q_l},\quad 
C_\mathrm{iter}^{\RTMLMC}=\sum_{l=0}^LC_lq_l
\]
and using the relation \eqref{eq:series_limit_1}.

\subsubsection*{RU-MLMC.}
By direct calculation,
\begin{equation}
\mathrm{Var}(v^{\RUMLMC}(x))\le \sum_{l=0}^{\infty}\frac{V_l}{q_l},\quad 
C_\mathrm{iter}^{\RUMLMC}=\sum_{l=0}^{\infty}C_lq_l.\label{Eq:var:C:RUMLMC}
\end{equation}
To ensure $\mathrm{Var}(v^{\RUMLMC}(x))<\infty$, it requires $q_l > 2^{-bl}$.
To ensure $C_\mathrm{iter}^{\RUMLMC}<\infty$, it requires $q_l<2^{-cl}$.
Combining those relations, we find the per-iteration cost, and its variance is bounded only if $b>c$.
In such case, we substitute the expression of $q_l, V_l, C_l$ into \eqref{Eq:var:C:RUMLMC} and utilize \eqref{eq:series_limit_1} to obtain the desired result.

\subsubsection*{RR-MLMC.} 
The expected cost to generate such an estimator is 
\begin{equation}
\begin{aligned}
        C_\mathrm{iter}^{\RRMLMC}
    = &
        \EE_L \sum_{l=0}^L C_l = \sum_{L=0}^\infty q_L \sum_{l=0}^LC_l
    \leq 
        \sum_{L=0}^\infty q_L \sum_{l=0}^L M_c 2^{cl} 
    = 
        M_c \sum_{L=0}^\infty q_L \frac{2^{c(L+1)}-1}{2^{c}-1}.
\end{aligned}\label{Eq:RRMLMC:C}
\end{equation}
To ensure $C_\mathrm{iter}^{\RRMLMC}$ is finite, it requires $q_L < 2^{-cL}$.
Without loss of generality, we assume that $q_l = 2^{\alpha l}(1-2^\alpha)$ for a constant $0<\alpha<-c$ and $l\in\NN$ so that $\sum_{l=0}^\infty q_l = 1$. As a result, we have $\sum_{l^\prime=0}^{l-1} q_{l^\prime} = 1-2^{\alpha l}$.
The variance of the estimator is 
\begin{equation}
\begin{multlined}
        \mathrm{Var}(v^{\RRMLMC}(x)) 
    \leq 
        \sum_{L=0}^\infty q_L \Big(\sum_{l=0}^L p_l^2 V_l\Big)
    =
        \sum_{L=0}^\infty q_L \Big(\sum_{l=0}^L \frac{1}{(1-\sum_{l^\prime=0}^{l-1} q_{l^\prime})^2}\mathbb{I}\{l\leq L\} V_l\Big) \\
    = 
        \sum_{L=0}^\infty  \sum_{l=0}^L q_L \frac{1}{(1-\sum_{l^\prime=0}^{l-1} q_{l^\prime})^2}\mathbb{I}\{l\leq L\} V_l
    = 
        \sum_{l=0}^\infty V_l  \frac{\sum_{L=l}^\infty   q_L}{(1-\sum_{l^\prime=0}^{l-1} q_{l^\prime})^2}
    = 
        \sum_{l=0}^\infty   \frac{V_l}{1-\sum_{l^\prime=0}^{l-1} q_{l^\prime}}. 
\end{multlined}
\label{Eq:var:RRMLMC}
\end{equation}
For the relation above, we abuse the use of notation to denote $\sum_{l^\prime=0}^{l-1}q_{l^\prime}=0$ when $l=0$.
To obtain a bounded variance, it requires that 
\begin{align*}
        \frac{2^{-bl}}{1-\sum_{l^\prime=0}^{l-1} q_{l^\prime}} < 1,
\end{align*}
which equivalently requires $\alpha >-b$.
Therefore, to ensure both the expected per-iteration cost and the variance are finite, it requires $c<b$. 
In this case, we substitute $q_l\propto 2^{-(b+c)l/2}$ into \eqref{Eq:RRMLMC:C} and \eqref{Eq:var:RRMLMC} to obtain the desired results.
\QED
\end{proof}

\section{Total Cost Analysis}
\label{app_sec:total_cost}
This section discusses the (expected) total cost of  $L$-SGD, V-MLMC, RT-MLMC, RU-MLMC, and RR-MLMC with SGD framework in Algorithm \ref{alg:SGD-framework}. 

For the convex objective, we only show the proof of the smooth case, which essentially relies on Lemma \ref{lm:SGD}\ref{Lemma:EC1:b}.
We find the convergence results in Lemma \ref{lm:SGD}\ref{Lemma:EC1:b}~(smooth case) only differ from Lemma \ref{lm:SGD}\ref{Lemma:EC1:d}~(Lipscthiz continuous case) up to some constants.
Thereby, one can analyze the convex and Lipscthiz continuous case by replacing \ref{lm:SGD}\ref{Lemma:EC1:b} with  Lemma \ref{lm:SGD}\ref{Lemma:EC1:d} to obtain similar results.
We will see that V-MLMC requires a large mini-batch to achieve a reduced total cost in a smooth setting. For the Lipschitz continuous case, V-MLMC cannot achieve a reduced total cost compared to $L$-SGD.

\subsection{(Expected) Total Costs of Biased Gradient Methods }\label{Appendix:EC:4.1}
In this subsection, let $\mathtt{A}$ refer to the $L$-SGD, V-MLMC, or RT-MLMC estimator.
Recall that $\mathbf{V}(v^{\texttt{A}})$ denotes the upper bound on the variance of the gradient estimator $\mathtt{A}$, and  $\hat{x}_T^{\texttt{A}}$ denotes the output uniformly selected from $\{x_1^{\texttt{A}},...,x_T^{\texttt{A}}\}$, where $x_t^{\texttt{A}}$ refers to the $t$-iteration point of SGD framework using estiamtor $\texttt{A}$.
Sometimes, we omit the dependence on the notation $\texttt{A}$ for simplicity of presentation.

In the (strongly) convex case, recall from the decomposition \eqref{eq:error_decomposition} that we shall take $L = \lceil 1/a \log(4M_a/\eps)  \rceil$.
It remains to apply Lemma~\ref{lm:SGD} in Appendix \ref{app_sec:SGD} to show $\EE F^L(\hat x^{\texttt{A}}_T) -  F^L(x^L)\leq \eps/2$. 
In the nonconvex smooth case,  by Assumption~\ref{assumption:gradient_bias}, we have the decomposition
\begin{equation}
\label{eq:truncate_nonconvex_error}
\begin{aligned}
        \EE \|\nabla F(\hat x^{\texttt{A}}_T)\|_2^2
    \leq & 
        2\EE \|\nabla F^L(\hat x^{\texttt{A}}_T)\|_2^2+ 2\EE \|\nabla F(\hat x^{\texttt{A}}_T)-\nabla F^L(\hat x^{\texttt{A}}_T)\|_2^2
    \leq  
        2\EE \|\nabla F^L(\hat x^{\texttt{A}}_T)\|_2^2 + 2M_a 2^{-aL}.
\end{aligned}   
\end{equation}
Therefore, we take $L = \lceil 1/a \log(4M_a/\eps^2)\rceil$ so that $2M_a 2^{-aL}\leq \eps^2/2$.
To ensure $\EE \|\nabla F(\hat x^{\texttt{A}}_T)\|_2^2\leq \eps^2$, it remains to apply Lemma~\ref{lm:SGD}\ref{Lemma:EC1:c} to show $\EE \|\nabla F^L(\hat x^{\texttt{A}}_T)\|_2^2\leq \eps^2/4$.

When the variance of the gradient estimator is  $\cO(1)$ or even larger, the following lemma shows the relationship between the variance of the gradient estimator, the (expected) cost to construct the gradient estimator, and the (expected) total cost for achieving $\eps$-optimality in the (strongly) convex case and $\eps$-stationarity in the nonconvex smooth case, respectively.

\begin{lemma}
\label{thm:biased}
Suppose $\mathbf{V}(v^{\texttt{A}})\ge \rho_0$, with $\rho_0>0$ being a constant independent of the precision $\eps>0$.
With properly selected hyper-parameters as in Table~\ref{tab:results:Lemma:EC4}, the following results hold.
\begin{enumerate}
    \item\label{thm:biased:a}
Suppose Assumption~\ref{assumption:general} holds and $F^L$ is $\mu$-strongly convex, and let $x_T^{\texttt{A}}$ be the estimated solution.
The (expected) total cost of $\texttt{A}$ for finding an $\eps$-optimal solution satisfies 
\[
C\le 4C_\mathrm{iter}^{\texttt{A}} S_F\max\{\mathbf{V}(v^{\texttt{A}}), \mu^2(1+2S_F^2/\mu^2) \|x_1-x^L\|_2^2\}\mu^{-2} \eps^{-1}.
\]
    \item\label{thm:biased:b}
Suppose Assumption~\ref{assumption:general} holds, and $F^L$ is convex and smooth. 
Let $\hat{x}_T^{\texttt{A}}$ be the estimated solution.
The (expected) total cost of $\texttt{A}$ for finding an $\eps$-optimal solution satisfies 
    \begin{align*}
    C \leq 8C_\mathrm{iter}^{\texttt{A}}\mathbf{V}(v^{\texttt{A}})(1+\|x_1-x^L\|_2^2)^2\eps^{-2}.
    \end{align*} 
    \item\label{thm:biased:c}
Suppose Assumptions \ref{assumption:general}\ref{assumption:general:II}\ref{assumption:general:III} and \ref{assumption:gradient_bias} hold and $F^L$ is $S_F$-smooth. Let $\hat{x}_T^{\texttt{A}}$ be the estimated solution.
The (expected) total cost of $\texttt{A}$ for finding an $\eps$-stationary point satisfies
\[
C \leq 32C_\mathrm{iter}^{\texttt{A}}\mathbf{V}(v^{\texttt{A}})(2(F^L(x_1)-F^L(x^L))+S_F)^2\eps^{-4}.
\]
\end{enumerate}

\begin{table}[!ht]
\caption{Hyper-parameters used in Lemma~\ref{thm:biased}.
}
\centering
\small{
\begin{tabular}{c|c|c|c}
	\toprule%
& \textbf{Level } $L$ & 
\textbf{Step Size } $\gamma_t$ & \textbf{Iteration } $T$
\\[8pt]
\hline
\mbox{Case~\ref{thm:biased:a}}&$\lceil 1/a \log(4M_a/\eps)  \rceil$&$\frac{1}{\mu(t+2S_F^2/\mu^2)}$& $\left\lceil 
\frac{2S_F\max\{\mathbf{V}(v^{\texttt{A}}), \mu^2(1+2S_F^2/\mu^2) \|x_1-x^L\|_2^2\}}{\mu^2 \eps}
\right\rceil$
\\[10pt]\hline
\mbox{Case~\ref{thm:biased:b}}& $\lceil 1/a \log(4M_a/\eps)  \rceil$ & $\frac{1}{\sqrt{T\mathbf{V}(v^{\texttt{A}})}}$ & $\left\lceil 
4\mathbf{V}(v^{\texttt{A}})(1+\|x_1-x^L\|_2^2)^2\eps^{-2}
\right\rceil$
\\[10pt]\hline
\mbox{Case~\ref{thm:biased:c}}& $\lceil 1/a \log(4M_a/\eps^2)  \rceil$&$\frac{1}{\sqrt{T\mathbf{V}(v^{\texttt{A}})}}$&$\left\lceil 
16\mathbf{V}(v^{\texttt{A}})(2(F^L(x_1)-F^L(x^L))+S_F)^2\eps^{-4}
\right\rceil$\\[10pt]
\bottomrule
\end{tabular}}
\label{tab:results:Lemma:EC4}
\end{table}
\end{lemma}

This lemma applies to $L$-SGD and RT-MLMC since the variance of these estimators are both $\cO(1)$ or higher.
On the contrary, V-MLMC has to use large mini-batches. Thus its variance is very small and is $\tilde{\cO}(\eps)$ or $\tilde{\cO}(\eps^2)$.
When showing Lemma~\ref{thm:biased}, we implicitly assumes the step size $\gamma_t\le 1/(2S_F)$ in case \ref{thm:biased:b} and $\gamma_t\le 1/S_F$ in case \ref{thm:biased:c}. 
Those requirements automatically hold since the precision $\eps$ is sufficiently small.

\begin{proof}{Proof of Lemma~\ref{thm:biased}}
We prove this lemma in three separate cases.
\subsubsection*{Strongly convex case.}
By Lemma \ref{lm:SGD}\ref{Lemma:EC1:a2}, it holds that
\begin{align*}
        \EE[F^L(x_T^{\texttt{A}})-F^L(x^L)]
    \leq 
        \frac{S_F\max\{\mathbf{V}(v^{\texttt{A}}), \mu^2(1+2S_F^2/\mu^2) \|x_1-x^L\|_2^2\}}{\mu^2 (T+2S_F^2/\mu^2)}.
\end{align*}
Based on the previous argument, it remains to ensure $\EE[F^L(x_T^A)-F^L(x^L)]\le \eps/2$.
Thus, we take $T$ to be the smallest integer such that
\[
T\ge \hat{T}:=2S_F\max\{V(v^A), \mu^2(1+2S_F^2/\mu^2) \|x_1-x^L\|_2^2\}\mu^{-2}\eps^{-1}.
\]
Correspondingly, the total cost $C=TC_\mathrm{iter}^{\texttt{A}}\le 2\hat{T}C_\mathrm{iter}^{\texttt{A}}$.
    
\subsubsection*{Convex case.}
By Lemma~\ref{lm:SGD}\ref{Lemma:EC1:b}, it holds that 
\begin{align*}
        \EE[F^L(x_T^{\texttt{A}})-F^L(x^L)]
    \leq  
        \frac{\sqrt{\mathbf{V}(v^{\texttt{A}})}(1+\|x_1-x^L\|_2^2)}{\sqrt{T}}.
\end{align*}
To ensure $\EE[F^L(x_T^{\texttt{A}})-F^L(x^L)]\le \eps/2$, we take $T$ to be the smallest integer such that
\[
T\ge \hat{T}:=4\mathbf{V}(v^{\texttt{A}})(1+\|x_1-x^L\|_2^2)^2\eps^{-2}.
\]
Correspondingly, the total cost $C=TC_\mathrm{iter}^{\texttt{A}}\le 2\hat{T}C_\mathrm{iter}^{\texttt{A}}$.
    
\subsubsection*{Nonconvex smooth case.}
By Lemma~\ref{lm:SGD}\ref{Lemma:EC1:c}, it holds that
\begin{align*}
\EE \|\nabla F^L(\hat x_T^{\texttt{A}})\|_2^2 \leq \frac{\sqrt{\mathbf{V}(v^{\texttt{A}})}(2(F^L(x_1)- F^L(x^L))+S_F)}{\sqrt{T}}.
\end{align*}
Based on the previous argument, it remains to ensure $\EE \|\nabla F^L(\hat x_T^{\texttt{A}})\|_2^2\le \eps^2/4$.
Thus, we take $T$ to be the smallest integer such that
\[
T\ge \hat{T}:=16\mathbf{V}(v^{\texttt{A}})(2(F^L(x_1)-F^L(x^L))+S_F)^2\eps^{-4}.
\]
Correspondingly, the total cost $C=TC_\mathrm{iter}^{\texttt{A}}\le 2\hat{T}C_\mathrm{iter}^{\texttt{A}}$.\QED
\end{proof}

\subsubsection{Total Cost of \texorpdfstring{$L$-SGD}{}}
\label{appendix:L-SGD}
In this subsection, we show the total cost of $L$-SGD when the batch size $n_L=1$ using Lemma~\ref{thm:biased}.
For (strongly) convex case, we take $L = \lceil 1/a \log(4M_a/\eps)  \rceil$, which implies 
\[
C_\mathrm{iter}^{\LSGD}\le M_c 2^{cL} = \cO(\epsilon^{-c/a}).
\]
For nonconvex case, we take $L = \lceil 1/a \log(4M_a/\eps^2) \rceil$, which implies $C_\mathrm{iter}^{\LSGD}=\cO(\eps^{-2c/a})$.
The total cost in Theorem~\ref{thm:L-SGD} can be derived by substituting the expressions of $C_\mathrm{iter}^{\LSGD}$ and $\mathbf{V}(v^{\LSGD})$ into Lemma~\ref{lm:variance_cost}.

\begin{remark}
As mentioned in Remark~\ref{rem:V-MLMC}, $L$-SGD cannot use large mini-batch sizes, i.e., $n_L=\cO(\eps^{-1})$, to reduce the total cost as V-MLMC does. We use the convex case for demonstration. If $n_L=\cO(\eps^{-1})$, then $V(v^{L\mathrm{-SGD}}) = \cO(\eps)$. Using Lemma \ref{lm:SGD}\ref{Lemma:EC1:c} with stepsizes $\gamma_t = 1/(2S_F)$, we know that 
\begin{align*}
\EE F^L(\hat x_T^{\LSGD})-F^L(x^L) \leq  \cO(\eps) + \frac{2S_F\|x_1-x^L\|_2^2}{T}.   
\end{align*}
Therefore the iteration complexity $T=\cO(\eps^{-1})$. The total cost would be $Tn_LC_L=\cO(\eps^{-2-c/a})$, which remains the same compared to without using batch size, i.e., $n_L=1$. The reason is that using large batch sizes for $L$-SGD only reduces the iteration complexity, while using large batch sizes for V-MLMC simultaneously reduces the iteration complexity and the per-iteration cost.
\end{remark}

\subsubsection{Expected Total Cost of RT-MLMC}
\label{appendix:RT-MLMC}

According to Lemma~\ref{lm:variance_cost_RT-MLMC}, one can check that $\textbf{V}(v^{\RTMLMC})$ is at least $\cO(1)$, indicating Lemma~\ref{thm:biased} is applicable.
After calculating upper bounds of $C_\mathrm{iter}^{\RTMLMC}$ and $V(v^{\RTMLMC})$ by specifying $L$ in Table~\ref{tab:assumptions:RTMLMC}, it holds that
\begin{equation}
\label{eq:upper_bound_RTMLMC}
C_\mathrm{iter}^{\RTMLMC}\textbf{V}(v^{\RTMLMC}) = 
\begin{cases}
\cO(1) & \text{ if } c<b;\\[10pt]
\cO(\log^2(1/\eps))  &
\text{ if } c=b;\\[10pt]
\cO(\eps^{-(c-b)/a})&
\text{ if } c>b\text{ and }L=\lceil \frac{\log(4M_a\eps^{-1})}{a}  \rceil;\\[10pt]
\cO(\eps^{-2(c-b)/a})&
\text{ if } c>b\text{ and }L=\lceil \frac{\log(4M_a\eps^{-2})}{a}  \rceil.\\[10pt]
\end{cases}   
\end{equation}
Hence, the total cost in Theorem~\ref{thm:RT-MLMC} can be derived by substituting the expressions of $C_\mathrm{iter}^{\RTMLMC}\textbf{V}(v^{\RTMLMC})$ into Lemma~\ref{lm:variance_cost}.

\subsubsection{Total Cost of V-MLMC}
\label{appendix:V_MLMC}
In this part, we consider the total cost of V-MLMC.

\begin{theorem}[Full version of Theorem \ref{thm:V-MLMC}]
\label{thm:V-MLMC-full}
For V-MLMC with $n_l = \lceil 2^{-(b+c)l/2} N\rceil$ for some $N>0$, with properly chosen hyper-parameters as in Table~\ref{tab:assumptions:VMLMC:appendix}, the following results hold. %
\begin{enumerate}
    \item\label{thm:V-MLMC-full:a}
Suppose $F^L$ is $\mu$-strongly convex and $S_F$-smooth, and Assumption~\ref{assumption:general} holds, and let $x_T^{\VMLMC}$ as the estimated solution.
The total cost for finding $\eps$-optimal solution of $F$ is
\[
\begin{cases}
    \cO(\log(\eps^{-1})\eps^{-\max\{1,c/a\}}), & \text{ if } c<b,\\[10pt]
    \cO(\log^3(\eps^{-1})\eps^{-\max\{1,c/a\}}), & \text{ if } c=b,\\[10pt]
    \cO(\log(\eps^{-1})\eps^{-\max\{1+{(c-b)/a},c/a\}}), & \text{ if } c> b.
    \end{cases}
\] 
     \item\label{thm:V-MLMC-full:b}
Suppose $F^L$ is convex and $S_F$-smooth, and Assumption~\ref{assumption:general} holds, and let $\hat{x}_T^{\VMLMC}$ as the estimated solution. 
The total cost for finding $\eps$-optimal solution of $F$ is
\[
\begin{cases}
        \cO(\eps^{-1-\max\{1,c/a\}}),  & \text{ if } c<b,\\[10pt]
        \cO(\log^2(\eps^{-1})\eps^{-1-\max\{1,c/a\}}),   & \text{ if } c=b,\\[10pt]
        \cO(\eps^{-1-\max\{1+(c-b)/a,c/a\}}), & \text{ if } c>b.
    \end{cases}
\] 
     \item\label{thm:V-MLMC-full:c}
Suppose $F^L$ is $S_F$-smooth, and Assumptions \ref{assumption:general}\ref{assumption:general}\ref{assumption:general:II}\ref{assumption:general:III} and \ref{assumption:gradient_bias} hold.
Let $\hat{x}_T^{\VMLMC}$ as the estimated solution, then the total cost for finding $\eps$-stationary point of $F$ is
\[
\begin{cases}
    \cO(\eps^{-2-2\max\{1,c/a\}}  ),& \text{ if } c<b,\\[10pt]
    \cO(\log^2(\eps^{-1})\eps^{-2-2\max\{1, c/a\}}),   & \text{ if } c=b,\\[10pt]
    \cO(\eps^{-2-2\max\{1+(c-b)/a,c/a\}} ),& \text{ if } c>b.
    \end{cases}
\]
\end{enumerate}

\begin{table}[!ht]
\caption{Hyper-parameters of V-MLMC in Theorem~\ref{thm:V-MLMC}}
\centering
\small{
\begin{tabular}{c|c|c|c|c}
	\toprule%
 & \bf{Level $L$} & \bf{Step Size $\gamma_t\equiv \gamma$} & $\begin{array}{c}
    \textbf{Multiplicative Factor}\\
    \textbf{$N$ of Batch Sizes}
\end{array}$& \bf{Iteration $T$} \\[8pt]
\hline
\mbox{Case~\ref{thm:V-MLMC-full:a}}
& $\lceil \frac{\log(4M_a\eps^{-1})}{a}  \rceil$&  $
1/S_F
$
& 
$\begin{cases}
    \cO(\eps^{-1}), &    \text{ if } c<b \\[10pt]
    \cO((L+1)\eps^{-1}), &    \text{ if } c=b \\[10pt]
    \cO(2^{-(b-c)(L+1)/2}\eps^{-1}), &    \text{ if } c>b
    \end{cases}$
&
$\Bigg\lceil 
\frac{2\log(4[F^L(x_1)-F^L(x^L)]\eps^{-1})}{\log(S_F/(S_F-\mu))}
\Bigg\rceil$
\\[10pt]
\hline
\mbox{Case~\ref{thm:V-MLMC-full:b}} &$\lceil \frac{\log(4M_a\eps^{-1})}{a}  \rceil$ &  $1/(2S_F)$& 
$\begin{cases}
        \cO(\eps^{-1}), & \text{ if } c<b\\[10pt]
        \cO((L+1)\eps^{-1}),& \text{ if } c=b\\[10pt]
        \cO(2^{-(b-c)(L+1)/2}\eps^{-1}), & \text{ if } c>b
    \end{cases}$
&
$\Bigg\lceil8S_F\|x_1-x^L\|_2^2\eps^{-1} \Bigg\rceil$
\\[10pt]
\hline
\mbox{Case~\ref{thm:V-MLMC-full:c}} &$\lceil \frac{\log(4M_a\eps^{-2})}{a}  \rceil$ &  $1/S_F$& 
$\begin{cases}
        \cO(\eps^{-2}), & \text{ if } c<b\\[10pt]
        \cO((L+1)\eps^{-2}),& \text{ if } c=b\\[10pt]
        \cO(2^{-(b-c)(L+1)/2}\eps^{-2}),& \text{ if } c>b
    \end{cases}$
&
$\Bigg\lceil 16S_F(F^L(x_1)-F^L(x^L))\eps^{-2} \Bigg\rceil$
\\[8pt]
\bottomrule
\end{tabular}}
\label{tab:assumptions:VMLMC:appendix}%
\end{table}

\end{theorem}
\begin{remark}
When $a\geq \min\{b,c\}$, it holds that $\max\{1, c/a\}=1$ when $c<b$, and $\max\{1+(c-b)/a,c/a\}=1$ when $c>b$.
As a result, when this additional condition holds, the total cost of V-MLMC is the same as that of RT-MLMC. 
\end{remark}

\begin{proof}{Proof of Theorem~\ref{thm:V-MLMC}.}
Recall the discussion from the beginning of Appendix~\ref{Appendix:EC:4.1}, we need to specify $L=\lceil \frac{\log(4M_a\eps^{-1})}{a}  \rceil$ in (strongly) convex case and $L = \lceil \frac{\log(4M_a\eps^{-2})}{a}  \rceil$ in nonconvex case.

\subsubsection*{Strongly convex case.} By Lemma \ref{lm:SGD}\ref{Lemma:EC1:a1}, it holds that 
\begin{align*}
    \EE F^L(x_T) -F^L(x^L) \leq (1-\gamma \mu)^{T-1}[F^L(x_1)-F^L(x^L)]+\frac{S_F\gamma \mathbf{V}(v^\mathrm{\VMLMC})}{2\mu}.
\end{align*}
To ensure $\EE F^L(x_T) -F^L(x^L)\le \eps/2$, it suffices to make sure that 
$(1-\gamma \mu)^{T-1}[F^L(x_1)-F^L(x^L)]\leq \eps/4$ and $\frac{S_F\gamma \mathbf{V}(v^\mathrm{\VMLMC})}{2\mu}\le \eps/4$.
We take step size $\gamma=1/S_F$, then it suffices to ensure
\[
\mathbf{V}(v^\mathrm{\VMLMC})\le \frac{\eps\mu}{2},\quad 
T\ge \hat{T}:=
\frac{2\log(4[F^L(x_1)-F^L(x^L)]\eps^{-1})}{\log(S_F/(S_F-\mu))}.
\]
Since $T\in\mathbb{N}$, we specify $T=\lceil \hat{T}\rceil$.
Now it remains to specify $N$ such that the variance bound provided in Lemma~\ref{lm:variance_cost} satisfies the requirement, and one can also bound the per-iteration cost by Lemma~\ref{lm:variance_cost}.
\begin{itemize}
    \item 
When $c<b$, it holds that
\[
\mathbf{V}(v^\mathrm{\VMLMC}) = M_b\frac{1-2^{-(b-c)(L+1)/2}}{1-2^{-(b-c)/2}}N^{-1}\le \frac{\eps\mu}{2},
\]
and it suffices to take $N = 2M_b(1-2^{-(b-c)/2})^{-1}\eps^{-1}\mu^{-1}=\cO(\eps^{-1})$. 
In this case, the per-iteration cost $C_{\mathrm{iter}}^{\VMLMC}=\cO(\eps^{-\max\{1, c/a\}})$ and total cost 
\[
C\le TC_{\mathrm{iter}}^{\VMLMC}=\cO(\log(\eps^{-1})\eps^{-\max\{1, c/a\}}).
\]
    \item
When $c=b$, it suffices to take $N = 2M_b(L+1)\eps^{-1}\mu^{-1}=\cO(\log(\eps^{-1})\eps^{-1})$ to ensure $\mathbf{V}(v^\mathrm{\VMLMC})\le \frac{\eps\mu}{2}$.
Then, the per-iteration cost $C_{\mathrm{iter}}^{\VMLMC}=\cO(\log^2(\eps^{-1})\eps^{-1} + \eps^{-c/a})$ and total cost 
\[
C\le TC_{\mathrm{iter}}^{\VMLMC}=\cO(\log^3(\eps^{-1})\eps^{-1} + \log(\eps^{-1})\eps^{-c/a})=\cO(\log^3(\eps^{-1}) \eps^{-\max\{1, c/a\}}).
\]
    \item
When $c>b$, it suffices to take $N = 2M_b2^{-(b-c)(L+1)/2}(2^{-(b-c)/2}-1)^{-1}\mu^{-1}\eps^{-1}=\cO(\eps^{-1}\eps^{-(c-b)/(2a)})$ to ensure $\mathbf{V}(v^\mathrm{\VMLMC})\le \frac{\eps\mu}{2}$.
Then, the per-iteration cost $C_{\mathrm{iter}}^{\VMLMC}=\cO(\eps^{-1-(c-b)/a} + \eps^{-c/a})$ and total cost
\[
C\le TC_{\mathrm{iter}}^{\VMLMC}=\cO(\log(\eps^{-1})\eps^{-\max\{1 + (c-b)/a, c/a\}}).
\]
\end{itemize}

\subsubsection*{Convex case.} By Lemma \ref{lm:SGD}\ref{Lemma:EC1:b}, with stepsizes $\gamma_t \equiv \gamma = \frac{1}{2S_F}$, it holds that 
\begin{align*}
        \EE [F^L(\hat x_T)- F^L(x^L)] 
    \leq  
        \frac{\mathbf{V}(v^\mathrm{\VMLMC})}{2S_F}+\frac{2S_F\|x_1-x^L\|_2^2}{T}.
\end{align*}
To ensure $\EE [F^L(\hat x_T)- F^L(x^L)] \le \eps/2$, it suffices to take 
\begin{align*}
    T\geq \hat{T}:=8S_F\|x_1-x^L\|_2^2\eps^{-1};\qquad 
    \mathbf{V}(v^\mathrm{\VMLMC})\leq \frac{\eps S_F}{2}.
\end{align*}
We specify $T=\lceil \hat{T}\rceil$. 
By a similar argument as in the strongly convex case, we specify
\begin{align*}
N =
\begin{cases}
2M_b(1-2^{-(b-c)/2})^{-1}\eps^{-1}S_F^{-1}=\cO(\eps^{-1}), & \text{ if } c<b,\\[10pt]
2M_b(L+1)\eps^{-1}S_F^{-1}=\cO(\log(\eps^{-1})\eps^{-1}),& \text{ if } c=b,\\[10pt]
2M_b2^{-(b-c)(L+1)/2}(2^{-(b-c)/2}-1)^{-1}\eps^{-1}S_F^{-1}=\cO(\eps^{-1-(c-b)/(2a)}),& \text{ if } c>b.
\end{cases}
\end{align*}
Consequently, $C_\mathrm{iter}^\mathrm{V-MLMC}$ is bounded by
\[
\begin{array}{llllll}
\cO(\eps^{-\max\{1, c/a\}}),    & \quad\text{if $c<b$;} & \qquad \cO(\eps^{-\max\{1, c/a\}})\log^2(\eps^{-1})), & \quad\text{if $c=b$;}
 & \qquad\cO( \eps^{-\max\{1+ (c-b)/a, c/a\}} ), &\quad \text{if $c>b$.}
\end{array}
\] 
and one can bound the total cost by $TC_\mathrm{iter}^\mathrm{V-MLMC}$ to obtain the desired result.

\subsubsection*{Nonconvex smooth case.} By Lemma \ref{lm:SGD}\ref{Lemma:EC1:d}, with stepsizes $\gamma_t \equiv \gamma = \frac{1}{S_F}$, it holds that
\begin{align*}
\EE \|\nabla F^L(\hat x_T)\|_2^2 \leq \frac{2S_F(F^L(x_1)-F^L(x^L))}{T}+\mathbf{V}(v^\mathrm{\VMLMC})\leq \frac{\eps^2}{4}.
\end{align*}
To ensure $\EE \|\nabla F^L(\hat x_T)\|_2^2\le \frac{\eps^2}{4}$, 
it suffices to have
\begin{align*}
T\geq \hat{T}:=16S_F(F^L(x_1)-F^L(x^L))\eps^{-2},\qquad 
\mathbf{V}(v^\mathrm{\VMLMC}) \leq \frac{\eps^2}{8}.
\end{align*}
We specify $T=\lceil \hat{T}\rceil$. 
By a similar argument as in the strongly convex case, we specify
\begin{align*}
N = \begin{cases}
8M_b(1-2^{-(b-c)/2})^{-1}\eps^{-2}=\cO(\eps^{-2}), & \text{ if } c<b,\\[10pt]
8M_b(L+1)\eps^{-2}=\cO(\eps^{-2}\log(\eps^{-1})),& \text{ if } c=b,\\[10pt]
8M_b2^{-(b-c)(L+1)/2}(2^{-(b-c)/2}-1)^{-1}\eps^{-2}=\cO(\eps^{-2-(b-c)/(a)}),& \text{ if } c>b.\\
\end{cases}
\end{align*}
Consequently, $C_\mathrm{iter}^{\VMLMC}$ is bounded by
\[
\begin{array}{llllll}
\cO(\eps^{-2\max\{1, c/a\}}),    & \quad\text{if $c<b$;} & \quad \cO(\eps^{-2\max\{1, c/a\}})\log^2(\eps^{-1})), & \quad\text{if $c=b$;}
 & \quad\cO( \eps^{-2\max\{1+ (c-b)/a, c/a\}} ), &\quad \text{if $c>b$.}
\end{array}
\] 
and one can bound the total cost by $TC_\mathrm{iter}^{\VMLMC}$ to obtain the desired result.

\QED
\end{proof}

\subsection{Expected Total Cost of Unbiased Gradient Methods: RU-MLMC and RR-MLMC}
\label{appendix:unbiased}

In this subsection, we derive the expected total cost of RU-MLMC (see \eqref{alg:RU-MLMC}) and RR-MLMC (see \eqref{alg:RR-MLMC}). 
The upper bounds on the variance and per-iteration cost of these two estimators are in Lemma~\ref{lm:variance_cost}.
Both estimators are unbiased and are applicable only when $c<b$, namely, the increase rate of the cost to generate a gradient estimator per level is strictly smaller than the decrease rate of the variance of $\nabla H^l(x,\zeta^l)$.

\begin{proof}{Proof of Theorem~\ref{thm:unbiased_MLMC}.}
Since the gradient estimator $\texttt{A}\in\{\text{RU-MLMC, RR-MLMC}\}$ is unbiased, Lemma \ref{lm:SGD} can be applied directly.
\subsubsection*{Strongly convex case.} By Lemma \ref{lm:SGD}\ref{Lemma:EC1:a2}, to ensure that $x_T^{\texttt{A}}$ is an $\eps$-optimal solution, it suffices to take $T$ as the smallest integer that satisfies 
\begin{align*}
T \geq \hat{T}:=S_F\max\{\mathbf{V}(v^{\texttt{A}}),\mu^2(1+2S_F^2/\mu^2)\|x_1-x^*\|_2^2)\}\mu^{-2}\eps^{-1}.
\end{align*}
Therefore, the expected total cost 
\[
C = T^* C_\mathrm{iter}^{\texttt{A}}\le 2\hat{T}C_\mathrm{iter}^{\texttt{A}}=\cO(\eps^{-1}),
\]
where the last equality is based on Lemma~\ref{lm:variance_cost} so that $C_\mathrm{iter}^{\texttt{A}}=\cO(1)$ and $\mathbf{V}(v^{\texttt{A}})=\cO(1)$.

\subsubsection*{Convex case.} 
To apply Lemma \ref{lm:SGD}\ref{Lemma:EC1:b}, it requires that $\gamma_t = \frac{1}{\sqrt{\mathbf{V}(v^\texttt{A})T}}\leq \frac{1}{2S_F}.$
Recall that $\mathbf{V}(v^\texttt{A})=\cO(1)$, and therefore the requirement is satisfied for large enough $T$.
By Lemma \ref{lm:SGD}\ref{Lemma:EC1:b}, to ensure that $\hat{x}_T^{\texttt{A}}$ is an $\eps$-optimal solution, it suffices to take $T$ as the smallest integer that satisfies 
\begin{align*}
T \geq \hat{T}:=\eps^{-2} \mathbf{V}(v^{\texttt{A}})(1+\|x_1-x^*\|_2^2)^2.
\end{align*}
Therefore, the expected total cost $C = T^* C_\mathrm{iter}^{\texttt{A}}\le 2\hat{T}C_\mathrm{iter}^{\texttt{A}}=\cO(\eps^{-2}),$ where the last equality is because $C_\mathrm{iter}^{\texttt{A}}=\cO(1)$ and $\mathbf{V}(v^{\texttt{A}})=\cO(1)$.

\subsubsection*{Nonconvex smooth case.} 
Similar to the convex case, one can verify that the stepsizes satisfy the requirement of Lemma \ref{lm:SGD}\ref{Lemma:EC1:c} for large enough $T$.
To ensure that $\hat{x}_T^{\texttt{A}}$ is an $\eps$-stationary point, it suffices to take $T$ as the smallest integer that satisfies 
\begin{align*}
T \geq \hat{T}:=\mathbf{V}(v^{\texttt{A}})(2(F(x_1)-F(x^*))+S_F)^2\eps^{-4}.
\end{align*}
Therefore, the expected total cost $C = T^* C_\mathrm{iter}^{\texttt{A}}\le 2\hat{T}C_\mathrm{iter}^{\texttt{A}}=\cO(\eps^{-4}),$ where the last equality is because $C_\mathrm{iter}^{\texttt{A}}=\cO(1)$ and $\mathbf{V}(v^{\texttt{A}})=\cO(1)$.
\QED
\end{proof}

\section{(Expected) Total Cost of Variance Reduced Methods}
\label{app_sec:variance_reduction}
In this section, we demonstrate the proof of Theorem \ref{thm:VR}.

\begin{proof}{Proof of  Theorem~\ref{thm:VR}.}
Recall that the variance of $L$-SGD with $n_L=1$, RT-MLMC, RU-MLMC, and RR-MLMC remains $\cO(1)$ when $b > c$. Thus $v^{L\mathrm{-SGD}}(x)$, $v^{\RTMLMC}(x)$, $v^\mathrm{RU-MLMC}(x)$, and $v^\mathrm{RR-MLMC}(x)$  satisfies the variance assumption of classical variance reduction methods.

As VR RU-MLMC and VR RR-MLMC have unbiased gradient estimators of $F$, the expected total cost of VR RU-MLMC and VR RR-MLMC equals the oracle complexity of classical variance reduction methods for finding approximate stationary point of $F$, which is $\cO(\eps^{-3})$~\citep{wang2018spiderboost}, times the expected cost to query the stochastic oracle, which is $\cO(1)$. Thus, the expected total cost of VR RU-MLMC and VR RR-MLMC turns out to be $\cO(\eps^{-3})$ when $b>c$.

As for $L$-SGD and RT-MLMC, since they construct unbiased gradient estimators of $F^L$, we have
$$
\EE \|\nabla F(\hat x_T))\|_2^2\leq 2 \EE \|\nabla F^L(\hat x_T))\|_2^2+ 2 \EE \|\nabla F(\hat x_T))-F^L(\hat x_T)\|_2^2.
$$
Similar to the analysis in Appendix \ref{app_sec:total_cost}, we bound the second term via Assumption \ref{assumption:gradient_bias} with properly selected truncation levels $L = \lceil 1/a \log(4M_a/\eps^2) \rceil$. As for the first term, by a similar argument of VR RU-MLMC and VR RR-MLMC, the (expected) total cost of RT-MLMC is $\cO(\eps^{-3})$ when $b>c$ and the total cost of VR $L$-SGD is $\cO(\eps^{-3}\times \eps^{-2c/a}) = \cO(\eps^{-3-2c/a})$.

Next, we demonstrate the (expected) total cost of RT-MLMC when $b\leq c$ and V-MLMC. 

By \citet{zhang2021all}, for $\gamma\leq 1/(3S_F)$, it holds that
$$
\frac{1}{T} \sum_{t=1}^{T} \EE\|\nabla F^L(x_t)\|_2^2\leq \frac{2(F^L(x_0)-F^L(x^L))}{\gamma T}-\frac{2}{3T} \sum_{t=1}^{T} \EE\|m_t\|_2^2+\frac{2}{T}\sum_{t=1}^{T}\EE \|m_t- \nabla F^L(x_t)\|_2^2.
$$
and
$$
\frac{1}{T}\sum_{t=1}^{T}\EE \|m_t- \nabla F^L(x_t)\|_2^2\leq \frac{V}{D_1}+\frac{S_F^2 Q_E\gamma^2}{D_2} \frac{2}{3T} \sum_{t=1}^{T} \EE\|m_t\|_2^2,
$$
where $V$ is the upper bound on the variance of $v_k(x)$ for any $x$ and $k$. After $T$ iteration, the total copies of the gradient estimator $v$  used is $\cO(2TD_2 + T D_1/Q_E)$.
It turns out that when
$\gamma = \cO(1)$, $D_1 = \cO(VT)$, $Q_E=D_2 = \cO(\sqrt{T})$, it guarantees that $\EE \|\nabla F(\hat x_T)\|_2^2 \leq \cO(T^{-1}).$
Therefore, to make sure that $
\hat x_T$ is an $\eps$-stationary point, the total copies of $v$ used is $\cO(VT^{3/2})=\cO(V\eps^{-3}).$
Thus the (expected) total cost of a method $A$ is the cost to construct $v$ times the total copies of $v$ required, which is $\cO(V(v^A)C_\mathrm{iter}^\mathrm{A}\eps^{-3}).$ 

Therefore, the expected total cost of VR RT-MLMC is $\cO(\eps^{-3} \log^2(\eps^{-1}))$ for $b=c$ and $\cO(\eps^{-3-2c/a})$ for $b<c$.

For VR V-MLMC,  to guarantee that each $n_l$ is an integer, we still set $N$ as we did in Theorem \ref{thm:V-MLMC}. Plugging in $V^\mathrm{V-MLMC} C_\mathrm{iter}^\mathrm
{V-MLMC}$, the total cost of VR V-MLMC becomes
\begin{equation}
C=
\begin{cases}
\cO(\eps^{-1-2\max\{1,c/a\}}) & \text{ if } b<c;\\
\tilde \cO(\eps^{-1-2\max\{1,c/a}\}) & \text{ if } b=c;\\
\cO(\eps^{-1-2\max\{1+(c-b)/a, c/a}\}) & \text{ if } b>c.
\end{cases}\tag*{\QED}
\end{equation}
\end{proof}

\section{Proofs of Technical Results in Section~\ref{section:applications}}
\label{Appendix:applicaitons}

\subsection{Conditional Stochastic Optimization}\label{Appendix:applications:1}
Inspired by \citet{hu2020biased}, we make the following assumption on CSO.
Next, we derive the bias and variance error bounds of the oracle $\mathcal{SO}^l$ for CSO.
\begin{assumption}
\label{assumption:cso}
Suppose that $\sigma_g^2 := {\sup}_{x\in\RR^d, \xi} \EE_{\eta|\xi} \norm{ g_{\eta}(x,\xi)-\EE_{\eta|\xi} g_\eta(x,\xi)}_2^2<+\infty$; $f_\xi(\cdot)$ is $S_f$-smooth and $L_f$-Lipschitz continuous; and $g_\eta(\cdot,\xi)$ is $S_g$-smooth and $L_g$-Lipschitz continuous.
\end{assumption}

\begin{proposition}[Bias and Variance of the Oracle $\mathcal{SO}^l$ for CSO]
\label{prop:cso}
Suppose Assumption~\ref{assumption:cso} holds, then 
\begin{itemize}
    \item The functions $F$ and $F^l$ are $(S_gL_f+S_fL_g^2)$-smooth  for any $l\in\NN$. In addition,
    \begin{equation}
    \label{eq:cso_bias}
    \|\nabla F^l(x)-\nabla F(x)\|_2^2\leq L_g^2S_f^2\sigma_g^2 2^{-l}, \quad
    |F^l(x) - F(x)|\leq S_f \sigma_g^22^{-l}.  
    \end{equation}
    \item The variance of the  oracle $\mathcal{SO}^l$ satisfies 
    \begin{equation}
    \label{eq:cso_variance}
    \mathrm{Var}(h^l(x,\zeta^l))\leq L_f^2L_g^2, \quad \mathrm{Var}(H^l(x,\zeta^l))\leq L_g^2S_f^2\sigma_g^22^{-l}.
    \end{equation}
\end{itemize}
\end{proposition}
\begin{proof}{Proof of Proposition~\ref{prop:cso}.}
Note that
\citet{hu2020biased} has shown that $F$ and $F^l$ are $(S_gL_f+S_fL_g^2)$-smooth and that \eqref{eq:cso_bias} holds. We first show the variance of $h^l(x,\zeta^l)$.
\begin{multline*}
        \mathrm{Var}{(h^l(x,\zeta^l))} 
    \leq 
        \EE\|h^l(x,\zeta^l)\|_2^2 
    = 
        \EE \|\nabla \hat g_{1:{2^l}}(x,\zeta^l)^\top \nabla f_{\xi^l}(\hat g_{1:{2^l}}(x,\zeta^l))\|_2^2\\
    \leq  
        \EE [\|\nabla \hat g_{1:{2^l}}(x,\zeta^l)\|_2^2 \| \nabla f_{\xi^l}(\hat g_{1:{2^l}}(x,\zeta^l))\|_2^2]
    \leq  L_f^2 L_g^2.
\end{multline*}
Next, we show the variance of $H^l(x,\zeta^l)$. 
\begin{equation*}
\begin{split}
& 
    \mathrm{Var}(H^l(x,\zeta^l))
\leq 
    \EE\|H^l(x,\zeta^l)\|_2^2 \\
= &
    \EE \Big\|\nabla \hat g_{1:2^l}^l(x,\xi^l)^\top \nabla f_{\xi^l}(\hat g_{1:2^l}^l(x,\xi^l))
    -
    \frac{1}{2}\nabla \hat g_{1:2^{l-1}}^l(x,\xi^l)^\top \nabla f_{\xi_i^l}(\hat g_{1:2^{l-1}}^l(x,\xi^l)) \\
& - 
    \frac{1}{2}\nabla \hat g_{2^{l-1}+1:2^l}^l(x,\xi^l)^\top \nabla f_{\xi_i^l}(\hat g_{2^{l-1}+1:2^l}^l(x,\xi^l))\Big\|_2^2\\
= & 
    \EE \Big\|\frac{1}{2}\nabla \hat g_{1:2^{l-1}}^l(x,\xi^l)^\top \Big[\nabla f_{\xi^l}(\hat g_{1:2^l}^l(x,\xi^l)) - \nabla f_{\xi^l}(\hat g_{1:2^{l-1}}^l(x,\xi^l))\Big] \\
& +  
    \frac{1}{2}\nabla \hat g_{2^{l-1}+1:2^l}^l(x,\xi^l)^\top \Big[\nabla f_{\xi^l}(\hat g_{1:2^l}^l(x,\xi^l))- \nabla f_{\xi^l}(\hat g_{2^{l-1}+1:2^l}^l(x,\xi^l))\Big]\Big\|_2^2\\
\leq &
    \frac{L_g^2}{2}\EE \Big\|\nabla f_{\xi_i^l}(\hat g_{1:2^l}^l(x,\xi^l)) - \nabla f_{\xi_i^l}(\hat g_{1:2^{l-1}}^l(x,\xi^l))\Big\|_2^2 \\
& +  
    \frac{L_g^2}{2}\EE\Big\|\nabla f_{\xi^l}(\hat g_{1:2^l}^l(x,\xi^l))- \nabla f_{\xi^l}(\hat g_{2^{l-1}+1:2^l}^l(x,\xi^l))\Big\|_2^2\\
\leq & 
    \frac{L_g^2S_f^2}{2}\Big[\EE \Big\|\hat g_{1:2^l}^l(x,\xi^l) - \hat g_{1:2^{l-1}}^l(x,\xi^l)\Big\|_2^2 
    + \EE\Big\|\hat g_{1:2^l}^l(x,\xi^l)- \hat g_{2^{l-1}+1:2^l}^l(x,\xi^l)\Big\|_2^2\Big]\\
= &
    \frac{L_g^2S_f^2}{4}\EE \Big\|\hat g_{2^{l-1}+1:2^l}^l(x,\xi^l) - \hat g_{1:2^{l-1}}^l(x,\xi^l)\Big\|_2^2 
\leq  
    \frac{L_g^2S_f^2}{4}\frac{2\sigma_g^2}{2^{l-1}} = 
    \frac{L_g^2S_f^2\sigma_g^2}{2^{l}},
\end{split}
\end{equation*}
where the second equality holds as $\hat g_{1:2^l}^l(x,\xi^l) = \frac{1}{2} \Big[\hat g_{1:2^{l-1}}^l(x,\xi^l)+ \hat g_{1+2^{l-1}:2^l}^l(x,\xi^l)\Big],$
the second inequality holds by Lipschitz continuity of $g_\eta(\cdot,\xi)$, the third inequality holds by Lipschitz smoothness of $f_\xi$, and the last inequality uses the fact that $\hat g_{2^{l-1}+1:2^l}^l(x,\xi^l) $ and $\hat g_{1:2^{l-1}}^l(x,\xi^l)$ are independently identical distributed for a given $\xi^l$. The other inequalities hold by definition.
\QED
\end{proof}

\subsubsection{Special Case: Invariant Learning against Linear Transformation.}
\label{Appendix:CSO:special_case}
In this part, we demonstrate a special case of CSO problems when $g_\eta(x,\xi)$ is linear in $x$ for any realization of $\eta$ and $\xi$.   Such special cases cover the invariant learning problems with linear transformation mentioned in \citet{hu2020biased}. For such special cases, the parameters of $H^l$ are $a=1$, $b=2$, $c=1$ so that the unbiased MLMC methods, RU-MLMC and RR-MLMC, are applicable. The following assumption formally characterizes the special setting.

\begin{assumption} 
\label{assumption:cso_additional}
Suppose that the outer function $f_\xi$ is $S_f$ smooth and its gradient $\nabla f_\xi$ is $S_\nabla$ smooth;
the inner function $g_\eta(\cdot,\xi)$ is linear in $x$ for any given $\eta$ and $\xi$ and $\|\nabla g_\eta(x,\xi)\|_2\leq L_g$ for any $x, \eta,  \xi$;
and 
the inner function $g_\eta(x,\xi)$ is sub-Gaussian with proxy $\sigma_g^2$ for any given $\xi$ and $x$.
\end{assumption}

We still construct $h^l$ and $H^l$ as we did in \eqref{eq:cso_oracle_return}. Following a similar argument, it is obvious that  $F$ and $F^l$ are $L_gS_f$-smooth, $|F(x)-F^l(x)|\leq S_f\sigma_g^2 2^{-l}$, and that $\mathrm{Var}(h^l(x,\zeta^l))\leq L_f^2L_g^2$. The next proposition demonstrates the variance of $H^l$.

\begin{proposition}\label{Proposition:invariant:learning}
Under Assumption \ref{assumption:cso_additional}, 
the variance of $H^l$   satisfies
$
\mathrm{Var}(H^l(x,\zeta^l))\leq   \frac{24L_g^2 S_\nabla\sigma^4}{2^{2(l-1)}}.
$
\end{proposition}
\begin{proof}{Proof of Proposition~\ref{Proposition:invariant:learning}.}
Here,
\begin{equation*}
\begin{split}
& \mathrm{Var}(H^l(x,\zeta^l))
\leq 
\mathbb{E}\|H^l(x,\zeta^l)\|_2^2 \\
= &
\mathbb{E} \Big\|\nabla \hat g_{1:2^l}^l(x,\xi^l)^\top \nabla f_{\xi^l}(\hat g_{1:2^l}^l(x,\xi^l))
-
\frac{1}{2}\nabla \hat g_{1:2^{l-1}}^l(x,\xi^l)^\top \nabla f_{\xi_i^l}(\hat g_{1:2^{l-1}}^l(x,\xi^l)) \\
- &
\frac{1}{2}\nabla \hat g_{2^{l-1}+1:2^l}^l(x,\xi^l)^\top \nabla f_{\xi_i^l}(\hat g_{2^{l-1}+1:2^l}^l(x,\xi^l))\Big\|_2^2\\
= & 
\mathbb{E} \Big\|\frac{1}{2}\nabla \hat g_{1:2^{l-1}}^l(x,\xi^l)^\top \Big[ \nabla f_{\xi^l}(\hat g_{1:2^{l-1}}^l(x,\xi^l))-\nabla f_{\xi^l}(\hat g_{1:2^l}^l(x,\xi^l))\\
&\qquad\qquad\qquad\qquad\qquad\qquad-\nabla^2 f_{\xi^l}(\hat g_{1:2^l}^l(x,\xi^l)^\top (\hat g_{1:2^{l-1}}^l(x,\xi^l)-\hat g_{1:2^l}^l(x,\xi^l)) \Big] \\
+ & \frac{1}{2}\nabla \hat g_{2^{l-1}+1:2^l}^l(x,\xi^l)^\top \Big[\nabla f_{\xi^l}(\hat g_{2^{l-1}+1:2^l}^l(x,\xi^l))-\nabla f_{\xi^l}(\hat g_{1:2^l}^l(x,\xi^l))
\\
&\qquad\qquad\qquad\qquad\qquad\qquad-\nabla^2 f_{\xi^l}(\hat g_{1:2^l}^l(x,\xi^l)^\top (\hat g_{1+2^{l-1}:2^{l}}^l(x,\xi^l)-\hat g_{1:2^l}^l(x,\xi^l))\Big]\\
+ & 
\frac{1}{2}\nabla \hat g_{1:2^{l-1}}^l(x,\xi^l)^\top\nabla^2 f_{\xi^l}(\hat g_{1:2^l}^l(x,\xi^l)^\top (\hat g_{1:2^{l-1}}^l(x,\xi^l)-\hat g_{1:2^l}^l(x,\xi^l)) \\
+ &
\frac{1}{2}\nabla \hat g_{2^{l-1}+1:2^l}^l(x,\xi^l)^\top\nabla^2 f_{\xi^l}(\hat g_{1:2^l}^l(x,\xi^l)^\top (\hat g_{1+2^{l-1}:2^{l}}^l(x,\xi^l)-\hat g_{1:2^l}^l(x,\xi^l))
\Big\|_2^2\\
\leq &
3\mathbb{E} \Big[\Big\|\frac{1}{2}\nabla \hat g_{1:2^{l-1}}^l(x,\xi^l)\|^2  \frac{S_\nabla}{2} \|\hat g_{1:2^{l-1}}^l(x,\xi^l)-\hat g_{1:2^l}^l(x,\xi^l))\Big\|_2^4\Big]\\
+ & 3\mathbb{E} \Big[\Big\|\frac{1}{2}\nabla \hat g_{1:2^{l-1}}^l(x,\xi^l)\|^2  \frac{S_\nabla}{2} \|\hat g_{1+2^{l-1}:2^{l}}^l(x,\xi^l)-\hat g_{1:2^l}^l(x,\xi^l))\Big\|_2^4\Big]\\
+ &
3\mathbb{E} \Big\|\frac{1}{4}(\nabla \hat g_{1:2^{l-1}}^l(x,\xi^l)-\nabla \hat g_{2^{l-1}+1:2^l}^l(x,\xi^l))^\top\nabla^2 f_{\xi^l}(\hat g_{1:2^l}^l(x,\xi^l)^\top(\hat g_{1+2^{l-1}:2^{l}}^l(x,\xi^l)-\hat g_{1:2^{l-1}}^l(x,\xi^l)) \Big\|_2^2\\
\leq &
\frac{3L_g^2}{4}\frac{S_\nabla}{2}\frac{2}{4} \mathbb{E} \|\hat g_{1+2^{l-1}:2^{l}}^l(x,\xi^l)-\hat g_{1:2^{l-1}}^l(x,\xi^l))\|_2^4
\\
\leq &
\frac{3L_g^2}{4}\frac{S_\nabla}{2}\frac{2}{4} 2 (4\mathbb{E} \|\hat g_{1+2^{l-1}:2^{l}}^l(x,\xi^l)-\hat g_{1:2^{l-1}}^l(x,\xi^l))\|_2^2)^2\\
\leq &
6L_g^2 S_\nabla \frac{4\sigma^4}{2^{2(l-1)}}.
\end{split}
\end{equation*}
The first inequality uses Cauchy-Schwarz inequality and the definition of $\hat g_{1:2^{l}}^l(x,\xi^l)$. The second inequality uses linearity of $g_\eta(\cdot,\xi)$, which naturally implies that $g$ is Lipschitz continuous.
The third inequality holds as $\hat g_{1+2^{l-1}:2^{l}}^l(x,\xi^l)-\hat g_{1:2^{l-1}}^l(x,\xi^l))$ is a zero-mean sub-Gaussian random variable and the moment definition of sub-Gaussian random variable. ($X$ is a zero-mean sub-Gaussian random variable if it holds $\mathbb{E}[X^{2q}]\leq q!(4\mathrm{Var}(X))^q$ for any $q\in \NN$, see \citet{rivasplata2012subgaussian} for details). \QED
\end{proof}
The proposition implies that CSO problems with linear sub-Gaussian inner function have $b=2$ under Assumption \ref{assumption:cso_additional}. As $c=1$, it holds that $b>c$. Therefore, the unbiased MLMC methods (RU-MLMC and RR-MLMC) are applicable, and we summarize the sample complexity in the following corollary. Note that the sample complexity of $L$-SGD stays the same as in Corollary \ref{cor:cso} since the additional linearity assumption does not affect the approximation error and the cost. Such a result further motivates \citep{goda2023constructing} to study other situations when one can obtain $b>c$ for CSO problems after the preliminary version of the manuscript is released.

\begin{corollary}[Sample Complexity for Solving CSO for linear $g_\eta(x,\xi)$]
\label{cor:cso_special_case}
Under Assumption \ref{assumption:cso_additional}, the sample complexity of V-MLMC, RT-MLMC, RU-MLMC, and RR-MLMC for finding $\eps$-optimal solution is $\cO(\eps^{-1})$ for strongly convex CSO problems and $\cO(\eps^{-2})$ for convex CSO problems; the sample complexity for finding $\eps$-stationary point of these methods is $\cO(\eps^{-4})$ for nonconvex smooth CSO problems. The corresponding VR MLMC methods achieve $\cO(\eps^{-3})$ sample complexity for nonconvex smooth CSO problems.
\end{corollary}

\subsection{Joint Pricing and Staffing in Stochastic Systems}\label{Appendix:applications:2}
Inspired by \citep[Assumptions~1 and 2]{chen2023online}, we consider the following assumptions for Problem~\eqref{eq:pricing}.
\begin{assumption}\label{Assumption:pricing}
\begin{enumerate}
    \item\label{Assumption:pricing:I}
The arrival rate $\lambda(p)$ is continuously differentiable and non-increasing in $p$.
    \item\label{Assumption:pricing:II}
The staffing cost $c(\mu)$ is continuously differentiable and non-decreasing in $\mu$.
    \item\label{Assumption:pricing:III}
The domain of \eqref{eq:pricing} is $\mathcal{B}=[\underline{\mu}, \overline{\mu}]\times[\underline{p}, \overline{p}]$, where the lower bounds $\underline{\mu}$ and $\underline{p}$ satisfy that $\lambda(\underline{p}) < \underline{\mu}$.
    \item\label{Assumption:pricing:IV}
Let $\phi_V(\theta)=\log\mathbb{E}[\exp(\theta V_n)]$ and $\phi_U(\theta)=\log\mathbb{E}[\exp(\theta U_n)]$ be the cummulant generating functions of $V$ and $U$, respectively.
There exists a sufficiently small constant $\eta>0$ such that $\phi_V(\eta)<\infty, \phi_U(\eta)<\infty$. 
In addition, there exist constants $\theta\in(0, \eta/(2\overline{\mu}))$, $\mathfrak{b}\in(0,
\frac{\underline{\mu} - \lambda(\underline{p})}{\underline{\mu} +\lambda(\underline{p})})$, and $\mathfrak{c}>0$ such that $\phi_U(-\theta)<-(1-\mathfrak{b})\theta - \mathfrak{c}$ and $\phi_V(\theta) < (1+\mathfrak{b})\theta - \mathfrak{c}$.
\end{enumerate}
\end{assumption}
The proof in this section largely relies on the following technical lemma with its analysis adopted from \cite[Lemmas~2 and 6]{chen2023online}.
\begin{lemma}\label{Lemma:pricing}
Under Assumption~\ref{Assumption:pricing}, it holds that
\[
\Big| 
\mathbb{E}[W_n(\mu,p) - W_{\infty}(\mu,p)]
\Big|\le A\cdot e^{-\mathfrak{c}n},\qquad
\Big| 
\mathbb{E}[X_n(\mu,p) - X_{\infty}(\mu,p)]
\Big|\le A\cdot e^{-0.25\mathfrak{c} n},
\]
and
\[
\mathbb{E}[W_n(\mu,p) - W_{\infty}(\mu,p)]^2
\le A\cdot e^{-\mathfrak{c}n},\qquad
\mathbb{E}[X_n(\mu,p) - X_{\infty}(\mu,p)]^2
\le A\cdot e^{-0.25\mathfrak{c} n},
\]
where $A>0$ is a constant depending on $\eta,\theta,\mathfrak{b},\mathfrak{c},\overline{\mu}, \underline{\mu}, \overline{p}, \underline{p}$ only.
\end{lemma}

\begin{proposition}[Bias and Variance of the Oracle $\mathcal{SO}^l$ in \eqref{eq:pricing}]\label{Prop:pricing}
\label{prop:pricing}
Under Assumption~\ref{Assumption:pricing}, the bias and variance of the oracle $\mathcal{SO}^l$ satisfy
\begin{equation}
    \label{eq:pricing_bias}
    \|\nabla F^l(x)-\nabla F(x)\|_2^2\leq
    \mathcal{O}(e^{-0.25\mathfrak{c}\cdot 2^{l}}).
    \end{equation}
and
\begin{equation}    \label{eq:pricing_var}
    \mathrm{Var}(h^l(x,\zeta^l))=\cO(1), \quad \mathrm{Var}(H^l(x,\zeta^l))= \mathcal{O}(e^{-0.125\mathfrak{c}\cdot 2^{l}}).
    \end{equation}
    where $\mathcal{O}(\cdot)$ hides constants depending on $\eta,\theta,\mathfrak{b},\mathfrak{c},\overline{\mu}, \underline{\mu}, \overline{p}, \underline{p}$ only.
\end{proposition}
Although the bias in \eqref{eq:pricing_bias} and the variance $\mathrm{Var}(H^l(x,\zeta^l))$ decays in exponential of exponential rate, 
as discussed in \citep[Section~3.1]{chen2023online}, the exponent $\mathfrak{c}$ is a sufficiently small number such that, in practice, the decaying rate is not observed to be fast if the level $l$ is not chosen to be too large.
Inspired by this observation, we use a conservative error bound on bias and variance:
\begin{equation}
\|\nabla F^l(x)-\nabla F(x)\|_2^2\leq
    M_1\cdot 2^{-l},\quad 
    \mathrm{Var}(H^l(x,\zeta^l))\le M_1\cdot 2^{-l},
    \label{Eq:bias:variance}
\end{equation}
where $M_1>0$ is a constant independent of $l$. 
We are ready to show Corollary~\ref{result:pricing:capacity:formula}.

\proof{Proof of Corollary~\ref{result:pricing:capacity:formula}.}
For the case where $F(\cdot)$ is $S_F$-smooth and satisfies PL condition with parameter $\mu$, recall from \eqref{eq:compare_gradient_2} that if we take step size $\gamma_t=\frac{2}{\mu(t + 2S_F/\mu - 1)}$, it holds that
\[
\mathbb{E}[F(x_T) - F(x^*)]\le 
\frac{2\max\Big\{S_F\textbf{V}(v^{\texttt{A}}),
\mu^2(1+\mathfrak{a})(F(x_1) - F(x^*))
\Big\}}{\mu^2(t+\mathfrak{a})} + 
\frac{\cO(2^{-L})}{2\mu}.
\]
In order to reach $\epsilon$-optimal solution, it suffices to take $L=\cO(\log\frac{1}{\eps})$.
Moreover, we take $T={\cO}(\eps^{-1})$ for $L$-SGD and $T=\tilde\cO(\eps^{-1})$ for RT-MLMC.

For the case where $F(\cdot)$ is $S_F$-smooth and $L_F$-Lipschitz continuous, suppose we take step size $\gamma_t\in(0,1/(2S_F)]$, by \citep[Lemma~1]{hu2023contextual}, it holds that
\begin{align*}
\mathbb{E}\|\nabla F(x_T)\|^2
&\le \frac{2(F(x_1) - F(x^*))}{\gamma T} + \frac{2}{T}\sum_{t=1}^T\Big[ 
L_F\Big\|
\nabla F^L(x_t) - \nabla F(\theta_t)
\Big\|
+
S_F\gamma\mathbb{E}\|v(\theta_t) - \nabla F(x_t)\|^2
\Big]\\
&\le 
\frac{2(F(x_1) - F(x^*))}{\gamma T} +
2L_F\cdot \cO(2^{-L/2}) + 2S_F\gamma\cdot (2\cO(2^{-L}) + 2\textbf{V}(v^{\texttt{A}})),
\end{align*}
where the last inequality is because of \eqref{eq:UBSR_bias} and the relation
\[
\mathbb{E}\|v(x_t) - \nabla F(x_t)\|^2
\le 2\|\mathbb{E}v(x_t)-\nabla F(x_t)\|^2
+
2\mathbb{E}\|\mathbb{E}v(\theta_t)-v(\theta_t)\|^2=2\cO(2^{-L}) + 2\textbf{V}(v^{\texttt{A}}).
\]
In order to reach $\epsilon$-optimal solution, it suffices to take $L=\cO(\log\frac{1}{\eps^2})$.
For $L$-SGD, we take $\gamma=\cO(T^{-1/2})$ and $T=\cO(\eps^{-4})$.
For RT-MLMC, we take $\gamma=\cO(T^{-1/2})$ and $T=\tilde{\cO}(\eps^{-4})$.
\QED
\endproof

\proof{Proof of Proposition~\ref{Prop:pricing}.}
For the bias bound, one can check that 
\begin{align*}
&\left| 
\frac{\partial F}{\partial\mu}(x) - \frac{\partial F^l}{\partial\mu}(x)
\right|^2
=
\big(h_0\lambda'(p)\big)^2\Big[ 
\EE[\hat{g}_{-m:2^l}(x,\zeta^l)] - \EE[W_{\infty}(\mu,p)+X_{\infty}(\mu,p)]
\Big]^2\\
\le&
\frac{\big(h_0\lambda'(p)\big)^2}{m}\sum_{j\in[2^l-m:2^l]}\Big( 
|\mathbb{E}[W_j(\mu,p) - W_{\infty}(\mu,p)]| + |\mathbb{E}[X_j(\mu,p) - X_{\infty}(\mu,p)]|
\Big)^2\\
\le& 
\frac{2\big(h_0\lambda'(p)\big)^2}{m}\sum_{j\in[2^l-m:2^l]}\Big( 
|\mathbb{E}[W_j(\mu,p) - W_{\infty}(\mu,p)]|^2 + |\mathbb{E}[X_j(\mu,p) - X_{\infty}(\mu,p)]|^2
\Big)\\
\le &\left[ 
2\big(h_0\lambda'(p)\big)^2\cdot A\cdot e^{\mathfrak{c}m}
\right]e^{-\mathfrak{c}\cdot 2^l} + \left[ 
2\big(h_0\lambda'(p)\big)^2\cdot A\cdot e^{0.25\mathfrak{c} m}
\right]e^{-0.25\mathfrak{c}\cdot 2^l}
\le4\big(h_0\lambda'(p)\big)^2Ae^{\mathfrak{c} m}e^{-0.25\mathfrak{c}\cdot 2^l},
\end{align*}
where the third inequality is by Lemma~\ref{Lemma:pricing}.
Similarly, 
\begin{align*}
\left| 
\frac{\partial F}{\partial\mu}(x) - \frac{\partial F^l}{\partial\mu}(x)
\right|^2
&=
\big(h_0\frac{\lambda(p)}{\mu}\big)^2\Big[ 
\EE[\hat{g}_{-m:2^l}(x,\zeta^l)] - \EE[W_{\infty}(\mu,p)+X_{\infty}(\mu,p)]
\Big]^2\\&\le 4\big(h_0\frac{\lambda(p)}{\mu}\big)^2Ae^{\mathfrak{c} m} e^{-0.25\mathfrak{c}\cdot 2^l}.
\end{align*}
Therefore, it holds that $\|\nabla F^l(\mu,p)-\nabla F(\mu,p)\|_2^2\leq A'\cdot e^{-0.25\gamma\cdot 2^l}$, where the constant 
\begin{equation}
A'=4Ae^{\mathfrak{c} m}\cdot\max_{(\mu,p)\in\mathcal{B}}~\left\{
\big(h_0\lambda'(p)\big)^2+\big(h_0\frac{\lambda(p)}{\mu}\big)^2
\right\}.
\label{Eq:A'}
\end{equation}
For the variance bound, it is easy to show
\begin{align*}
\mathrm{Var}(h^l(x,\zeta^l))&\le \mathbb{E}\|h^l(x,\zeta^l)\|_2^2\le 
B_0
+
4A'\mathbb{E}[\hat{g}_{-m:2^l}(x,\zeta^l)^2]\end{align*}
where the constants
\[
\begin{aligned}
    B_0&=\max_{(\mu,p)\in\mathcal{B}}~\left\{
2\left[ 
(\lambda(p) + p\lambda'(p))^2 + (c'(\mu))^2
\right]
+
\frac{4}{\underline{\mu}}\left[ 
(h_0\lambda'(p))^2 + (\frac{h_0\lambda(p)}{\mu})^2
\right]
\right\}.
\end{aligned}
\]
Denote by $\hat{W}_j(\mu,p)$ and $\hat{X}_j(\mu,p)$ the waiting and busy time of the stochastic system with the steady-state initial state, which turns out to be \emph{time-stationary}.
Then we have
\begin{align*}
\mathbb{E}[\hat{g}_{-m:2^l}(x,\zeta^l)^2]&=\frac{1}{m^2}\mathbb{E}\left[ 
\sum_{n\in[2^l-m:2^l]}W_j(\mu,p) + X_j(\mu,p)
\right]^2
\\
&\le \frac{1}{m^2}\mathbb{E}\left[ 
\sum_{n\in[2^l-m:2^l]}\hat{W}_j(\mu,p) + \frac{\lambda(\underline{p})}{\lambda(\overline{p})}\hat{X}_j(\mu,p)
\right]^2\\
&\le \frac{2}{m^2}\mathbb{E}\left[ 
\sum_{n\in[2^l-m:2^l]}\hat{W}_j(\mu,p)
\right]^2 + \frac{2}{m^2}\mathbb{E}\left[ 
\sum_{n\in[2^l-m:2^l]}\frac{\lambda(\underline{p})}{\lambda(\overline{p})}\hat{X}_j(\mu,p)
\right]^2\\
&\le 2\mathbb{E}[\hat{W}_0(\mu,p)^2] + 2\mathbb{E}[(\frac{\lambda(\underline{p})}{\lambda(\overline{p})}\hat{X}_0(\mu,p))^2]\le 4B_1,
\end{align*}
where the first inequality is because  $\hat{W}_j(\mu,p)$ and $\frac{\lambda(\underline{p})}{\lambda(\overline{p})}\hat{X}_j(\mu,p)$ dominates $W_j(\mu,p)$ and $X_j(\mu,p)$, respectively~\cite[Lemma~1]{chen2023online}, and the constant
\begin{equation}
    B_1=\max_{(\mu,p)\in\mathcal{B}}~\left\{ 
\max\{\mathbb{E}[W_{\infty}(\mu,p)^2], \mathbb{E}[\frac{\lambda(\underline{p})^2}{\lambda(\overline{p})^2}X_{\infty}(\mu,p)^2]\}
    \right\}
\end{equation}
Finally, we find
\begin{align*}
&\mathrm{Var}(H^l(x,\zeta^l))\le \mathbb{E}\|H^l(x,\zeta^l)\|_2^2\le B_2
\mathbb{E}[\hat{g}_{-m:2^l}(x, \zeta^l) - \hat{g}_{-m:2^{l-1}}(x,\zeta^l)]^2\\
\le &2B_2\Big\{ 
\mathbb{E}[\hat{g}_{-m:2^l}(x, \zeta^l) - \EE[W_{\infty}(\mu,p)+X_{\infty}(\mu,p)]]^2 + \mathbb{E}[\hat{g}_{-m:2^{l-1}}(x, \zeta^l) - \EE[W_{\infty}(\mu,p)+X_{\infty}(\mu,p)]]^2\Big\},
\end{align*}
in which for the second inequality, we take $B_2=\max_{(\mu,p)\in\mathcal{B}}~\left\{
\big(h_0\lambda'(p)\big)^2+\big(h_0\frac{\lambda(p)}{\mu}\big)^2
\right\}.$
In detail,
\begin{align*}
&\mathbb{E}[\hat{g}_{-m:2^l}(x, \zeta^l) - \EE[W_{\infty}(\mu,p)+X_{\infty}(\mu,p)]]^2\\
\le&\frac{2}{m}\sum_{n\in[2^l-m]}\Big(\mathbb{E}[X_n(\mu,p) - X_{\infty}(\mu,p)]^2 + \mathbb{E}[W_n(\mu,p) - W_{\infty}(\mu,p)]^2\Big)\le 4Ae^{\mathfrak{c} m}\cdot e^{-0.25\mathfrak{c}\cdot 2^l}.
\end{align*}
and similarly, $\mathbb{E}[\hat{g}_{-m:2^{l-1}}(x, \zeta^l) - \EE[W_{\infty}(\mu,p)+X_{\infty}(\mu,p)]]^2\le 4Ae^{\mathfrak{c} m}\cdot e^{-0.25\mathfrak{c}\cdot 2^{l-1}}.$ Combining those results, we obtain the desired bound on $\mathrm{Var}(H^l(x,\zeta^l))$.
\QED
\endproof

\subsection{UBSR Optimization}\label{Appendix:applications:3}
We adopt the technical conditions in \cite[Assumptions~1 to 9]{hegde2021online} to consider UBSR optimization in Corollary~\ref{cor:UBSR_O}.
Lemma~1 in \cite{hegde2021online} provided the bias and the mean-square error regarding the estimator \eqref{eq:UBSR}:
\[
\mathbb{E}
\left| 
h^l(\theta, \zeta^l)
-
\nabla\mathcal{SR}_{\lambda}(\theta)
\right|=\cO(2^{-l/2}),\quad 
\mathbb{E}
\left| 
h^l(\theta, \zeta^l)
-
\nabla\mathcal{SR}_{\lambda}(\theta)
\right|^2=\cO(2^{-l}),
\]
based on which we obtain the bias and variance of the oracle $\mathcal{SO}^l$ in \eqref{eq:UBSR}.
\begin{proposition}[Bias and Variance of the Oracle $\mathcal{SO}^l$ in \eqref{eq:UBSR}]
\label{prop:UBSR}
Under the technical conditions in \cite[Assumptions~1 to 9]{hegde2021online}, the bias and variance of the oracle $\mathcal{SO}^l$ satisfy
\begin{equation}
    \label{eq:UBSR_bias}
    |\nabla\mathcal{SR}_{\lambda}^l(\theta)-\nabla\mathcal{SR}_{\lambda}(\theta)|^2=
    \mathcal{O}(2^{-l})
    \end{equation}
    and
    \begin{equation}
    \label{eq:UBSR_var}
    \mathrm{Var}(h^l(\theta,\zeta^l))= \cO(1), \quad \mathrm{Var}(H^l(\theta,\zeta^l))= \mathcal{O}(2^{-\ell}).
    \end{equation}
\end{proposition}

\proof{Proof of Proposition~\ref{prop:UBSR}}
By Jensen's inequality,
\[
|\nabla\mathcal{SR}_{\lambda}^l(\theta)-\nabla\mathcal{SR}_{\lambda}(\theta)|
=
|\mathbb{E}h^l(\theta, \zeta^l)
-\nabla\mathcal{SR}_{\lambda}(\theta)|
\le 
\mathbb{E}|h^l(\theta, \zeta^l)
-\nabla\mathcal{SR}_{\lambda}(\theta)|=\cO(2^{-l/2}),
\]
which implies the relation \eqref{eq:pricing_bias}.
Next, one can see 
\[
\mathrm{Var}(h^l(\theta,\zeta^l)) \le
\mathbb{E}|h^l(\theta,\zeta^l)|^2\le 
2\mathbb{E}
\left| 
h^l(\theta, \zeta^l)
-
\nabla\mathcal{SR}_{\lambda}(\theta)
\right|^2
+
2
\left|
\nabla\mathcal{SR}_{\lambda}(\theta)
\right|^2 = \cO(1),
\]
and
\begin{align*}
\mathrm{Var}(H^l(x,\zeta^l))&\le \mathbb{E}|H^l(x,\zeta^l)|^2
=
 \mathbb{E}\Big|
\frac{1}{2}(\hat{g}_{1:2^l}(\theta, \xi^l) - \hat{g}_{1:2^{l-1}}(\theta, \xi^l))
+
\frac{1}{2}(\hat{g}_{1:2^l}(\theta, \xi^l) - \hat{g}_{2^{l-1}+1:2^l}(\theta, \xi^l))
 \Big|^2\\
 &\le 
 \frac{1}{2}\mathbb{E}\Big|\hat{g}_{1:2^l}(\theta, \xi^l) - \hat{g}_{1:2^{l-1}}(\theta, \xi^l)\Big|^2
 +
 \frac{1}{2}\mathbb{E}\Big|\hat{g}_{1:2^l}(\theta, \xi^l) - \hat{g}_{2^{l-1}+1:2^l}(\theta, \xi^l)\Big|^2.
\end{align*}
Specifically, we find
\begin{align*}
&\mathbb{E}\Big|\hat{g}_{1:2^l}(\theta, \xi^l) - \hat{g}_{1:2^{l-1}}(\theta, \xi^l)\Big|^2
\le 
2\mathbb{E}\Big|\hat{g}_{1:2^l}(\theta, \xi^l) - \nabla\mathcal{SR}_{\lambda}(\theta)\Big|^2
+
2\mathbb{E}\Big|\hat{g}_{1:2^{l-1}}(\theta, \xi^l) - \nabla\mathcal{SR}_{\lambda}(\theta)\Big|^2\\
\le&
2\mathbb{E}\Big|h^l(\theta, \zeta^l) - \nabla\mathcal{SR}_{\lambda}(\theta)\Big|^2
+
2\mathbb{E}\Big|h^{l-1}(\theta, \zeta^l)- \nabla\mathcal{SR}_{\lambda}(\theta)\Big|^2
=2\cO(2^{-l}) + 2\cO(2^{-l-1})=\cO(2^{-l}),
\end{align*}
and the term $\mathbb{E}\Big|\hat{g}_{1:2^l}(\theta, \xi^l) - \hat{g}_{2^{l-1}+1:2^l}(\theta, \xi^l)\Big|^2$ can be bounded using the similar manner.
Thus, we conclude the relation \eqref{eq:UBSR_var} also holds.
\QED
\endproof

\begin{corollary}[Full Version of Corollary~\ref{cor:UBSR_O}.]
Let $\mathtt{A}\in\{\text{V-MLMC}, \text{RT-MLMC}\}$ and $\mathfrak{a}=2S_F/\mu-1$.
    Under technical conditions in \cite[Assumptions~1 to 9]{hegde2021online}, with properly chosen hyper-parameters, it holds that 
    \begin{enumerate}
        \item 
    Suppose $\mathcal{SR}_{\lambda}(\theta)$ is $S_F$-smooth and satisfies PL condition with parameter $\mu$, then the total cost of $\mathtt{A}$ for finding $\eps$-optimal solution is $\tilde{\cO}(\eps^{-1})$.
       \item
    Suppose $\mathcal{SR}_{\lambda}(\theta)$ is $S_F$-smooth, then the total cost of $\mathtt{A}$ for finding $\eps$-stationary point is $\tilde{\cO}(\eps^{-4})$.
    \end{enumerate}
\end{corollary}

\proof{Proof of Corollary~\ref{cor:UBSR_O}.}
For the case where $\mathcal{SR}_{\lambda}(\theta)$ is $S_F$-smooth and satisfies PL condition with parameter $\mu$, recall from \eqref{eq:compare_gradient_2} that if we take step size $\gamma_t=\frac{2}{\mu(t + 2S_F/\mu - 1)}$, it holds that
\[
\mathbb{E}[\mathcal{SR}_{\lambda}(\theta_T) - \mathcal{SR}_{\lambda}(\theta^*)]\le 
\frac{2\max\Big\{S_F\textbf{V}(v^{\texttt{A}}),
\mu^2(1+\mathfrak{a})(\mathcal{SR}_{\lambda}(\theta_1) - \mathcal{SR}_{\lambda}(\theta^*))
\Big\}}{\mu^2(t+\mathfrak{a})} + 
\frac{\cO(2^{-L})}{2\mu}.
\]
In order to reach $\epsilon$-optimal solution, it suffices to take $L=\cO(\log\frac{1}{\eps})$ and $T=\tilde{\cO}(\eps^{-1})$.

For the case where $\mathcal{SR}_{\lambda}(\theta)$ is $S_F$-smooth and $L_F$-Lipschitz continuous, suppose we take step size $\gamma_t\in(0,1/(2S_F)]$, by \citep[Lemma~1]{hu2023contextual}, it holds that
\begin{align*}
\mathbb{E}|\nabla \mathcal{SR}_{\lambda}(\theta_T)|^2
&\le \frac{2(\mathcal{SR}_{\lambda}(\theta_1) - \mathcal{SR}_{\lambda}(\theta^*))}{\gamma T}\\
&\qquad\qquad+ \frac{2}{T}\sum_{t=1}^T\Big[ 
L\Big| 
\nabla \mathcal{SR}_{\lambda}^L(\theta_t) - \nabla \mathcal{SR}_{\lambda}(\theta_t)
\Big|
+
S\gamma\mathbb{E}|v(\theta_t) - \nabla \mathcal{SR}_{\lambda}(\theta_t)|^2
\Big]\\
&\le 
\frac{2(\mathcal{SR}_{\lambda}(\theta_1) - \mathcal{SR}_{\lambda}(\theta^*))}{\gamma T} +
2L_F\cdot \cO(2^{-L/2}) + 2S_F\gamma\cdot (2\cO(2^{-L}) + 2\textbf{V}(v^{\texttt{A}})),
\end{align*}
where the last inequality is because of \eqref{eq:UBSR_bias} and the relation
\[
\mathbb{E}|v(\theta_t) - \nabla \mathcal{SR}_{\lambda}(\theta_t)|^2
\le 2|\mathbb{E}v(\theta_t)-\nabla\mathcal{SR}_{\lambda}(\theta)|^2
+
2|\mathbb{E}v(\theta_t)-v(\theta_t)|^2=2\cO(2^{-L}) + 2\textbf{V}(v^{\texttt{A}}).
\]
In order to reach $\epsilon$-optimal solution, it suffices to take $L=\cO(\log\frac{1}{\eps^2})$.
For RT-MLMC, we take $\gamma=\cO(T^{-1/2})$ and $T=\tilde{\cO}(\eps^{-4})$.
For V-MLMC, we take $\gamma=1/(2S_F)$ and $T=\tilde{\cO}(\eps^{-2})$.
\QED
\endproof

\subsection{Proof of Technical Results in Section~\ref{Sec:statistical:inference}}

\begin{assumption}\label{Assume:CSO:special}
The decision set $\cX\subseteq \mathbb{R}^d$ has a finite diameter; %
the functions $f_{\xi}(\cdot)$ and $g_{\eta}(\cdot, \xi)$ are Lipschitz continuous; 
and $\sigma^2:=\max_{x\in\cX}~\mathbb{E}_{\eta\mid\xi}\|g_{\eta}(x,\xi) - \mathbb{E}_{\eta\mid\xi}g_{\eta}(x,\xi)\|^2_2<\infty$.
\end{assumption}
\begin{proof}{Proof of Theorem~\ref{Theorem:uniform:converge}.}
To make the statement precise, denote by $D_{\cX}$ the diameter of $\cX$, $L_f$ the Lipschitz continuous constant of $f_{\xi}(\cdot)$, and $L_g$ the Lipschitz continuous constant of $g_{\xi}(\cdot,\xi)$.
Take an $\upsilon$-net on $\cX$, which is denoted as $\{x_k\}_{k=1}^J$, to ensure
\begin{equation}
\begin{split}
& 
\PP\Big(\sup_{x\in \cX} |\hat F(x)-F(x)|\geq \eps\Big)   
\leq \sum_{k=1}^J \PP\Big(|\hat F(x_k)-F(x_k)|\geq \frac{\eps}{2}\Big).    
\end{split}
\end{equation}
We follow the similar argument as in \citep[Theorem~4.1]{hu2020sample}, and the fact that $\hat{F}$ is $2(L+1)L_fL_g$ Lipschitz continuous. In this case, 
\[
J\leq \cO(1)\cdot \left(\frac{8D_\cX L_f L_g(L+1)}{\eps}\right)^d.
\]
It suffices to provide an upper bound on $J\cdot\PP(|\hat F(x)-F(x)|\geq \frac{\eps}{2})$ for any $x\in\cX$.
Define the random variable $Z_{li}(x)=V^l(x, \zeta_i^l)$, where $V^l(x, \zeta_i^l)$ is defined in \eqref{Eq:MLMC:estimator}.
Assume one can ensure
\begin{equation}\label{Eq:Vl:z}
\left|\EE\sum_{l=0}^L\frac{1}{n_l}\sum_{i=1}^{n_l}Z_{li}(x) -F(x)\right|\leq \frac{\eps}{4},
\end{equation}
then
\begin{equation}
\begin{split}
& 
J\cdot\PP(|\hat F(x)-F(x)|\geq \frac{\eps}{2}) 
= 
J\cdot\PP\left(\left|\sum_{l=0}^L \frac{1}{n_l}\sum_{i=1}^{n_l}Z_{li}(x)-F(x)\right|\geq \frac{\eps}{2}\right) \\
= & 
J\cdot\PP\left(\left|\sum_{l=0}^L \frac{1}{n_l}\sum_{i=1}^{n_l}Z_{li}(x)- \EE \sum_{l=0}^L \frac{1}{n_l}\sum_{i=1}^{n_l}Z_{li}(x)+ \EE \sum_{l=0}^L \frac{1}{n_l}\sum_{i=1}^{n_l}Z_{li}(x) -F(x)\right|\geq \frac{\eps}{2}\right) \\
\leq & 
J\cdot\PP\left(\left|\sum_{l=0}^L \frac{1}{n_l}\sum_{i=1}^{n_l}Z_{li}(x)- \EE \sum_{l=0}^L \frac{1}{n_l}\sum_{i=1}^{n_l}Z_{li}(x)\right|\geq \frac{\eps}{4}\right) \\
\leq &
J\cdot\PP\left(\sum_{l=0}^L\left| \frac{1}{n_l}\sum_{i=1}^{n_l}Z_{li}(x)- \EE Z_{li}(x)\right|\geq \frac{\eps}{4}\right) \\
\leq & 
J\cdot\sum_{l=0}^L\PP\left(\left| \frac{1}{n_l}\sum_{i=1}^{n_l}Z_{li}(x)- \EE Z_{li}(x)\right|\geq \frac{\eps}{4(L+1)}\right) \\
\leq &2J\cdot\sum_{l=0}^L  \exp\Big(-\frac{n_l \eps^2}{16(2+\bar{\epsilon}) 
\bV_l (L+1)^2}\Big),
\end{split}
\end{equation}
where the last inequality uses Cramer's large deviation Theorem, and the fact that $Z_{li}(x)$ is sub-Gaussian and $\bV_l$ is the variance of $Z_{li}(x)$.
Consequently, we take $n_l$ to ensure $\exp\Big(-\frac{n_l \eps^2}{16(2+\bar{\epsilon}) 
\bV_l (L+1)^2}\Big)\le \frac{\alpha}{2J(L+1)}$, i.e., 
$$
n_l = \left\lceil
\frac{16(2+\bar{\epsilon}) 
\bV_l (L+1)^2}{\eps^2}\Big(d\log\Big(\frac{8D_\cX L_f L_g(L+1)}{\eps}\Big)+\log\Big(\frac{2(L+1)}{\alpha}\Big) + \cO(1)\Big)
\right\rceil.
$$
To finish the proof, it remains to
\begin{enumerate}
    \item Obtain an upper bound on the $\bV_l$:
    \begin{equation}
\begin{split}
\bV_l 
\le &
\EE|Z_{li}|_2^2 
\leq 
\EE\Big|\frac{1}{2}f_{\xi_i^l}\Big(
\hat{g}_{1:m_l}(x, \zeta_i^l)
\Big)
-
\frac{1}{2}f_{\xi_i^l}\Big(\hat{g}_{1:m_{l-1}}(x, \zeta_i^l)\Big) \\
+ &
\frac{1}{2}f_{\xi_i^l}\Big(\hat{g}_{1:m_l}(x, \zeta_i^l)\Big)
-
\frac{1}{2}f_{\xi_i^l}\Big(\hat{g}_{m_{l-1}+1:m_l}(x, \zeta_i^l)\Big) 
\Big|^2\\
\leq &
\frac{L_f^2}{8}\EE \Big|\frac{1}{m_{l-1}}\sum_{j=1}^{m_{l-1}} g_{\eta_{ij}^l}(x,\xi_i^l)
-
\frac{1}{m_{l-1}}\sum_{j=1+m_{l-1}}^{m_{l}} g_{\eta_{ij}^l}(x,\xi_i^l)\Big|_2^2 
= 
\frac{L_f^2\sigma^2}{4m_{l-1}}.
\end{split}
\end{equation}
    \item Ensure \eqref{Eq:Vl:z}: it suffices to make  
    $m_L\ge\frac{16L_f^2\sigma^2}{\eps^2}$ since
$$
|\EE\sum_{l=0}^L \frac{1}{n_l}\sum_{i=1}^{n_l} Z_{li}(x) -F(x)|\leq \frac{L_f\sigma}{\sqrt{m_L}}\leq  \frac{\eps}{4}.
$$
    \item Select appropriate $L$, $n_l$, $m_l$ and compute total cost: 
    we take $m_l = 2^l, l=0,\ldots,L$, with $L=\lceil\log(16L_f^2\sigma^2/\eps^2)\rceil$, then
    \begin{equation*}
\begin{split}
 &\sum_{l=0}^L n_lm_l 
\\
=&  
\sum_{l=0}^L\left\lceil
\frac{32(2+\bar{\epsilon}) 
 (L+1)^2L_f^2\sigma^2}{\eps^2}\Big(d\log\Big(\frac{8D_\cX L_f L_g(L+1)}{\eps}\Big)+\log\Big(\frac{2(L+1)}{\alpha}\Big) + \cO(1)\Big)
\right\rceil \\
\leq &
\frac{32(2+\bar{\epsilon})(\log(16L_f^2\sigma^2/\eps^2)+2)^3L_f^2 \sigma^2}{\eps^2}\Big(d\log\Big(\frac{8D_\cX L_f L_g(\log(16L_f^2\sigma^2/\eps^2)+2)}{\eps}\Big)\\
&\qquad\qquad+\log\Big(\frac{2(\log(16L_f^2\sigma^2/\eps^2)+2)}{\alpha}\Big) + \cO(1)\Big)\\
= & \tilde{\cO}(d\eps^{-2}).
\end{split}
\end{equation*}
\end{enumerate}

It completes the proof.
\QED
\end{proof}

\section{Experiment Implementation Details}

\subsection{Synthetic Problem with Biased Oracles}
We specify the optimal choice of hyperparameters for various gradient methods using grid search. Specifically, for three biased gradient methods ($L$-SGD, V-MLMC, RT-MLMC), we consider the level $L\in \{0,\ldots,10\}$.
We stop the optimization process if the number of generated samples exceeds the total budget $\texttt{4e+4}$.
The mini-batch size for V-MLMC is chosen as $n_l=2^{L-l}$ to avoid additional costs from rounding to integer numbers.
The stepsize $\gamma$ is selected from the set $\{\texttt{1e-1}, \texttt{1e-2}, \texttt{1e-3}, \texttt{5e-4}, \texttt{1e-4}\}$.
We evaluate the performance of the various gradient methods by comparing the ground truth objective value of \eqref{Eq:SDRO} to the number of generated samples and produce an error bar based on $10$ independent trials. 
\subsection{Joint Pricing and Staffing}

We specify $\lambda(p)=\chi\cdot e^{a-p}/(1 + e^{a-p})$ with $a=0.1, \chi=10$, $c(\mu)=c_0\mu^2$ with $c_0=0.1$, and $h_0=1$;
The initial guess is specified as $\mu=9, p=9$.
One can apply the Pollaczek–Khinchine formula~\citep{haigh2002probability} to obtain the closed-form reformulation of the objective function:
\[
F(\mu,p)= h_0\Big[
(\lambda(p)/\mu) + \frac{(\lambda(p)/\mu)^2(1 + \text{Var}(U_n))}{2(1-\lambda(p)/\mu)}
\Big] + c(\mu) - p\lambda(p).
\]
We model (normalized) inter-arrival time $U_n$ to follow the exponential distribution, and the (normalized) service time $V_n$ to follow the Erlang distribution (with variance $0.1$), exponential distribution, or hyper-exponential distribution (i.e., a mixture of $10$ equal-weight exponential distributions with rate parameter $i^2\cdot 0.155$).

\subsection{Contrastive Learning}
\label{Appendix:con:learn}

The widely adopted \emph{linear evaluation protocol} is employed to examine the performance of obtained representation: we train a linear classifier based on neural network outputs and use the classification accuracy on testing data as a proxy for representation quality.
For the data augmentation step, we use a similar setup as in \cite{chen2020simple} that uses random flip (i.e., random horizontal/left-to-right flip with $0.5$ probability), color jittering, and color dropping.
The chosen neural network architecture is ResNet-18~\citep{he2016deep} followed by a two-layer ReLU neural network. 

\end{APPENDICES}
\end{APPENDICES}

\end{document}